\documentclass[10pt]{article}

\usepackage[left=1.5in, right=1.5in, bottom=1in, top=1in, footskip=0.3in, 
]{geometry}

\usepackage{xcolor}
\definecolor{LinkColor}{rgb}{0,0,1}
\definecolor{LinkColor2}{rgb}{1,0,0}
\definecolor{lbcolor}{rgb}{0.85,0.85,0.85}
\definecolor{FrameColor}{rgb}{0.85,0.85,0.85}

\usepackage[ngerman, english]{babel}
\usepackage[utf8]{inputenc}
\usepackage{enumerate}
\usepackage[normalem]{ulem}
\usepackage{setspace}

\def\pskip{\\[-3mm]}

\usepackage{multicol}
\usepackage[babel,french=guillemets,german=swiss]{csquotes}

\usepackage{array}
\usepackage{booktabs}
\usepackage{framed}
\usepackage{rotating}
\usepackage{longtable}
\usepackage{multirow}
\usepackage{tabularx}
\newcolumntype{L}[1]{>{\raggedright\arraybackslash}p{#1}} 
\newcolumntype{C}[1]{>{\centering\arraybackslash}p{#1}} 
\newcolumntype{R}[1]{>{\raggedleft\arraybackslash}p{#1}} 

\usepackage{graphicx,import}
\usepackage{caption}

\usepackage{amssymb}
\usepackage{dsfont}
\usepackage{nicefrac}
\usepackage{empheq}
\allowdisplaybreaks

\usepackage[%
pdftitle={Titel},%
pdfauthor={Autor},%
pdfcreator={LaTeX, LaTeX with hyperref and KOMA-Script},
pdfsubject={Betreff}, 
pdfkeywords={Keywords}
]{hyperref} 

\hypersetup{%
	colorlinks=true,
	linkcolor=LinkColor,%
	anchorcolor=LinkColor,%
	citecolor=LinkColor2,%
	filecolor=LinkColor,%
	menucolor=LinkColor,%
	urlcolor=LinkColor,%
}

\usepackage[savemem]{listings}
\usepackage{paralist}
{\begin{list}{$\diamondsuit$}{}}%
	{\end{list}}

\usepackage{amsthm}
\usepackage{upgreek}

\newtheoremstyle{tstyle}
{15pt}	
{5pt}	
{\itshape}	
{}	
{\bfseries}	
{.}	
{0.5em}	
{}	

\theoremstyle{tstyle}

\newtheorem{theorem}{Theorem}[section]
\newtheorem{lemma}[theorem]{Lemma}

\newtheorem{remark}[theorem]{Remark}
\newtheorem{assumptions}[theorem]{Assumptions}

\newtheoremstyle{cstyle}
{15pt}	
{5pt}	
{}	
{}	
{\bfseries}	
{}	
{0.2222em}	
{}	

\theoremstyle{cstyle}

\makeatletter
\g@addto@macro{\thm@space@setup}{\thm@headpunct{}}
\renewenvironment{proof}[1][\proofname.]{\par
	\pushQED{\qed}%
	\normalfont \topsep6\p@\@plus6\p@\relax
	\trivlist
	\item[\hskip\labelsep
	\bfseries
	#1\@addpunct{\,}]\ignorespaces
}{%
	\popQED\endtrivlist\@endpefalse
}
\makeatother

\makeatletter
\g@addto@macro{\thm@space@setup}{\thm@headpunct{}}
\newenvironment{sketch-proof}[1][Sketch of the proof]{\par
	\pushQED{\qed}%
	\normalfont \topsep6\p@\@plus6\p@\relax
	\trivlist
	\item[\hskip\labelsep
	\bfseries
	#1\@addpunct{\,}]\ignorespaces
}{%
	\popQED\endtrivlist\@endpefalse
}
\makeatother

\newcounter{subeq}

\usepackage{mathtools}  
\usepackage{amsmath}
\usepackage{fancyhdr}
\usepackage{tikz}
\usepackage{pdflscape}
\usepackage{subfig}
\usepackage{float}
\usepackage{hyperref}
\usepackage{tikz}
\usepackage{pgfplots}
\pgfplotsset{compat=1.16} 
\usepackage{algorithmic}

\usepackage[section]{placeins}         

\numberwithin{equation}{section}
\usepackage[only,llbracket,rrbracket]{stmaryrd}

\makeatletter
\newcommand{\mylabel}[2]{%
   \protected@write \@auxout {}{\string \newlabel {#1}{{#2}{\thepage}{#2}{#1}{}}}%
   \hypertarget{#1}{#2}%
}
\newcommand{\mylabelHIDE}[2]{%
   \protected@write \@auxout {}{\string \newlabel {#1}{{#2}{\thepage}{#2}{#1}{}}}%
   \hypertarget{#1}{}%
}
\makeatother

\renewcommand{\phi}{\varphi}

\newcommand{\bI}{\ensuremath{\mathbf{I}}}

\newcommand{\bn}{\ensuremath{\mathbf{n}}}
\newcommand{\bu}{\ensuremath{\mathbf{u}}}

\newcommand{\bw}{\ensuremath{\mathbf{w}}}
\newcommand{\bX}{\ensuremath{\mathbf{X}}}
\newcommand{\bz}{\ensuremath{\mathbf{z}}}

\newcommand{\bbB}{\ensuremath{\mathbb{B}}}
\newcommand{\bbC}{\ensuremath{\mathbb{C}}}
\newcommand{\bbD}{\ensuremath{\mathbb{D}}}

\newcommand{\bbG}{\ensuremath{\mathbb{G}}}
\newcommand{\bbI}{\ensuremath{\mathbb{I}}}
\newcommand{\bbN}{\ensuremath{\mathbb{N}}}
\newcommand{\bbO}{\ensuremath{\mathbb{O}}}
\newcommand{\bbP}{\ensuremath{\mathbb{P}}}
\newcommand{\bbR}{\ensuremath{\mathbb{R}}}
\newcommand{\bbS}{\ensuremath{\mathbb{S}}}
\newcommand{\bbT}{\ensuremath{\mathbb{T}}}
\newcommand{\bbU}{\ensuremath{\mathbb{U}}}

\newcommand{\D}{\ensuremath{\mathbb{D}}}

\newcommand{\I}{\ensuremath{\mathbb{I}}}  
    
\newcommand{\mycal}[1]{\ensuremath{\mathcal{#1}}}
\newcommand{\calA}{\mycal{A}}
\newcommand{\calB}{\mycal{B}}

\newcommand{\calG}{\ensuremath{\mathcal{G}}}
\newcommand{\calH}{\ensuremath{\mathcal{H}}}
\newcommand{\calI}{\mycal{I}}

\newcommand{\calK}{\mycal{K}}

\newcommand{\calP}{\ensuremath{\mathcal{P}}}
\newcommand{\calS}{\ensuremath{\mathcal{S}}}
\newcommand{\calT}{\ensuremath{\mathcal{T}}}

\newcommand{\calV}{\ensuremath{\mathcal{V}}}

\DeclareMathOperator{\Div}{div}
\DeclareMathOperator{\diam}{diam}

\DeclareMathOperator{\trace}{Tr}

\newcommand{\abs}[1]{\left\lvert{#1}\right\rvert}
\newcommand{\norm}[1]{\|{#1}\|}
\newcommand{\nnorm}[1]{\left\|{#1}\right\|}
\newcommand{\dualp}[2]{\left<{#1\, ,\, #2}\right>}
\newcommand{\skp}[2]{\left({#1\, ,\, #2}\right)}

\newcommand{\skpp}[2]{\left\llangle {#1\, , \, #2}\right\rrangle }  

\newcommand{\dv}[1]{\,{\mathrm d}#1}
\newcommand{\dx}{\dv{x}}

\newcommand{\ds}{\dv{s}}

\newcommand{\ddv}[2]{ \frac{{\mathrm d}#1}{{\mathrm d}#2} }

\newcommand{\ddt}{\ddv{}{t}}

\def\softd{{\leavevmode\setbox1=\hbox{d}%
\hbox to 1.05\wd1{d\kern-0.4ex{\char039}\hss}}}

\usepackage{MnSymbol}   
\usepackage{todonotes}

\newcommand{\jump}[1]{\llbracket{#1}\rrbracket}
\newcommand{\jumpp}[1]{\left\llbracket{#1}\right\rrbracket}
\newcommand{\cblue}[1]{{{#1}}}
\newcommand{\cred}[1]{{{#1}}}

\newcommand{\cmagenta}[1]{{{#1}}}
\newcommand{\cteal}[1]{{{#1}}}

\begin{document}

%
%
	
\begin{center}	
	\LARGE{Parametric finite element approximation of two-phase Navier--Stokes flow with viscoelasticity} 
\end{center}
\bigskip
\begin{center}	
	\normalsize{Harald Garcke}\\[1mm]
	\textit{University of Regensburg, Germany}\\[1mm]
	\texttt{Harald.Garcke@ur.de}
\end{center}
\medskip
\begin{center}	
	\normalsize{Robert N\"urnberg}\\[1mm]
	\textit{University of Trento, Italy}\\[1mm]
	\texttt{Robert.Nurnberg@unitn.it}
\end{center}
\medskip
\begin{center}	
	\normalsize{Dennis Trautwein}\\[1mm]
	\textit{University Regensburg, Germany}\\[1mm]
	\texttt{Dennis.Trautwein@ur.de}
\end{center}
\bigskip

\begin{abstract}
\footnotesize
In this work, we present a parametric finite element approximation of two-phase Navier--Stokes flow with viscoelasticity. The free boundary problem is given by the viscoelastic Navier--Stokes equations in the two fluid phases, connected by jump conditions across the interface. The elasticity in the fluids is characterised using the Oldroyd-B model with possible stress diffusion. The model was originally introduced to approximate fluid-structure interaction problems between an incompressible Newtonian fluid and a hyperelastic neo-Hookean solid, which are possible limit cases of the model.

We approximate a variational formulation of the model with an unfitted finite element method that uses piecewise linear parametric finite elements. The two-phase Navier--Stokes--Oldroyd-B system in the bulk regions is discretised in a way that guarantees unconditional solvability and stability for the coupled bulk-interface system. Good volume conservation properties for the two phases are observed in the case where the pressure approximation space is enriched with the help of an XFEM function. We show the applicability of our method with some numerical results.
\pskip

\pskip
	
\end{abstract}



%
%

%
%
\section{Introduction}

In the field of fluid mechanics, the modelling and simulation of viscoelastic materials that consist of two or more components have gained increasing importance in recent years. This area of study presents challenges arising not only from the inherent complexities of viscoelasticity, but also from the need to accurately describe the individual phases and the transitions between them.
From the modelling point of view, viscoelastic fluids require constitutive models which account for both elastic and viscous behaviour, such as the Oldroyd-B model, the FENE-P model or the Giesekus model \cite{alves_2021_viscoelastic_review}. The viscous behaviour refers to their resistance to flow, which is typical for liquids, while the elastic behaviour is related to that of solids, allowing them to return to their original shape after deformation. 
Many examples described in the literature are in the field of polymeric fluids \cite{barrett_sueli_2007, lukacova_2015, masmoudi_2018_handbook, sethian_2007_viscoelastic_jetting}, or have applications in biomedicine \cite{ambrosi_2009, GKT_2022_viscoelastic, lowengrub_2021_viscoelastic}. 

Compared to a Newtonian fluid, viscoelastic fluids have an additional stress tensor that accounts for the elastic effects. This extra stress tensor is usually described with the help of an additional tensor $\bbB$. In the context of polymeric fluids, 
the tensor $\bbB$ is often referred to as the \textit{structure or conformation tensor}. This tensor describes the molecular configurations of the polymers, capturing how the polymer chains stretch and orient under flow \cite{masmoudi_2018_handbook}. 
The polymer molecules are often modelled as chains of beads and springs, where spring forces are responsible for the movement of the molecules \cite{lukacova_2015}.
Assuming linear spring forces, using Hooke's law, leads to the Oldroyd-B model, while nonlinear descriptions of the elastic spring laws lead to nonlinear variants such as the FENE-P model.
Other models like the Giesekus model are based on different relaxation behaviours of the material.
Here, the relaxation refers to the process by which the polymer chains in the viscoelastic fluid return to their equilibrium state after being deformed by an applied stress \cite{bird_curtiss_armstrong_hassager_1987}. 
In this paper, we are going to consider the Oldroyd-B model to describe the viscoelastic effects.

The numerical treatment of viscoelastic models typically requires further attention, especially at high Weissenberg numbers, which correspond to large relaxation parameters. At high Weissenberg numbers, numerical schemes may suffer from accuracy and convergence issues in benchmark problems, which is often referred to as the \textit{high Weissenberg number problem}. Some of the reasons have been identified as the lack of existence of solutions to the mathematical system of equations or limitations of numerical schemes. These issues are discussed in, e.g., \cite{ashby_pryer_2024_oldroyd_lie_derivative, boyaval_LM_2009_oldroyd, fattal_kupferman_2005_log_formulation} and in the review article \cite{alves_2021_viscoelastic_review}. A promising strategy for addressing instabilities at high Weissenberg numbers is to understand the dissipation laws of the free energy at the continuous level, such as those for the Oldroyd-B or the FENE-P model \cite{Hu_Lelievre_2007}, and then construct numerical schemes that incorporate similar dissipation laws at the discrete level. 
Some examples for Oldroyd-B fluids have been studied in \cite{barrett_boyaval_2009, boyaval_LM_2009_oldroyd, lukacova_NS_2016_oldroyd}, while other viscoelastic models are considered in \cite{barrett_2018_fene-p, GT_2023_DCDS}. In many cases, energy estimates form the foundation for a rigorous convergence analysis of numerical schemes using compactness arguments \cite{barrett_boyaval_2009, barrett_2018_fene-p, GKT_2022_viscoelastic, GT_2023_DCDS, sieber_2020}.

In systems with two or more phases, accurately capturing the interface dynamics and interactions between viscoelastic fluids and other phases adds another layer of complexity.
The strategies for describing the interfaces can generally be divided into methods with an implicit or explicit description of the interface. The implicit methods usually use an indicator function to capture the interface. A representative is the level set method, in which a level set function is used to describe the interface \cite{osher_sethian_1988_levelset, sussman_osher_1994_levelset}. Here, some examples of models with viscoelasticity are \cite{pillapakkam_2007_viscoelastic, pillapakkam_2001_viscoelastic, sethian_2007_viscoelastic_jetting}.
Another implicit interface method is the phase-field method, which is also called the diffuse interface method in this context, as the sharp interface is approximated with the help of a smooth order parameter with a thin transition layer between the phases. This method has been used in, e.g., \cite{GKT_2022_viscoelastic, mokbel_abels_aland_2018, schmeller_2023_gradient_flows, schmeller_2023_sharp_interface, yue_feng_liu_shen_2004, yue_feng_liu_shen_2005}.
Other methods that are applied for viscoelastic two-phase models are based on the volume of fluid method \cite{bonito_2006_viscoelastic, renardy_2005_viscoelastic, stewart_sussman_2008_viscoelastic}, which goes back to the work \cite{hirt_nichols_1981_VOF} and uses a
characteristic function of one of the phases to evolve the interface, or the Lattice Boltzmann method 
used in \cite{phanthien_2019_viscoelastic, phanthien_2020_rising_bubble}.

The methods using an explicit interface description can be categorised into fitted methods, where the surface mesh aligns with the phases of the bulk mesh, and unfitted methods, where the lower dimensional surface mesh is entirely independent of the bulk mesh, as discussed in \cite{agnese_nürnberg_2020_ALE, baensch_2001, BGN_2013_stokes, BGN_2015_navierstokes, ganesan_2008_moving_mesh, ganesan_2007, GNZ_2023, hughes_1981_ALE, perot_2003_moving_mesh, unverdi_tryggvason_1992}.
\cblue{The idea of \textit{parametric finite element methods} originates from the seminal work of Dziuk \cite{dziuk_1991}. These methods track the position of the interface using a moving, lower-dimensional mesh, which evolves \textit{parametrically} at each time step through a finite element based motion. They are highly efficient, as they require fewer degrees of freedom than implicit methods like the phase-field or level set approaches, which approximate the interface as a level set of a function defined on the entire domain.}
On the other hand, special care must be taken to ensure an accurate representation of complex interface geometries and dynamics, especially when dealing with large deformations.
In particular, it is essential to ensure a good approximation of the interface mesh and to maintain a good mesh quality during the evolution. In parametric methods, the interface mesh is transported with the help of the fluid velocity \cite{baensch_2001, ganesan_2007, unverdi_tryggvason_1992}. In many cases, this can lead to mesh degeneracies, such as the coalescence of mesh points or the formation of very small angles in the polyhedral interface mesh. Consequently, frequent re-parametrisation of the interface mesh is often necessary in practice. These challenges have been addressed in \cite{BGN_2015_navierstokes}, which is based on previous work on geometric evolution equations such as \cite{BGN_2007, BGN_2008_pfem_three_dimensions}. Here, the interface is moved in the normal direction with the normal component of the fluid velocity, while the tangential degrees of freedom are implicitly used to maintain a good mesh quality, without the need for re-parametrisations or other manipulations. This strategy can be applied to both unfitted and fitted parametric methods \cite{GNZ_2023}.

Another difficulty for parametric methods is to find a good approximation of the pressure, which is discontinuous across the interface due to capillary effects. Here, a fitted bulk mesh that is adapted to the interface can be used \cite{ganesan_2007}. 
In the case of an unfitted bulk mesh, where the interface and bulk meshes are totally independent, it is helpful to use an extended finite element method (XFEM) or a cut finite element method (cutFEM) to accurately capture jumps of physical quantities across the interface \cite{claus_kerfriden_2019_cutfem, zahedi_2019_cutfem, gross_reusken_2007}.
We achieve this with the XFEM approach of \cite{BGN_2013_stokes} by enriching the pressure space with exactly one additional degree of freedom, together with a refinement of the bulk mesh close to the interface. The XFEM enrichment has the additional benefit of exact volume conservation of the phases on the semi-discrete level, whereas the fully discrete scheme shows good volume conservation properties in practice.

In this work, we combine the ideas of the parametric method from \cite{BGN_2015_navierstokes} for the two-phase Navier--Stokes flow and the finite element method from \cite{barrett_boyaval_2009} for the one-phase viscoelastic Oldroyd-B model. This allows us to construct a numerical scheme for viscoelastic two-phase flow with unconditional stability and existence of discrete solutions and good parametric mesh properties.
\cblue{A key challenge in this work is the numerical analysis of the nonlinear coupling between the viscoelastic system and the interface dynamics. Unlike previous works on two-phase flows \cite{BGN_2013_stokes, BGN_2015_navierstokes}, which resulted only in linear systems and used arguments from linear algebra for existence and uniqueness, our approach requires a fixed-point framework similar to \cite{barrett_boyaval_2009}. However, directly applying the methods from \cite{barrett_boyaval_2009} is not possible because the available \textit{a priori} estimates (see Lemma \ref{lemma:stability_delta}, below) do not provide explicit information about the discrete curvature, which is one of the unknowns in our discrete scheme.
To overcome this, we use a Schur complement approach.
This allows the fixed-point argument to focus only on the bulk quantities. Careful handling of these arguments is essential, as the proof of the main result involves the limit passing in a regularized discrete scheme, which is presented in Section \ref{sec:stability_proof}.}

The remainder of this work is organised as follows. In Section \ref{sec:model}, we present a mathematical model for viscoelastic two-phase flow in a bounded domain with a moving free interface, and we derive the necessary jump and continuity conditions at the interface using an energy identity. We introduce a weak formulation and discuss key properties, such as an energy dissipation law and phase volume conservation. In Section \ref{sec:fem}, we develop a stable fully discrete numerical scheme using the parametric finite element method in space and a semi-implicit discretisation in time. We prove the unconditional stability and existence of solutions. The stability proof is based on an approximative solution layer with a regularised elastic energy and a limit passage on the fully discrete level.
Using a semi-discrete variant, we discuss further properties of the numerical method, such as beneficial parametric mesh quality and an XFEM approach for good phase volume conservation, which we also observe in practice. We also discuss an extension of the scheme to a variable shear modulus. In Section \ref{sec:numerics}, we present \cred{numerical results} 
in two dimensions, demonstrating the practicability of the introduced scheme.

\section{The two-phase Navier--Stokes problem with viscoelasticity}
\label{sec:model}

\begin{figure}
\centering
\def\svgwidth{0.5\linewidth}
\begingroup%
  \makeatletter%
  \providecommand\color[2][]{%
    \errmessage{(Inkscape) Color is used for the text in Inkscape, but the package 'color.sty' is not loaded}%
    \renewcommand\color[2][]{}%
  }%
  \providecommand\transparent[1]{%
    \errmessage{(Inkscape) Transparency is used (non-zero) for the text in Inkscape, but the package 'transparent.sty' is not loaded}%
    \renewcommand\transparent[1]{}%
  }%
  \providecommand\rotatebox[2]{#2}%
  \newcommand*\fsize{\dimexpr\f@size pt\relax}%
  \newcommand*\lineheight[1]{\fontsize{\fsize}{#1\fsize}\selectfont}%
  \ifx\svgwidth\undefined%
    \setlength{\unitlength}{361.03479767bp}%
    \ifx\svgscale\undefined%
      \relax%
    \else%
      \setlength{\unitlength}{\unitlength * \real{\svgscale}}%
    \fi%
  \else%
    \setlength{\unitlength}{\svgwidth}%
  \fi%
  \global\let\svgwidth\undefined%
  \global\let\svgscale\undefined%
  \makeatother%
  \begin{picture}(1,0.70758922)%
    \lineheight{1}%
    \setlength\tabcolsep{0pt}%
    \put(0,0){\includegraphics[width=\unitlength,page=1]{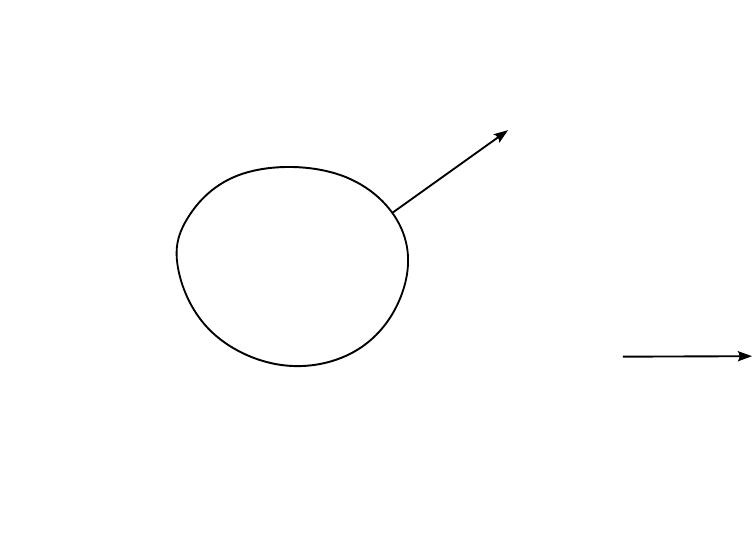}}%
    \put(0.32287917,0.50849018){\color[rgb]{0,0,0}\makebox(0,0)[lt]{\lineheight{1.25}\smash{\begin{tabular}[t]{l}$\Gamma(t)$\end{tabular}}}}%
    \put(0.60849288,0.43783838){\color[rgb]{0,0,0}\makebox(0,0)[lt]{\lineheight{1.25}\smash{\begin{tabular}[t]{l}$\pmb \nu(t)$\end{tabular}}}}%
    \put(0.86554513,0.26346371){\color[rgb]{0,0,0}\makebox(0,0)[lt]{\lineheight{1.25}\smash{\begin{tabular}[t]{l}$\mathbf{n}$\end{tabular}}}}%
    \put(0.01020741,0.56410967){\color[rgb]{0,0,0}\makebox(0,0)[lt]{\lineheight{1.25}\smash{\begin{tabular}[t]{l}$\partial\Omega$\end{tabular}}}}%
    \put(0.33640821,0.34163172){\color[rgb]{0,0,0}\makebox(0,0)[lt]{\lineheight{1.25}\smash{\begin{tabular}[t]{l}$\Omega_-(t)$\end{tabular}}}}%
    \put(0.39353104,0.09961178){\color[rgb]{0,0,0}\makebox(0,0)[lt]{\lineheight{1.25}\smash{\begin{tabular}[t]{l}$\Omega_+(t)$\end{tabular}}}}%
    \put(0,0){\includegraphics[width=\unitlength,page=2]{figures/img_domain_square.pdf}}%
  \end{picture}%
\endgroup%

\caption{A schematic sketch of the domain $\Omega$ in two space dimensions.}
\label{fig:domain}
\end{figure}

We now present the mathematical model for viscoelastic two-phase flow.
Let $\Omega\subset\bbR^{d}$, $d\in\{2,3\}$, be a bounded domain with boundary $\partial\Omega$ and let $t\in(0,T)$, where $T>0$. 
We assume that the domain $\Omega$ contains two distinct incompressible immiscible viscoelastic fluid phases that occupy time-dependent regions $\Omega_+(t)$ and $\Omega_-(t)\coloneqq \Omega \setminus \overline{\Omega_+(t)}$ for all $t\in[0,T]$.
The two phases are separated by a time-evolving interface $(\Gamma(t))_{t\in[0,T]}$ which is contained within $\Omega$ and does not intersect the boundary $\partial\Omega$. 
We assume that the interface $(\Gamma(t))_{t \in [0, T]}$ is an orientable and smooth $(d-1)$-dimensional hypersurface without boundary, globally parametrised by the mapping $\bX(\cdot, t) \colon\, \Upsilon \to \mathbb{R}^d$, where $\Upsilon \subset \mathbb{R}^d$ is a predefined reference manifold. Therefore, $\Gamma(t) = \bX(\Upsilon, t)$.
Then, $\calV \coloneqq \partial_t \bX \cdot \pmb\nu$ denotes the normal velocity of the evolving hypersurface, where $\pmb\nu$ is the unit normal on $\Gamma(t)$ pointing into $\Omega_+(t)$. Moreover, we denote the outer unit normal on $\partial\Omega$ by $\bn$.
We refer to Figure \ref{fig:domain} for a schematic sketch in two dimensions.  

In the subdomains $\Omega_\pm(t)$, we study the incompressible Navier--Stokes equations for the velocity $\bu(\cdot,t)\colon\, \Omega \to \bbR^d$ and the pressure $p(\cdot,t)\colon\, \Omega\to\bbR$, along with an additional stress tensor $\bbT_e(\cdot,t)\colon\, \Omega\to\bbR^{d\times d}$ that incorporates fluid elasticity effects:
\begin{subequations} 
\label{eq:system_bulk}
\begin{align}
    \label{eq:u}
    \rho_\pm \partial_t \bu 
    + \rho_\pm (\bu \cdot\nabla) \bu
    + \nabla p &= \Div(2 \mu_\pm \bbD(\bu)) + \Div(\bbT_e) 
    + \cred{\rho_\pm} \mathbf{f}_1 + \mathbf{f}_2
    \qquad &&\text{in } \Omega_\pm(t),
    \\
    \label{eq:div}
    \Div \bu &= 0
    \qquad &&\text{in } \Omega_\pm(t).
\end{align}
Here, we denote the mass densities and viscosities of the phases $\Omega_\pm(t)$ by $\rho_\pm>0$ and $\mu_\pm>0$, respectively. 
For future reference, we define $\rho(\cdot,t) \coloneqq \rho_+ (1-\chi_{|\Omega_-(t)}) + \rho_- \; \chi_{|\Omega_-(t)}$ and $\mu(\cdot,t) \coloneqq \mu_+ (1-\chi_{|\Omega_-(t)}) + \mu_- \; \chi_{|\Omega_-(t)}$, where $\chi_{|\Omega_-(t)}$ is the characteristic function of $\Omega_-(t)$.
Moreover, the symmetrised velocity gradient is denoted by $\bbD(\bu)\coloneqq \frac12 (\nabla\bu + (\nabla\bu)^\top)$, and $\mathbf{f}_1(\cdot,t), \mathbf{f}_2(\cdot,t)\colon\, \Omega\to\bbR^d$ describe body acceleration and external forces, respectively. 
%

As indicated before, in this paper we are going to consider the Oldroyd-B model to describe the viscoelastic effects in the two phases. In particular, the elastic stress tensor is given by 
\begin{align*}
    \bbT_e = G (\bbB-\bbI),
\end{align*}
where $G>0$ is the (constant) shear modulus, $\bbI\in\bbR^{d\times d}$ denotes the identity matrix and $\bbB(\cdot,t)\colon\, \Omega\to\bbR^{d\times d}_{\mathrm{S}}$ is a tensor which satisfies 
\begin{align}
    \label{eq:B}
    \partial_t \bbB + (\bu\cdot\nabla)\bbB
    + \frac{1}{\lambda_\pm} (\bbB-\bbI)
    &= 
    \nabla\bu \bbB + \bbB (\nabla\bu)^\top
    + \alpha \Delta \bbB
    \qquad &&\text{in } \Omega_\pm(t).
\end{align}
\end{subequations}
Here, $\alpha\geq0$ is a constant that accounts for possible stress diffusion, $\lambda_\pm>0$ are the relaxation parameters for the viscoelastic fluid phases and $\bbR^{d\times d}_\mathrm{S}$ is the set of all symmetric $(d\times d)$-matrices. Similarly to before, we define $\lambda(\cdot,t) \coloneqq \lambda_+ (1-\chi_{|\Omega_-(t)}) + \lambda_- \; \chi_{|\Omega_-(t)}$. 

We point out that the Oldroyd-B equation \eqref{eq:B} with $\alpha=0$ serves as a model that approximates fluid-structure interaction problems \cite{mokbel_abels_aland_2018}. Sending the relaxation parameter $\lambda_\pm$ to zero in one phase results in a representation of a Newtonian fluid where $\bbB=\bbI$, indicating that the stress tensor does not have an elastic contribution. In contrast, increasing the relaxation parameter to infinity in another phase suggests that \eqref{eq:B} governs the dynamics of the left Cauchy--Green tensor in a solid, which may exhibit viscoelastic properties if $\mu>0$, or be entirely elastic when $\mu=0$.
In this context, the tensor $\bbB$ can be interpreted as the left Cauchy--Green tensor associated with the purely elastic part of the total mechanical response of the fluid, also known as the \textit{elastic Cauchy--Green tensor} \cite{malek_2018_Oldroyd_diffusive}. 

The case $\alpha>0$ can be seen as a mathematical regularisation, improving the regularity of the solutions. On the other hand, there is a physical justification for the additional diffusive term $\alpha\Delta\bbB$ in \eqref{eq:B}, taking the macroscopic closure of specific particle-based systems \cite{barrett_sueli_2007}.

\subsection{Conditions at the free interface and on the boundary}
We now derive conditions which need to hold on the boundary $\partial\Omega$ and at the interface $\Gamma(t)$ using an energy identity.
For future reference, we recall the following transport identities from \cite[Thm.~33, Thm.~35]{BGN_2019_handbook}:
\begin{align}
\label{eq:transport_bulk}
    \ddt \int_\Omega b \dx 
    &= \int_{\Omega} \partial_t b \dx 
    - \int_{\Gamma(t)} \calV \, \jump{b}_-^+ \ds,
    \\
    \label{eq:transport_interface}
    \ddt \int_{\Gamma(t)} 1 \ds 
    &= - \int_{\Gamma(t)} \kappa \calV \ds,
    \\
    \label{eq:transport_subdomain}
    \ddt \int_{\Omega_{-}(t)} 1 \dx
    &= \int_{\Gamma(t)} {\calV} \ds,
\end{align}
where $\kappa$ denotes the mean curvature of the interface $\Gamma(t)$, i.e., the sum of the principal curvatures of $\Gamma(t)$, with the convention that $\kappa$ is negative where $\Omega_-(t)$ is locally convex. We also define the jump across the interface in the normal direction by $\jump{b}_-^+ \coloneqq b_+ - b_-$, where $b_\pm \coloneqq b_{|\Omega_\pm(t)}$ for scalar, vector or matrix valued functions $b\colon\, \Omega\times(0,T)\to X$ with $X\in\{\bbR, \bbR^d, \bbR^{d\times d}\}$.

Let us introduce the energy of the two-phase viscoelastic model \eqref{eq:system_bulk} with free interface $\Gamma(t)$:
\begin{align}
    \label{eq:energy}
    E(\bu,\bbB,\Gamma(t)) = 
    \int_\Omega \Big( \frac{\rho}{2} \abs{\bu}^2
    + \frac{G}{2} \trace(\bbB - \ln\bbB - \bbI) \Big) \dx
    + \int_{\Gamma(t)} \gamma \ds ,
\end{align}
which consists of the kinetic, the elastic and the interfacial energy with surface tension $\gamma>0$.
The elastic energy density $\frac{G}{2} \trace(\bbB - \ln\bbB - \bbI)$ has been introduced in \cite{Hu_Lelievre_2007} and will ensure that $\bbB$ remains positive definite, provided that it is positive definite at the initial time. We point out the identity $\trace\ln\bbB = \ln\det\bbB$ for positive definite symmetric matrices by applying the logarithmic function to the spectrum of $\bbB$. From a physical point of view, the positive definiteness of $\bbB$ is an important property to avoid compression to a single point or (local) self-interpenetration within the material \cite{kruzik_roubicek_2019}. 

Taking the time derivative of the energy \eqref{eq:energy} and using the transport identities \eqref{eq:transport_bulk} and \eqref{eq:transport_interface}, we obtain
\begin{align*}
    \ddt E(\bu, \bbB, \Gamma(t))  
    & = 
    \int_{\Omega} \partial_t \big( \frac{\rho}{2} \abs{\bu}^2
    + \frac{G}{2} \trace(\bbB - \ln\bbB - \bbI) \big) \dx
    \\
    &\quad
    - \int_{\Gamma(t)} \calV \, \jumpp{\frac{\rho}{2} \abs{\bu}^2 + \frac{G}{2} \trace(\bbB-\ln\bbB-\bbI)}_-^+ \ds
    - \int_{\Gamma(t)} \gamma \kappa \calV \ds.
\end{align*}
Taking the dot product of \eqref{eq:u} with $\bu$, integrating over $\Omega$ and using that $\rho$ is constant in the subdomains $\Omega_\pm(t)$, we obtain
\begin{align*}
    \int_{\Omega} 
    \partial_t \big(\frac{\rho}{2} \abs{\bu}^2 \big)  \dx
    =
    \int_{\Omega} \Big(  
    - \bu\cdot\nabla \big( \frac{\rho}{2} \abs{\bu}^2 \big) 
    + \Div( 2 \mu \bbD(\bu) - p \bbI + \bbT_e) \cdot\bu
    + ( \rho \mathbf{f}_1 + \mathbf{f}_2 ) \cdot \bu \Big) \dx.
\end{align*}
Integrating by parts over the subdomains $\Omega_\pm(t)$ yields
\begin{align*}
    -\int_{\Omega} 
    \bu\cdot\nabla \big( \frac{\rho}{2} \abs{\bu}^2 \big) \dx
    &=  \int_{\Omega} \frac{\rho}{2} \abs{\bu}^2 \Div(\bu)  \dx
    - \int_{\partial\Omega} \frac{\rho}{2} \abs{\bu}^2 \bu\cdot\bn  \ds
    \\
    &\quad
    + \int_{\Gamma(t)} \jumpp{\frac{\rho}{2} \abs{\bu}^2 \bu\cdot\pmb\nu }_-^+ \ds,
\end{align*}
and 
\begin{align*}
    &\int_{\Omega} \Div( 2 \mu \bbD(\bu) - p \bbI + \bbT_e) \cdot\bu \dx
    \\
    &=
    -\int_{\Omega}
    \big( 2 \mu \abs{\bbD(\bu)}^2  - p \Div(\bu) + \bbT_e : \nabla\bu \big) \dx
    + \int_{\partial\Omega} \bu^\top \big(2 \mu \bbD(\bu) - p \bbI + \bbT_e\big) \bn \ds
    \\ 
    &\quad
    - \int_{\Gamma(t)} \jumpp{ \bu^\top \big(2 \mu \bbD(\bu) - p \bbI + \bbT_e\big) \,\pmb\nu}_-^+ \ds,
\end{align*}
where we used $\bbD(\bu) : \nabla\bu = \bbD(\bu) : \bbD(\bu) = \abs{\bbD(\bu)}^2$ and $p\bbI : \nabla\bu = p \Div(\bu)$. 
By taking the Frobenius inner product of \eqref{eq:B} with $\frac{G}{2} (\bbI - \bbB^{-1})$ and integrating over $\Omega$, we get
\begin{align*}
    \int_{\Omega} \partial_t \big( \frac{G}{2} \trace(\bbB-\ln\bbB-\bbI) \big)\dx
    &= - \int_{\Omega} \Big( \bu\cdot\nabla \big( \frac{G}{2} \trace(\bbB-\ln\bbB-\bbI) \big)
    - \bbT_e : \nabla\bu \Big) \dx
    \\
    &\quad
    - \int_{\Omega} \Big( \frac{G}{2\lambda} \trace(\bbB + \bbB^{-1} - 2 \bbI ) 
    - \frac{\alpha G}{2} \Delta\bbB : (\bbI - \bbB^{-1}) \Big) \dx,
\end{align*}
where we used the symmetry of $\bbB$ and $(\bbI-\bbB^{-1})$ and elementary properties of the Frobenius inner product to write $\frac{G}{2}(\nabla \bu \bbB + \bbB (\nabla\bu)^\top) : (\bbI-\bbB^{-1}) = G (\bbB-\bbI) : \nabla\bu = \bbT_e : \nabla\bu$ and $(\bbB-\bbI) : (\bbI-\bbB^{-1}) = \trace(\bbB+\bbB^{-1}-2\bbI)$. Again, from integration by parts over the subdomains $\Omega_\pm(t)$, we have
\begin{align*}
    & -\int_{\Omega} \bu\cdot\nabla \big( \frac{G}{2} \trace(\bbB-\ln\bbB-\bbI) \big) \dx
    \\
    &= \int_{\Omega} \frac{G}{2} \trace(\bbB-\ln\bbB-\bbI) \Div(\bu)  \dx
    - \int_{\partial\Omega} \frac{G}{2} \trace(\bbB-\ln\bbB-\bbI) \bu\cdot\bn  \ds
    \\
    &\quad
    + \int_{\Gamma(t)} \jumpp{\frac{G}{2} \trace(\bbB-\ln\bbB-\bbI) \bu\cdot\pmb\nu }_-^+ \ds,
\end{align*}
and
\begin{align*}
    & \int_{\Omega} \frac{\alpha G}{2} \Delta\bbB : (\bbI - \bbB^{-1}) \dx
    \\
    &= \int_{\Omega} \frac{\alpha G}{2} \nabla\bbB : \nabla \bbB^{-1} \dx
    + \int_{\partial\Omega} \frac{\alpha G}{2} (\bn\cdot\nabla) \bbB : (\bbI-\bbB^{-1}) \ds
    \\
    &\quad
    - \int_{\Gamma(t)}  \jumpp{\frac{\alpha G}{2} (\pmb\nu\cdot\nabla) \bbB :  (\bbI-\bbB^{-1}) }_-^+ \ds.
\end{align*}
Collecting the computations from above and using \eqref{eq:div}, i.e., $\Div(\bu)=0$, we get the general energy identity
\begin{align}
    \label{eq:energy_formal} \nonumber
    &\ddt E(\bu, \bbB, \Gamma(t)) 
    \\ \nonumber
    &= - \int_{\Omega}
    \Big( 2 \mu \abs{\bbD(\bu)}^2 + \frac{G}{2\lambda} \trace(\bbB+\bbB^{-1}-2\bbI) 
    - \frac{\alpha G}{2} \nabla\bbB : \nabla\bbB^{-1} \Big) \dx
    + \int_\Omega  ( \rho \mathbf{f}_1 + \mathbf{f}_2 ) \cdot\bu \dx
    \\ \nonumber
    &\quad 
    + \int_{\partial\Omega}  \Big(
    \bu^\top \big( 2\mu\bbD(\bu) - p\bbI  + \bbT_e \big) \bn  
    + \frac{\alpha G}{2} (\bn\cdot\nabla)\bbB : (\bbI-\bbB^{-1})
    \Big)\dx
    \\ \nonumber
    &\quad
    - \int_{\partial\Omega} \big( \frac{\rho}{2} \abs{\bu}^2 + \frac{G}{2} \trace(\bbB-\ln\bbB-\bbI) \big) \bu\cdot\bn  \ds
    \\ \nonumber
    &\quad
    + \int_{\Gamma(t)} \Big( - \gamma \kappa \calV  - \jumpp{ \bu^\top \big(2 \mu \bbD(\bu) - p \bbI + \bbT_e\big) \,\pmb\nu}_-^+   \Big) \ds 
    \\
    &\quad
    + \int_{\Gamma(t)} 
   \Big( \jumpp{\big( \frac{\rho}{2} \abs{\bu}^2 
   + \frac{G}{2} \trace(\bbB-\ln\bbB-\bbI) \big) \, (\bu\cdot\pmb\nu - \calV) }_-^+ 
   - \jumpp{\frac{\alpha G}{2} (\pmb\nu\cdot\nabla) \bbB :  (\bbI-\bbB^{-1}) }_-^+
   \Big) \ds.
\end{align}
If $\bbB$ is positive definite, we point out that it holds
\begin{align}
    \label{eq:nabla_ln_formal}
    -\nabla\bbB:\nabla\bbB^{-1} \geq \frac{1}{d} \abs{\nabla\trace\ln\bbB}^2, 
    \qquad 
    \trace(\bbB+\bbB^{-1}-2\bbI) = \abs{\bbB^{1/2}-\bbB^{-1/2}}^2 \geq 0,
\end{align}
see \cite[Lem.~3.1]{barrett_lu_sueli_2017} for the first statement in \eqref{eq:nabla_ln_formal}.
Hence, the energy identity \eqref{eq:energy_formal} motivates the application of specific conditions at the boundary $\partial\Omega$ and the interface $\Gamma(t)$ to ensure that all integrals across $\partial\Omega$ and $\Gamma(t)$ on the right-hand side of \eqref{eq:energy_formal} are zero. At the boundary $\partial\Omega$, we set
\begin{subequations}
\label{eq:bc_all}
\begin{align}
    \label{eq:bc_T}
    (2 \mu\bbD(\bu) - p \bbI + \bbT_e) \bn 
    &= \mathbf{0}, \quad \bu\cdot\bn=0
    \qquad && \text{on } \partial_{\mathrm{N}}\Omega \times (0,T),
    \\
    \label{eq:bc_u}
    \bu&= \mathbf{0} \qquad && \text{on } \partial_{\mathrm{D}}\Omega \times (0,T), 
    \\
    \label{eq:bc_B}
    \alpha (\bn\cdot\nabla)\bbB &= \bbO \qquad &&\text{on } \partial\Omega \times (0,T),
\end{align}
\end{subequations}
where $\mathbf{0}\in\bbR^d$ and $\bbO\in\bbR^{d\times d}$ denote the zero vector and zero matrix, respectively.
Here, $\partial_{\mathrm{D}}\Omega \subset \partial\Omega$ is assumed to be closed and have positive surface measure $\calH^{d-1}(\partial_{\mathrm{D}}\Omega) > 0$, and we define $\partial_\mathrm{N}\Omega = \partial\Omega \setminus \partial_{\mathrm{D}}\Omega$.
At the interface $\Gamma(t)$, we prescribe the jump conditions
\begin{subequations}
\label{eq:jump_all}
\begin{align} \label{eq:jump}
    \jump{\bu}_{-}^+ &= \mathbf{0},
    \qquad - \jump{2 \mu \bbD(\bu) - p \bbI + \bbT_e}_{-}^+ \,\pmb\nu = \gamma \kappa \,\pmb\nu,
    \\
    \label{eq:jump2}
    \alpha \jump{\bbB}_-^+ &= \bbO,
    \qquad \alpha \jumpp{ G (\pmb\nu \cdot \nabla) \bbB}_-^+ = \bbO,
\end{align}
\end{subequations}
and, for the interface evolution, we assume that the normal velocity satisfies
\begin{align}
	\label{eq:X}
	\calV = \bu\cdot\pmb\nu \qquad \text{ on } \Gamma(t).
\end{align}
Note that in the case $\alpha=0$, the conditions \eqref{eq:bc_B} and \eqref{eq:jump2} are not needed.
We close the viscoelastic two-phase system \eqref{eq:system_bulk}, \eqref{eq:X}, \eqref{eq:bc_all}, \eqref{eq:jump_all} by fixing the initial conditions 
\begin{align*}
    \Gamma(0) = \Gamma_0, \quad \bu(\cdot,0) = \bu_0(\cdot), \quad \bbB(\cdot,0) = \bbB_0(\cdot) \quad \text{in } \Omega,
\end{align*}
where $\Gamma_0$, $\bu_0$ and $\bbB_0$ are given initial data. We assume throughout that $\bbB_0$ is symmetric, i.e., $\bbB_0 \colon\, \Omega\to \bbR^{d\times d}_\mathrm{S}$, as \eqref{eq:B} preserves the symmetry of $\bbB$ for symmetric initial data.

We remark that in the absence of viscoelastic effects, our model describes the two-phase Navier--Stokes flow.
Similarly to, e.g., \cite{baensch_2001, BGN_2015_navierstokes, ganesan_2007, gross_reusken_2011}, we have to compute surface tension forces and consequently curvature quantities, see \eqref{eq:jump}. Therefore, it is essential to find a suitable description for the mean curvature of the interface. 
In this context, the equation 
\begin{align}
    \label{eq:kappa}
    \kappa \,\pmb\nu &= \Delta_s \mathbf{id}  \qquad \text{on } \Gamma(t),
\end{align}
represents the relationship between the mean curvature vector $\kappa \,\pmb\nu$ and the Laplace--Beltrami operator $\Delta_s = \Div_s \circ \nabla_s$ applied to the identity map on the manifold $\Gamma(t)$, with $\Div_s$ and $\nabla_s$ denoting surface divergence and surface gradient on $\Gamma(t)$, respectively.
Initially introduced by \cite{dziuk_1991} for approximating curvature-driven interface evolution through linear finite elements, this idea was expanded in subsequent studies by \cite{baensch_2001, ganesan_2007, gross_reusken_2011} and used in, e.g., \cite{BGN_2013_stokes, BGN_2015_navierstokes, GNZ_2023}. Likewise, the approach of the current paper to approximate the mean curvature also depends on a weak formulation of the identity \eqref{eq:kappa}.

%
%
\subsection{Variational formulation of the model}
We now state a variational formulation, which is the starting point for the finite element approximation of the viscoelastic two-phase model. 
We multiply the equations \eqref{eq:system_bulk}, \eqref{eq:X}, \eqref{eq:kappa} with time-dependent test functions, integrate over $\Omega$ or, respectively, $\Gamma(t)$, integrate by parts and use the interface conditions \eqref{eq:jump_all}, the boundary conditions \eqref{eq:bc_all} and the transport identity \eqref{eq:transport_bulk}. We only note the calculation
\begin{align}
    \label{eq:B_weak_derivation} \nonumber
    &\int_\Omega \partial_t \bbB : \bbG \dx 
    + \int_\Omega (\bu\cdot\nabla)\bbB : \bbG 
    \\ \nonumber
    &= \ddt \int_\Omega \bbB : \bbG \dx
    - \int_\Omega \bbB : \partial_t \bbG \dx
    - \int_\Omega \bbB : (\bu\cdot\nabla)\bbG \dx
    \\
    &\quad
    - \int_\Omega \Div(\bu) \, \bbB : \bbG \dx
    + \int_{\partial\Omega} \bu\cdot\bn\, \bbB : \bbG \ds
    + \int_{\Gamma(t)} \big( \calV\, \jump{\bbB:\bbG}_-^+ - \jump{ \bu \cdot \pmb\nu\, \bbB:\bbG}_-^+ \big) \ds,
\end{align}
where $\bbG(\cdot,t) \in H^1(\Omega;\bbR^{d\times d}_\mathrm{S})$ is a time-dependent test function. Note that the last three terms in \eqref{eq:B_weak_derivation} vanish due to \eqref{eq:div}, \eqref{eq:bc_all}, \eqref{eq:jump} and \eqref{eq:X}. 

Subsequently, the variational formulation is as follows. Find time-dependent functions $\bu$, $p$, $\bbB$, $\bX$ and $\kappa$ such that $\bu(\cdot,t) \in H^1_\mathrm{D}(\Omega;\bbR^d)$, $p(\cdot,t) \in L^2(\Omega)$, $\bbB(\cdot,t) \in H^1(\Omega;\bbR^{d\times d}_\mathrm{S})$, $\bX(\cdot,t) \in H^1(\Upsilon;\bbR^d)$ and $\kappa(\cdot,t) \in L^2(\Upsilon)$ with
\begin{subequations}
\label{eq:system_weak}
\begin{align}
    \label{eq:u_weak} \nonumber
    0 &=  \ddt  \int_\Omega \frac12 \rho \bu \cdot \bw \dx
    + \int_\Omega \frac12 \rho \partial_t \bu\cdot\bw \dx
    - \int_\Omega \frac12 \rho \bu \cdot \partial_t \bw \dx
    + \int_\Omega \frac12 \rho (\bu\cdot\nabla)\bu \cdot\bw \dx
    \\ \nonumber
    &\quad 
    - \int_\Omega \frac12 \rho \bu \cdot (\bu\cdot\nabla)\bw \dx
    + \int_\Omega 2 \mu \D(\bu) : \D(\bw) \dx
    - \int_\Omega p \Div(\bw) \dx
    \\
    &\quad
    + \int_\Omega G(\bbB-\bbI) : \nabla\bw \dx
    - \int_{\Gamma(t)} \gamma \kappa \, \pmb\nu \cdot \bw \ds
    - \int_\Omega ( \rho \mathbf{f}_1 + \mathbf{f}_2 ) \cdot \bw ,
    \\ \label{eq:div_weak}
    0 &= \int_\Omega \Div(\bu) q \dx, 
    \\ \label{eq:B_weak} \nonumber
    0 &= \ddt \int_\Omega \bbB : \bbG \dx 
    - \int_\Omega \bbB : \partial_t \bbG \dx 
    - \int_\Omega \bbB : (\bu\cdot\nabla)\bbG \dx
    + \int_\Omega {\frac{1}{\lambda} (\bbB - \bbI)} : {\bbG} \dx
    \\
    &\quad
    - \int_\Omega 2 {\nabla\bu} : {(\bbG\bbB)} \dx
    + \int_\Omega \alpha {\nabla\bbB} : {\nabla\bbG} \dx,
    \\ \label{eq:X_weak} 
    0 &= \int_{\Gamma(t)} \partial_t\bX \cdot \pmb \nu \, \chi \ds
    - \int_{\Gamma(t)} \bu \cdot\pmb\nu \,\chi \ds,
    \\ \label{eq:kappa_weak}
    0 &= \int_{\Gamma(t)} \kappa \, \pmb\nu \cdot \pmb\zeta \ds
    + \int_{\Gamma(t)} \nabla_s \bX : \nabla_s \pmb\zeta \ds,
\end{align}
\end{subequations}
for a.e.~$t\in(0,T)$ and for all time-dependent test functions $\bw(\cdot,t)\in H^1_\mathrm{D}(\Omega;\bbR^d)$, $q(\cdot,t)\in L^2(\Omega)$, $\bbG(\cdot,t)\in H^1(\Omega;\bbR^{d\times d}_{\mathrm{S}})$, $\chi(\cdot,t)\in L^2(\Upsilon)$ and $\pmb\zeta(\cdot,t) \in H^1(\Upsilon; \bbR^d)$. 
Here we define 
\begin{align*}
    H^1_\mathrm{D}(\Omega;\bbR^d) &\coloneqq \{ \bz \in H^1(\Omega;\bbR^d) \mid 
    \bz = \mathbf{0} \text{ on } \partial_\mathrm{D}\Omega, \ 
    \bz\cdot\bn=\mathbf{0} \text{ on } \partial_\mathrm{N}\Omega\}.
\end{align*}


In the case without viscoelasticity, e.g., by setting $G=0$, the system \eqref{eq:system_weak} coincides with the weak formulation from \cite{BGN_2015_navierstokes} for two-phase Navier--Stokes flow.

Similarly to \cite{BGN_2013_stokes, BGN_2015_navierstokes}, we employ a slight abuse of notation in \eqref{eq:u_weak}, \eqref{eq:X_weak} and \eqref{eq:kappa_weak}, in the sense that here and throughout this paper, functions defined on the reference manifold $\Upsilon$ are identified with functions defined on $\Gamma(t)$.

\subsection{Energy identity and volume conservation}
For completeness, we now present an energy identity for sufficiently smooth solutions to the variational formulation \eqref{eq:system_weak} with $\bbB$ being positive definite, which follows directly from \eqref{eq:energy_formal} by applying the boundary conditions \eqref{eq:bc_all} and the interface conditions \eqref{eq:jump_all} and \eqref{eq:X}.
%
Let $(\bu,p,\bbB,\bX, \kappa)$ denote a smooth solution to \eqref{eq:system_weak} and assume that $\bbB$ is positive definite. Then, it holds for a.e.~$t\in(0,T)$ that
\begin{align}
    \label{eq:stability} \nonumber
    &\ddt \int_\Omega \Big( \frac{\rho}{2} \abs{\bu}^2
    + \frac{G}{2} \trace(\bbB - \ln\bbB - \bbI) \Big) \dx
    + \ddt \int_{\Gamma(t)} \gamma \ds
    \\ \nonumber
    &\quad
    + \int_{\Omega}
    \Big( 2 \mu \abs{\bbD(\bu)}^2 + \frac{G}{2\lambda} \trace(\bbB+\bbB^{-1}-2\bbI) 
    - \frac{\alpha G}{2} \nabla\bbB : \nabla\bbB^{-1} \Big) \dx
    \\
    &= \int_\Omega  ( \rho \mathbf{f}_1 + \mathbf{f}_2 ) \cdot\bu \dx .
\end{align}

A different method to obtain \eqref{eq:stability} involves applying energy estimates to the variational formulation \eqref{eq:system_weak}, similar to our approach in the finite element method. 
Specifically, \eqref{eq:stability} is derived by choosing the test functions $\bw=\bu$, $q=p$, $\bbG = \tfrac{G}{2} (\bbI - \bbB^{-1})$, $\chi=-\gamma\kappa$ and $\pmb\zeta = \gamma \partial_t \bX$ within \eqref{eq:system_weak}.
%
Let us only mention the geometric identity 
\begin{align}
    \label{eq:stability_formal_1}
    \ddt \int_{\Gamma(t)} 1 \ds 
    = \int_{\Gamma(t)}  {\nabla_s \bX} : {\nabla_s \partial_t\bX} \ds,
\end{align}
and the computation
\begin{align}
    \label{eq:stability_formal_2} \nonumber
	&\ddt \int_\Omega  \bbB : (\bbI-\bbB^{-1}) \dx 
	+ \int_\Omega \bbB : \partial_t \bbB^{-1} \dx 
	+ \int_\Omega  \bbB : (\bu\cdot\nabla)\bbB^{-1} \dx 
	\\
	&= \ddt \int_\Omega  \trace(\bbB-\bbI)\dx 
	- \int_\Omega  \partial_t \trace\ln\bbB \dx
	- \int_\Omega  (\bu\cdot\nabla) \trace\ln\bbB \dx,
\end{align}
which are based on the transport identities \eqref{eq:transport_bulk}--\eqref{eq:transport_interface}, integration by parts over the subdomains $\Omega_\pm(t)$ and \eqref{eq:kappa}.
Similarly, using \eqref{eq:transport_bulk}, integration by parts over the subdomains $\Omega_\pm(t)$, \eqref{eq:div_weak} with $q= \frac{1}{2}\trace\ln\bbB$ and \eqref{eq:X_weak} with $\chi= \jump{\frac12 \trace\ln\bbB}_-^+$, we note that
\begin{align}
    \label{eq:stability_formal_3} \nonumber
	&\int_\Omega  \partial_t \trace\ln\bbB \dx
	+ \int_\Omega  (\bu\cdot\nabla) \trace\ln\bbB \dx
	\\ \nonumber
	&= \ddt \int_\Omega  \trace\ln\bbB \dx 
	- \int_\Omega  \Div(\bu)  \trace\ln\bbB \dx
	+ \int_{\Gamma(t)} (\calV - \bu\cdot \pmb\nu) \jump{ \trace\ln\bbB}_-^+ \ds
	\\
	&= \ddt \int_\Omega  \trace\ln\bbB \dx .
\end{align}
Combining \eqref{eq:stability_formal_1}, \eqref{eq:stability_formal_2} and \eqref{eq:stability_formal_3} with the aforementioned testing procedure then allows to derive \eqref{eq:stability} from \eqref{eq:system_weak}.

Another important property of all smooth solutions to \eqref{eq:system_weak} is the volume conservation of the bulk regions $\Omega_{\pm}(t)$. 
Using the transport identity \eqref{eq:transport_subdomain}, \eqref{eq:X_weak} with $\chi=1$, the Gauss theorem and \eqref{eq:div_weak} with $q=\chi_{|\Omega_{-}(t)}$, it holds
\begin{align} \label{eq:volume_conservation}
    \ddt \mathcal{L}^d({\Omega_{-}(t)})
    = \int_{\Gamma(t)} {\calV} \ds
    = \int_{\Gamma(t)} {\bu} \cdot {\pmb\nu} \ds
    = \int_{\Omega} {\Div(\bu)} \,{\chi_{|\Omega_{-}(t)}} \dx
    = 0,
\end{align}
where $\mathcal{L}^d(\cdot)$ denotes the $d$-dimensional Lebesgue measure.
The volume of $\Omega_+(t)$ does not change due to the fact that $\Omega_+(t) = \Omega \setminus \overline{\Omega_-(t)}$.
%
%


%
%
\section{Finite element approximation}\label{sec:fem}
The main goal of this section is to introduce and study a fully discrete finite element approximation based on the weak formulation \eqref{eq:system_weak} for which the existence of at least one discrete solution with a positive definite tensor $\bbB$ is guaranteed. Moreover, we show that any solution to the discrete system satisfies a stability inequality.

In the remainder of this work, we consistently assume the following.
Let $T>0$ and let $\Omega\subset\bbR^d$, $d\in\{2,3\}$, be a bounded Lipschitz domain with a polygonal (or polyhedral, respectively) boundary $\partial\Omega$. Moreover, we assume that $\partial_\mathrm{D}\Omega \subset \partial\Omega$ is closed and has a positive surface measure, and we set $\partial_\mathrm{N}\Omega = \partial\Omega \setminus \partial_\mathrm{D}\Omega$.
We divide the time interval $[0,T]$ into equidistant subintervals $[t^{n-1},t^n]$ with $t^n = n \Delta t$ and $t^{N_T}=T$, where $\Delta t = \frac{T}{N_T}$, $N_T\in\bbN$ and $n \in \{0, \ldots, N_T\}$.
For any $n\in\{0,\ldots,N_T-1\}$, we require $\calT^n$ to be a regular partitioning of $\Omega$ into disjoint open simplices such that $\overline\Omega = \bigcup_{K\in\calT^n} \overline K$, and we set $h^n \coloneqq \max_{K\in\calT^n} \diam(K)$.
We always assume that $\partial_{\mathrm{D}}\Omega$ and $\calT^n$ are given such that $\partial_{\mathrm{D}}\Omega$ is exactly matched by the sides in $\calT^n$. 
Moreover, we assume that the family of triangulations is shape regular, i.e., there exists a constant $C>0$ which is independent of $h^n$, $n\in\{0,\ldots,N_T\}$, such that
\begin{align*}
    \sup_{K\in\calT^n} \frac{\diam(K)}{ r_K } \leq C,
\end{align*}
where $r_K$ denotes the diameter of the largest inscribed ball in the simplex $K\in\calT^n$.
In the case $\alpha>0$ we also require that $\calT^n$ consists of non-obtuse simplices only, i.e., all dihedral angles of any simplex in $\calT^n$ are less than or equal to $\frac\pi2$.


For $k\in\bbN$ and $n\in\{0,\ldots,N_T-1\}$, we introduce the finite element space of continuous and piecewise polynomial functions
\begin{align*}
     \calS_k^n &\coloneqq \left\{ q \in C(\overline\Omega) \mid q_{|K} \in \calP_k(K) \ \forall \,  K\in\calT^n \right\} ,
\end{align*}
and the finite element space of discontinuous and piecewise constant functions
\begin{align*}
    \calS_0^n &\coloneqq \left\{ q \in L^1(\Omega) \mid q_{|K} \in \calP_0(K) \ \forall \,  K\in\calT^n \right\} .
\end{align*}
The standard nodal interpolation operators are denoted by $\mathrm{I}^n_k \colon\, C(\overline\Omega) \to \calS_k^n$ and their definitions are naturally extended to vector and matrix valued functions, i.e., $\bI^n_k \colon\, C(\overline\Omega;\bbR^d) \to [\calS_k^n]^d$ and \cblue{$\bbI^n_k \colon\, C(\overline\Omega;\bbR^{d\times d}_{\mathrm{S}})\to [\calS_k^n]^{d \times d}_{\mathrm{S}}$,} respectively.
Moreover, we recall the standard projection operator $\mathrm{I}^n_0 \colon\, L^1(\Omega)\to \calS_0^n$ mapping to the space of discontinuous and piecewise constant functions.

For the approximation space of the tensor $\bbB$, we choose $\bbS^n = [\calS_1^n]^{d \times d}_{\mathrm{S}}$. 
Let $(\bbU^n, \bbP^n)$ with $\bbU^n \subset H^1_\mathrm{D}(\Omega;\bbR^d)$ be a pair of finite element spaces on $\calT^n$ that satisfy the inf-sup condition
\begin{align}
    \label{eq:LBB}
    \inf_{q\in \hat\bbP^n \setminus\{0\}} 
    \sup_{\bw\in \bbU^n \setminus\{0\}}
    \frac{\skp{\Div\bw}{q} }{ \norm{\bw}_{H^1} \norm{q}_{L^2} }
    \geq C_0 ,
\end{align}
where $C_0>0$ is a constant independent of $h^n$, and $\hat\bbP^n = \bbP^n \cap L^2_0(\Omega)$, where $L^2_0(\Omega) = \{ q \in L^2(\Omega) \mid \int_\Omega q\dx = 0\}$.
Possible examples, which satisfy the inf-sup condition, are the lowest order Taylor--Hood element for $d=2,3$, also called the P2--P1 element, or the P2--(P1+P0) element for $d=2$, i.e.,
\begin{subequations}
\begin{align}
    \label{eq:stable_P2_P1}
    (\bbU^n, \bbP^n) 
    &= ( [\calS_2^n]^d \cap H^1_\mathrm{D}(\Omega;\bbR^d), 
    \, \calS_1^n), 
    \\
    \label{eq:stable_P2_P1+P0}
    (\bbU^n, \bbP^n) 
    &= ( [\calS_2^n]^d \cap H^1_\mathrm{D}(\Omega;\bbR^d), 
    \, \calS_1^n + \calS_0^n),
\end{align}
\end{subequations}
see, e.g., \cite[Chap.~4.2.5]{ern_guermond_2004} for \eqref{eq:stable_P2_P1} and \cite{boffi_2012_infsup} for \eqref{eq:stable_P2_P1+P0}. Other stable pairs such as the MINI element \cite{girault_raviart_1986} ($d=2, 3$) or the P2--P0 element \cite{boffi_2012_infsup} ($d=2$) are not considered in this work for the following reasons: 
The MINI element requires higher order quadrature rules to ensure exact integration, while the P2--P0 element does not contain P1 functions in the pressure space, which however are used in the stability estimates when dealing with the term $\bbB : (\bu\cdot\nabla)\bbG$ from \eqref{eq:B_weak}.
In principle, $\bbP^n = \calS_0^n$ would be possible if the approximation space of the tensor $\bbB$ is $\bbS^n = [\calS_0^n]^{d \times d}_{\mathrm{S}}$. Then, a different approximation of expressions with spatial derivatives of $\bbB$ has to be applied. For more details, we refer to \cite{barrett_boyaval_2009}.

Similarly to \cite{BGN_2015_navierstokes}, we introduce the following parametric finite element spaces. Let $\Gamma^n\subset\bbR^d$, $d\in\{2,3\}$, be a $(d-1)$-dimensional polyhedral surface, i.e., a union of non-degenerate $(d-1)$-simplices without hanging vertices (see \cite[p.~164]{deckelnick_dziuk_elliott_2005} for $d=3$), approximating the closed surface $\Gamma(t^n)$, $n\in\{0,\ldots,N_T-1\}$. 
More precisely, let $\calT(\Gamma^n)$ be a family of mutually disjoint open $(d-1)$-simplices such that $\Gamma^n = \bigcup_{\sigma\in\calT(\Gamma^n)} \overline{\sigma}$.
Then, for $n\in\{0,\ldots,N_T-1\}$, we define
\begin{align*}
	\mathrm{V}(\Gamma^n) &\coloneqq 
	\{ \chi \in {\color{blue} C(\Gamma^n)} \mid \chi_{|\sigma} \in \calP_1(\sigma) \ \forall\, \sigma\in\calT(\Gamma^n) \} \subset H^1(\Gamma^n),
	\\ 
	\mathbf{V}(\Gamma^n) &\coloneqq [\mathrm{V}(\Gamma^n)]^d \subset H^1(\Gamma^n;\bbR^d).
\end{align*}

We let $\skp{\cdot}{\cdot}$ denote the $L^2$ inner product in $\Omega$. By $\skp{\cdot}{\cdot}_{\calT^n}^h$ we define the mass-lumped inner product in $\Omega$ with respect to the mesh $\calT^n$ as
\begin{align}
    \label{eq:def_mass_lumping}
	\skp{a}{b}_{\calT^n}^h \coloneqq 
	\frac{1}{d+1} \sum_{K\in\calT^n} \abs{K}
	\sum_{k=0}^d \lim_{K\ni \mathbf{p}\to \mathbf{p}^K_k} a(\mathbf{p}) \, b(\mathbf{p}),
\end{align}
where $a,b\in L^\infty(\Omega)$ are piecewise continuous functions with possible jumps across the sides of the elements $K\in\calT^n$.
Here, $\mathbf{p}_k^K$, $k\in\{0,\ldots,d\}$, are the vertices on $K\in\calT^n$.
Similarly, we let $\dualp{\cdot}{\cdot}_{\Gamma^n}$ denote the $L^2$ inner product on the current polyhedral surface $\Gamma^n$. Moreover, we define a mass-lumped inner product $\dualp{\cdot}{\cdot}_{\Gamma^n}^h$ of $\Gamma^n$ as
\begin{align}
    \label{eq:def_mass_lumping_interface}
	\dualp{a}{b}_{\Gamma^n}^h \coloneqq 
	\frac{1}{d} \sum_{\sigma\in\calT(\Gamma^n)} \abs{\sigma}
	\sum_{k=0}^{d-1} \lim_{\sigma \ni \mathbf{q}\to \mathbf{q}^\sigma_k} a(\mathbf{q}) \, b(\mathbf{q}),
\end{align}
where $a,b\in L^\infty(\Gamma^n)$ are piecewise continuous functions with possible jumps across the sides of the interface elements $\sigma \in \calT(\Gamma^n)$, where $\mathbf{q}_k^\sigma$, $k\in\{0,\ldots,d-1\}$, are the vertices of $\sigma \in \calT(\Gamma^n)$.
For continuous functions, the mass-lumped inner products \eqref{eq:def_mass_lumping}, \eqref{eq:def_mass_lumping_interface} can be expressed as 
\begin{alignat*}{2}
    \skp{a}{b}_{\calT^n}^h &= \skp{\mathrm{I}^n_1 [ a b]}{1}  \quad 
    &&\forall\, a,b\in C(\overline\Omega),
    \\
    \skp{a}{b}_{\calT^n}^h &= \dualp{\mathrm{I}^n_1 [ a b]}{1}_{\Gamma^n} \quad 
    &&\forall\, a,b\in C(\Gamma^n).
\end{alignat*}
Moreover, we naturally extend the definitions of \eqref{eq:def_mass_lumping}, \eqref{eq:def_mass_lumping_interface} to vector and matrix valued functions.

Given $\Gamma^n$, we denote the interior and exterior of $\Gamma^n$ by $\Omega_-^n$ and $\Omega_+^n$, respectively, such that $\Gamma^n = \partial\Omega_-^n = \overline {\Omega_-^n} \cap \overline {\Omega_+^n}$. Then, we partition the elements of the bulk triangulation $\calT^n$ into interior, exterior and interfacial elements, i.e.,
\begin{align*}
	\calT_-^n &\coloneqq \{ K \in \calT^n \mid K \subset \Omega_-^n \},
	\qquad
	\calT_+^n \coloneqq \{ K \in \calT^n \mid K \subset \Omega_+^n \},
	\\
	\calT_{\Gamma^n}^n &\coloneqq \{ K \in \calT^n \mid K \cap \Gamma^n \not= \emptyset \}.
\end{align*}
Note that in the case of a fitted bulk mesh it holds that $\calT_{\Gamma^n}^n = \emptyset$.
Moreover, we see that $\calT^n = \calT_-^n \cup \calT_+^n \cup \calT_{\Gamma^n}^n$ is a disjoint partition which can be computed with a practical algorithm from \cite{BGN_2013_stefanproblem}.
Here we assume that $\Gamma^n$ has no self-intersections, and for the numerical experiments in this paper this is always the case. Moreover, we define the piecewise constant unit normal $\pmb\nu^n$ to $\Gamma^n$ such that $\pmb\nu^n$ points into $\Omega_+^n$.

Similarly to \cite{BGN_2013_stefanproblem}, we define the following practical approximation of the mass density $\rho(\cdot,t) =\rho_+ (1-\chi_{|\Omega_-(t)}) + \rho_- \, \chi_{|\Omega_-(t)}$. For $n\in\{0,\ldots,N_T-1\}$, we define the discrete mass density $\rho^n\in\calS_0^n$ elementwise on any $K\in\calT^n$ as
\begin{align}
    \label{eq:def_density}
	(\rho^n)_{|K} = 
	\begin{cases}
		\rho_- & K\in\calT_-^n,\\
		\rho_+ & K\in\calT_+^n,\\
		\frac12(\rho_- + \rho_+) &K\in\calT_{\Gamma^n}^n.
	\end{cases}
\end{align}
We note that in the case of a fitted bulk mesh, \eqref{eq:def_density} reduces to $\rho^n = \rho_+ (1-\chi_{|\Omega_-^n}) + \rho_- \, \chi_{|\Omega_-^n} \in \calS_0^n$.
We define $\mu^n, \lambda^n \in\calS_0^n$ analogously to \eqref{eq:def_density} as the approximations of the shear viscosity $\mu(\cdot,t)$ and the relaxation parameter $\lambda(\cdot,t)$, respectively. For other possible definitions, we refer to, e.g., \cite{BGN_2015_navierstokes}.

\subsection{A stable fully discrete approximation}
Let the closed polyhedral hypersurface $\Gamma^0$ be an approximation of $\Gamma(0)$ and let $\bu^0 \in \bbU^0$ and $\bbB^0 \in \bbS^0$ positive definite be approximations of $\bu_0$ and $\bbB_0$, respectively. Then, for any $n\in\{0,\ldots,N_T-1\}$, find $(\bu^{n+1}, p^{n+1}, \bbB^{n+1}, \bX^{n+1}, \kappa^{n+1}) \in \bbU^n \times \bbP^n \times \bbS^n \times \mathbf{V}(\Gamma^n) \times \mathrm{V}(\Gamma^n)$ with $\bbB^{n+1}$ positive definite such that 
\begin{subequations}
\label{eq:system_FE}
\begin{align}
    \label{eq:u_FE} \nonumber
    0 &= \frac{1}{2\Delta t} \skp{\rho^n \bu^{n+1} - (\mathrm{I}^n_0\rho^{n-1}) \bI^n_2 \bu^n }{\bw}
    + \frac{1}{2\Delta t} \skp{\mathrm{I}^n_0\rho^{n-1} (\bu^{n+1} - \bI^n_2 \bu^n) }{\bw}
    \\ \nonumber
    &\quad
    + \frac12 \skp{\rho^n }{ [(\bI^n_2\bu^n \cdot\nabla) \bu^{n+1} ] \cdot \bw}
    - \frac12 \skp{\rho^n }{\bu^{n+1} \cdot [(\bI^n_2\bu^n  \cdot\nabla) \bw]}
    \\ \nonumber
    &\quad 
    + \skp{2 \mu^n \D(\bu^{n+1})}{\D(\bw)}
    - \skp{p^{n+1}}{\Div \bw}
    + G \skp{\bbB^{n+1} - \bbI}{\nabla\bw}
    \\
    &\quad
    - \gamma \dualp{\kappa^{n+1} \pmb\nu^n}{ \bw}_{\Gamma^n}
    - \skp{\rho^n \mathbf{f}_1^{n+1} + \mathbf{f}_2^{n+1}}{\bw} ,
    \\
    \label{eq:div_FE}
    0 &= \skp{\Div \bu^{n+1}}{q},
    \\
    \label{eq:B_FE} \nonumber
    0 &= 
    \frac{1}{\Delta t} \skp{\bbB^{n+1} - \bbB^n}{\bbG}_{\calT^n}^h
    - \sum_{i,j=1}^d \skp{[\bu^{n+1}]_i \, \mathbf{\Lambda}_{i,j}(\bbB^{n+1})}{\partial_{x_j} \bbG}
    \\
    &\quad 
     + \skp{\tfrac{1}{\lambda^n} (\bbB^{n+1} - \bbI) }{\bbG}_{\calT^n}^h 
    - 2 \skp{\nabla\bu^{n+1}}{ \bbI_1^n \left[ \bbG \bbB^{n+1} \right]}
    + \alpha \skp{\nabla\bbB^{n+1}}{\nabla\bbG},
    \\
    \label{eq:kappa_FE}
    0 &= \dualp{\kappa^{n+1} \pmb\nu^n}{\pmb\zeta}_{\Gamma^n}^h
    + \dualp{\nabla_s \bX^{n+1}}{\nabla_s \pmb\zeta}_{\Gamma^n},
    \\
    \label{eq:X_FE}
    0 &= \frac{1}{\Delta t} \dualp{(\bX^{n+1} - \mathbf{id}) \cdot \pmb\nu^n}{\chi}_{\Gamma^n}^h
    - \dualp{\bu^{n+1} \cdot \pmb\nu^n}{\chi}_{\Gamma^n},
\end{align}
\end{subequations}
for any $(\bw, q, \bbG, \pmb\zeta, \chi) \in \bbU^n \times \bbP^n \times \bbS^n \times \mathbf{V}(\Gamma^n) \times \mathrm{V}(\Gamma^n)$,
and set $\Gamma^{n+1} \coloneqq \bX^{n+1}(\Gamma^n)$.
Here we have defined $\mathbf{f}_i^{n+1} \coloneqq \bI_2^n \, \mathbf{f}_i(t^{n+1}) \in \bbU^n$, $i\in\{1,2\}$, where here and throughout we assume that $\mathbf{f}_i \in C([0,T]; C(\overline\Omega;\bbR^d))$, $i\in\{1,2\}$. 
Moreover, we define $\calT^{-1} \coloneqq \calT^0$ and $\rho^{-1} \coloneqq \rho^0$.

For $\bbB^{n+1}\in\bbS^n$ positive definite, the nonlinear quantity $\mathbf{\Lambda}_{i,j}(\bbB^{n+1}) \in [\calS_0^n]^{d\times d}$, $i,j\in\{1,\ldots,d\}$, is constructed similarly to \cite{barrett_boyaval_2009} such that the identity
\begin{align}
\label{eq:Lambda_kettenregel}
    \sum\limits_{j=1}^d 
    {\mathbf\Lambda}_{i,j} (\bbB^{n+1}) : \partial_{x_j} \bbI_1^n \big[ (\bbB^{n+1})^{-1} \big]
    = \partial_{x_i} \mathrm{I}_1^n\left[ \trace \ln (\bbB^{n+1})^{-1} \right] \, ,
\end{align}
holds on any simplex $K \in \calT^n$ and for each $i\in\{1,\ldots,d\}$. 
The identity \eqref{eq:Lambda_kettenregel} will be useful in the energy estimates on the fully discrete level, as we then have with integration by parts and \eqref{eq:div_FE} that
\begin{align}
    \label{eq:Lmbda_kettenregel2}
    \sum_{i,j=1}^d \skp{[\bu^{n+1}]_i \, \mathbf{\Lambda}_{i,j}(\bbB^{n+1})}{\partial_{x_j} \bbI_1^n\big[ (\bbB^{n+1})^{-1} \big]}
    = \skp{\bu^{n+1}}{\nabla \mathrm{I}_1^n \left[ \trace \ln (\bbB^{n+1})^{-1} \right] }
    = 0.
\end{align}
Without mesh adaption and if $\bu^n$ is divergence-free with respect to the discrete pressure space, meaning $\skp{\Div \bu^n}{q} = 0$ for all $q \in \bbP^{n-1}=\bbP^n$, the expression $[\bu^{n+1}]_i \mathbf{\Lambda}_{i,j}(\bbB^{n+1})$ in \eqref{eq:B_FE} can be substituted by $[\bu^n]_i \mathbf{\Lambda}_{i,j}(\bbB^{n+1})$. Thus, \eqref{eq:Lmbda_kettenregel2} remains correct when $\bu^{n+1}$ is replaced by $\bu^n$.

For the sake of completeness, we state the definition of $\mathbf{\Lambda}_{i,j}(\bbB^{n+1})$.
Let $\hat K$ denote the standard open reference simplex in $\bbR^d$ with vertices $\hat{\mathbf{p}}_0, \ldots, \hat{\mathbf{p}}_d$, where $\hat{\mathbf{p}}_0$ is the origin and $\hat{\mathbf{p}}_i$, $i\in\{1,\ldots,d\}$, is the $i$-th standard basis vector in the $\bbR^d$. Given a non-degenerate simplex $K\in\calT^n$ with vertices $\mathbf{p}_0^K, \ldots, \mathbf{p}_d^K$, there exists an invertible matrix $\calA_K \in \bbR^{d\times d}$ such that the affine transformation
\begin{align*}
	\calB_K \colon\, \hat K \to K, \quad \hat{\mathbf{x}}\mapsto \mathbf{p}_0^K + \calA_K \hat{\mathbf{x}}
\end{align*}
maps vertex $\hat{\mathbf{p}}_i$ to vertex $\mathbf{p}_i^K$, $i\in\{0,\ldots,d\}$.
Similarly to \cite{barrett_boyaval_2009}, we introduce the notation
\begin{align*}
	\hat\bbB(\hat{\mathbf{x}}) \coloneqq \bbB(\calB_K(\hat{\mathbf{x}})),
	\quad\quad 
	(\hat\calI_h \hat\bbB)(\hat{\mathbf{x}}) \coloneqq (\calI_h \bbB)(\calB_K(\hat{\mathbf{x}})),
	\quad\quad \forall \,  \hat{\mathbf{x}}\in\hat K, \, \bbB\in C(\overline K; \bbR^{d\times d}_{\mathrm{S}}),
\end{align*}
and we define $\bbB_j^K \coloneqq \bbB(\mathbf{p}_j^K)$, $j\in\{0,\ldots,d\}$, for any $\bbB\in\bbS^n$ and $K\in\calT^n$, where $\mathbf{p}_0^K, \ldots, \mathbf{p}_d^K$ denote the vertices of the simplex $K$. 
On the reference element $\hat K$, we define
\begin{align*}
	\hat{\mathbf\Lambda} _i(\hat\bbB) 
	\coloneqq 
	\begin{cases}
		\bbB_i^K & \text{ if } \bbB_i^K = \bbB_0^K,
		\\
		\bbB_i^K + \lambda_i(\hat\bbB)
		\big( \bbB_0^K - \bbB_i^K \big)
		& \text{ if } \bbB_i^K \not= \bbB_0^K,
	\end{cases}
\end{align*}
where $\lambda_i(\hat\bbB) \in \bbR$, $i\in\{1,\ldots,d\}$, are defined as
\begin{align*}
	\lambda_i(\hat\bbB) 
	&\coloneqq
	\frac{  \trace \ln (\bbB_i^K)^{-1}
		- \trace \ln (\bbB_0^K)^{-1} 
        + \bbB_i^K : 
		\big( - (\bbB_i^K)^{-1} 
		+ (\bbB_0^K)^{-1} \big)}
	{ \big(\bbB_i^K - \bbB_0^K \big)  
		: \big( - (\bbB_i^K)^{-1} 
		+ (\bbB_0^K)^{-1}  \big)  }\,.
\end{align*}
With the help of $\hat{\mathbf\Lambda} _i(\hat\bbB)$, we define ${\mathbf\Lambda} _{i,j}(\bbB)_{|K} \in \bbR^{d\times d}_{\mathrm{S}}$, $i,j\in\{1,\ldots,d\}$, on $K\in\calT^n$ by
\begin{align*}
	{\mathbf\Lambda} _{i,j}(\bbB) _{|K}
	\coloneqq
	\sum\limits_{m=1}^d [(\calA_K^\top)^{-1}]_{i,m} \,
	\hat{\mathbf\Lambda} _{m}(\hat\bbB) \, [\calA_K^\top]_{m,j} \, \in\bbR^{d\times d}_{\mathrm{S}},
\end{align*}
on any $K \in \calT^n$ and for each $i,j\in\{1,\ldots,d\}$. 
By construction, we see that the quantity $\{ [{\mathbf\Lambda} _{i,j}(\bbB)]_{m,n} \}_{i,j,m,n\in\{1,\ldots,d\}}$ is a fourth-order tensor that is constant on each element $K\in\calT^n$ for $\bbB\in\bbS^n$ positive definite.
Besides, ${\mathbf\Lambda} _{i,j}(\bbB)_{|K} \in \bbR^{d\times d}_{\mathrm{S}}$ continuously depends on $\bbB_{|K}$ for any $K\in\calT^n$ and for each $i,j\in\{1,\ldots,d\}$. 
We note that it follows from the construction that \eqref{eq:Lambda_kettenregel} is fulfilled and it holds $\lambda_i(\hat\bbB) \in [0,1]$, where the proof is similar to \cite[Lem.~3.6]{sieber_2020}, using that $-\ln(\cdot)$ is strictly convex as a function on the set of positive definite matrices. As the family of triangulations is shape regular, one can derive an error estimate for $\mathbf\Lambda_{i,j}(\bbB)$ to conclude that $\mathbf\Lambda_{i,j}(\bbB) \approx \delta_{i,j} \, \bbB$ for all $i,j\in\{1,\ldots,d\}$ and $\bbB\in\bbS^n$ positive definite, where $\delta_{i,j}$ denotes the Kronecker symbol. For similar approaches, we refer to \cite{barrett_boyaval_2009, barrett_2018_fene-p, GT_2023_DCDS}.

For the mathematical analysis of \eqref{eq:system_FE}, it is convenient to introduce a reduced version where the discrete pressure is eliminated. This can also be important if $(\bbU^n,\hat\bbP^n)$ does not satisfy the inf-sup condition \eqref{eq:LBB}.
Let
\begin{align*}
    \bbU^n_0 \coloneqq \{ \bu\in\bbU^n \mid \skp{\Div \bu}{q} = 0 \ \forall \, q\in \bbP^n \}.
\end{align*}
Then the discrete pressure can be eliminated from \eqref{eq:system_FE} to produce the following reduced variant.
For any $n\in\{0,\ldots,N_T-1\}$, find $(\bu^{n+1}, \bbB^{n+1}, \bX^{n+1}, \kappa^{n+1}) \in \bbU_0^n \times \bbS^n \times \mathbf{V}(\Gamma^n) \times \mathrm{V}(\Gamma^n)$ with $\bbB^{n+1}$ positive definite such that 
\begin{subequations}
\label{eq:system_FE2}
\begin{align}
    \label{eq:u_FE2} \nonumber
    0 &= \frac{1}{2\Delta t} \skp{\rho^n \bu^{n+1} - (\mathrm{I}^n_0\rho^{n-1}) \bI^n_2 \bu^n }{\bw}
    + \frac{1}{2\Delta t} \skp{\mathrm{I}^n_0\rho^{n-1} (\bu^{n+1} - \bI^n_2 \bu^n) }{\bw}
    \\ \nonumber
    &\quad
    + \frac12 \skp{\rho^n }{ [(\bI^n_2\bu^n \cdot\nabla) \bu^{n+1} ] \cdot \bw}
    - \frac12 \skp{\rho^n }{\bu^{n+1} \cdot [(\bI^n_2\bu^n  \cdot\nabla) \bw]}
    \\ \nonumber
    &\quad 
    + \skp{2 \mu^n \D(\bu^{n+1})}{\D(\bw)}
    + G \skp{\bbB^{n+1} - \bbI}{\nabla\bw}
    \\
    &\quad 
    - \gamma \dualp{\kappa^{n+1} \pmb\nu^n}{ \bw}_{\Gamma^n}
    - \skp{\rho^n \mathbf{f}_1^{n+1} + \mathbf{f}_2^{n+1}}{\bw} ,
    \\
    \label{eq:B_FE2} \nonumber
    0 &= 
    \frac{1}{\Delta t} \skp{\bbB^{n+1} - \bbB^n}{\bbG}_{\calT^n}^h
    - \sum_{i,j=1}^d \skp{[\bu^{n+1}]_i \, \mathbf{\Lambda}_{i,j}(\bbB^{n+1})}{\partial_{x_j} \bbG}
    \\
    &\quad 
     + \skp{\tfrac{1}{\lambda^n} (\bbB^{n+1} - \bbI) }{\bbG}_{\calT^n}^h 
    - 2 \skp{\nabla\bu^{n+1}}{ \bbI_1^n \left[ \bbG \bbB^{n+1} \right]}
    + \alpha \skp{\nabla\bbB^{n+1}}{\nabla\bbG},
    \\
    \label{eq:kappa_FE2}
    0 &= \dualp{\kappa^{n+1} \pmb\nu^n}{\pmb\zeta}_{\Gamma^n}^h
    + \dualp{\nabla_s \bX^{n+1}}{\nabla_s \pmb\zeta}_{\Gamma^n},
    \\
    \label{eq:X_FE2}
    0 &= \frac{1}{\Delta t} \dualp{(\bX^{n+1} - \mathbf{id}) \cdot \pmb\nu^n}{\chi}_{\Gamma^n}^h
    - \dualp{\bu^{n+1} \cdot \pmb\nu^n}{\chi}_{\Gamma^n},
\end{align}
\end{subequations}
for any $(\bw, \bbG, \pmb\zeta, \chi) \in \bbU_0^n \times \bbS^n \times \mathbf{V}(\Gamma^n) \times \mathrm{V}(\Gamma^n)$,
and set $\Gamma^{n+1} = \bX^{n+1}(\Gamma^n)$.


%
%
%

\begin{remark}
The discrete systems \eqref{eq:system_FE} and \eqref{eq:system_FE2} are nonlinear due to the coupling of $[\bu^{n+1}]_i$ with the nonlinear quantity $\mathbf{\Lambda}_{i,j}(\bbB^{n+1}) \in [\calS_0^n]^{d\times d}$, $i,j\in\{1,\ldots,d\}$, and the coupled term $\skp{\nabla\bu^{n+1}}{ \bbI_1^n \left[ \bbG \bbB^{n+1} \right]}$ in the viscoelastic equation. 
Notably, without viscoelastic effects, our approximations \eqref{eq:system_FE} and \eqref{eq:system_FE2} coincide with the scheme, and its pressure-free variant, from \cite{BGN_2015_navierstokes} for two-phase Navier--Stokes flow, which are both linear. This linearity is attributed to the explicit treatment of the interface $\Gamma^n$ and the approximation $\pmb\nu^n$ of the unit normal, as well as the lagging in the approximation $\rho^n$ of the density.

Furthermore, we aim to compute a positive definite $\bbB^{n+1} \in \bbS^n$, which is essential for our approximations to satisfy a discrete energy inequality based on \eqref{eq:stability}, where the positive definiteness of the tensor $\bbB$ is required. The identity \eqref{eq:Lambda_kettenregel}, involving the nonlinear quantity $\mathbf{\Lambda}_{i,j}(\bbB^{n+1}) \in [\calS_0^n]^{d\times d}$ for $i,j\in\{1,\ldots,d\}$, is crucial to derive this discrete energy inequality. To construct $\mathbf{\Lambda}_{i,j}(\bbB^{n+1}) \in [\calS_0^n]^{d\times d}$, $i,j\in\{1,\ldots,d\}$, the discrete tensor $\bbB^{n+1} \in \bbS^n$ must be positive definite. It is also worth noting that the test functions $\bbG \in \bbS^n$ in \eqref{eq:B_FE} and \eqref{eq:B_FE2}, respectively, only need to be symmetric to ensure the symmetry of $\bbB^{n+1}$.
\end{remark}

Before we prove the existence and stability for a solution to \eqref{eq:system_FE}, we introduce the notion of a vertex normal on $\Gamma^n$. This definition is combined with an assumption that is needed in the existence proof.

\begin{assumptions} \label{assumptions}
For $n\in\{0,\ldots,N_T-1\}$, we assume that $\calH^{d-1}(\sigma)>0$ for all interface elements $\sigma\in\calT(\Gamma^n)$, and that $\Gamma^n\subset\overline\Omega$. For $k\in\{1,\ldots,\calK^n_{\Gamma^n}\}$, let $\mathbf{q}_k^n$ denote the vertices of $\Gamma^n$. Define $\Xi_k^n \coloneqq \{\sigma \mid \mathbf{q}_k^n \in \overline{\sigma} \}$, and set
\begin{align*}
    \Theta_k^n \coloneqq \bigcup_{\sigma \in \Xi_k^n} \overline{\sigma},
    \qquad 
    \pmb\omega_k^n \coloneqq \frac{1}{\calH^{d-1}(\Theta_k^n)} \, \sum_{\sigma \in \Xi_k^n} \calH^{d-1}(\sigma) \, \pmb\nu^n_{| \sigma}.
\end{align*}
We further assume that 
\begin{align}
    \label{eq:assumptions_vertex_normal}
    \mathrm{dim} \; \mathrm{span}\{ \, \pmb\omega_1^n, \ldots, \, \pmb\omega_{\calK_{\Gamma^n}^n}^n \} = d,
\end{align}
for any $n\in\{0,\ldots,N_T-1\}$.

Moreover, in the case $\alpha>0$ we also require that the bulk mesh $\calT^n$ consists of only non-obtuse simplices, i.e., all dihedral angles of any simplex in $\calT^n$ are less than or equal to $\frac\pi2$.
\end{assumptions}

We stress that \eqref{eq:assumptions_vertex_normal} is a very mild assumption and is only rarely violated, usually occurring only in exceptional cases, such as surfaces with self-intersections.
More details can be found in \cite{BGN_2007}.
The assumption \eqref{eq:assumptions_vertex_normal} is used in the proof of Lemma \ref{lemma:eq:existence_FE_delta} to conclude the unique solvability of the linear subsystem \eqref{eq:kappa_FE}--\eqref{eq:X_FE} if the velocity $\bu^{n+1}$ is already given. This idea is related to the Schur complement approach from \cite{BGN_2015_navierstokes} which is used in Section \ref{sec:numerics} for the numerical computations. 
However, constructing a non-obtuse mesh $\calT^n$ in the case $d=3$ for a general polytope $\Omega$ is not straightforward. 
In the case $\alpha>0$, this assumption is needed for a discrete analogue of the inequality \eqref{eq:nabla_ln_formal}, as seen in \eqref{eq:nabla_ln}.
Note that in the case $\alpha=0$ the assumption of non-obtuse simplices is not needed at all. 

We now present our main result, which is valid with stress diffusion ($\alpha>0$) and without stress diffusion ($\alpha=0$). 

\begin{theorem} \label{thm:stability_FE}
Let Assumptions \ref{assumptions} hold true.
\begin{enumerate}
\item[(i)] There exists a solution $(\bu^{n+1}, \bbB^{n+1}, \bX^{n+1}, \kappa^{n+1}) \in \bbU_0^n \times \bbS^n \times \mathbf{V}(\Gamma^n) \times \mathrm{V}(\Gamma^n)$ to \eqref{eq:system_FE2} with $\bbB^{n+1}$ positive definite. 

\item[(ii)] If $(\bu^{n+1}, p^{n+1}, \bbB^{n+1}, \bX^{n+1}, \kappa^{n+1}) \in \bbU^n \times \bbP^n \times \bbS^n \times \mathbf{V}(\Gamma^n) \times \mathrm{V}(\Gamma^n)$ with $\bbB^{n+1}$ positive definite solves \eqref{eq:system_FE}, then $(\bu^{n+1}, \bbB^{n+1}, \bX^{n+1}, \kappa^{n+1})$ is a solution to \eqref{eq:system_FE2}.

\item[(iii)] If the pair $(\bbU^n,\bbP^n)$ satisfies the inf-sup condition \eqref{eq:LBB}, then there exists a solution $(\bu^{n+1}, p^{n+1}, \bbB^{n+1}, \bX^{n+1}, \kappa^{n+1}) \in \bbU^n \times \bbP^n \times \bbS^n \times \mathbf{V}(\Gamma^n) \times \mathrm{V}(\Gamma^n)$ to \eqref{eq:system_FE} with $\bbB^{n+1}$ positive definite.

\item[(iv)] All solutions to \eqref{eq:system_FE2} with $\bbB^{n+1}$ positive definite satisfy 
\begin{align}
    \label{eq:stability_FE} \nonumber
    & \frac12 \norm{\sqrt{\rho^n} \bu^{n+1}}_{L^2}^2 
    + \skp{W(\bbB^{n+1})}{1}_{\calT^n}^h 
    + \gamma \calH^{d-1}(\Gamma^{n+1})
    \\ \nonumber
    &\quad
    + \frac12 \norm{\sqrt{\mathrm{I}_0^n \rho^{n-1}} (\bu^{n+1} - \bI_2^n \bu^n)}_{L^2}^2
    + 2 \Delta t \norm{\sqrt{\mu^n} \D(\bu^{n+1})}_{L^2}^2 
    \\ \nonumber
    &\quad
    + \Delta t \skp{\tfrac{G}{2\lambda^n} }{\trace( \bbB^{n+1} + [\bbB^{n+1}]^{-1} - 2\bbI ) }_{\calT^n}^h
    + \Delta t \frac{\alpha G}{2d} \norm{\nabla \mathrm{I}_1^n \trace\ln\bbB^{n+1}}_{L^2}^2
    \\
    &\leq
    \frac12 \norm{\sqrt{\mathrm{I}_0^n \rho^{n-1}} \bI_2^n \bu^n}_{L^2}^2
    + \skp{W(\bbB^n)}{1}_{\calT^n}^h + \gamma \calH^{d-1}(\Gamma^n)
    + \Delta t \skp{\rho^n \mathbf{f}_1^{n+1} + \mathbf{f}_2^{n+1}}{\bu^{n+1}} ,
\end{align}
where $W(\bbB) = \tfrac{G}{2} \trace(\bbB - \ln\bbB - \bbI)$ is the elastic energy density. 
\end{enumerate}
\end{theorem}

\begin{proof}
The proofs of (i) and (iv) are presented in the next subsection. The statement in (ii) follows directly from the definition of the discrete function space $\bbU_0^n$, while (iii) follows from (i) and the inf-sup condition \eqref{eq:LBB}.
\end{proof}

The above theorem allows us to prove unconditional stability.

\begin{theorem} \label{thm:stability_FE_sum}
Let Assumptions \ref{assumptions} hold true.
In addition, we assume that
\begin{align*}
    & \int_\Omega \mathrm{I}_0^n \rho^{n-1} \abs{\bI_2^n \bu^n}^2 \dx
    + \int_\Omega \mathrm{I}_1^n  W(\bbB^n) \dx 
    \leq \int_\Omega \rho^{n-1} \abs{\bu^n}^2 \dx
    + \int_\Omega \mathrm{I}_1^{n-1} W(\bbB^n) \dx, 
\end{align*}
which, e.g., is fulfilled if no bulk mesh coarsening in time is performed.
Then, it holds for any $n\in\{0,\ldots,N_T-1\}$, that
\begin{align}
    \label{eq:stability_FE_sum} \nonumber
    & \frac12 \norm{\sqrt{\rho^n} \bu^{n+1}}_{L^2}^2 
    + \skp{W(\bbB^{n+1})}{1}_{\calT^n}^h
    + \gamma \calH^{d-1}(\Gamma^{n+1})
    \\ \nonumber
    &\quad
    + \frac12 \sum_{k=0}^n \norm{\sqrt{\rho^{k-1}} (\bu^{k+1} - \bI_2^k \bu^k)}_{L^2}^2
    + 2 \Delta t \sum_{k=0}^n \norm{\sqrt{\mu^k} \D(\bu^{k+1})}_{L^2}^2 
    \\ \nonumber
    &\quad
    + \Delta t \sum_{k=0}^n \skp{\tfrac{G}{2\lambda^k} }{\trace( \bbB^{k+1} + [\bbB^{k+1}]^{-1} - 2\bbI )}_{\calT^k}^h
    + \Delta t \frac{\alpha G}{2d} \sum_{k=0}^n \norm{\nabla \mathrm{I}_1^k \trace\ln\bbB^{k+1}}_{L^2}^2
    \\
    &\leq
    \frac12 \norm{\sqrt{\rho^{-1}} \bu^0}_{L^2}^2
    + \skp{W(\bbB^0)}{1}_{\calT^{-1}}^h + \gamma \calH^{d-1}(\Gamma^0)
    + \Delta t \sum_{k=0}^n \skp{\rho^k \mathbf{f}_1^{k+1} + \mathbf{f}_2^{k+1}}{\bu^{k+1}}.
\end{align}
\end{theorem}


\cblue{We note that, as is common with parametric finite element methods, there is no guarantee that $\Gamma^{n+1}$ remains free of self-intersections. However, no self-intersections were observed in any of the numerical experiments presented in this paper.}

\subsection{Proof of the stability and existence results} \label{sec:stability_proof}
We now present the proof of Theorem \ref{thm:stability_FE}, which is divided into the following steps.
First, we introduce cut-offs in $\bbB$ in certain terms in the system \eqref{eq:system_FE2}, such that $\bbB^{n+1}$ does not have to be positive definite. After that, we derive an energy estimate for the regularised problem, which holds uniformly in the regularisation parameter and in the discretization parameters. This allows us to control the negative eigenvalues of $\bbB^{n+1}$ uniformly in the regularisation parameter. Then, using the energy estimates and a fixed-point argument, we show the existence of at least one solution to the regularised system. Sending the regularisation parameter to zero and extracting converging subsequences of solutions finally proves Theorem \ref{thm:stability_FE}, i.e., the limit functions form a solution to \eqref{eq:system_FE2} with $\bbB^{n+1}$ positive definite. We reconstruct the pressure in the case when the inf-sup stability condition \eqref{eq:LBB} is met.

Given $\delta>0$, we define a regularisation of the logarithmic function by
\begin{align*}
	\bbR \ni s \mapsto g_\delta(s) \coloneqq 
	\begin{cases}
    \frac{s}{\delta} + \ln(\delta) - 1, 
    & s < \delta,
    \\
    \ln(s), 
    & s \geq \delta.
	\end{cases}
\end{align*}
We further define $\beta_\delta(s) \coloneqq g_\delta'(s)^{-1} = \max\{s,\delta\}$ for any $s\in\bbR$. In the following, we extend these functions to matrix valued functions by applying them to the spectrum. In particular, let $\bbB \in\bbR^{d\times d}_S$ be a symmetric matrix such that $\bbB = \mathbb{U}^\top \bbD \mathbb{U}$ with an orthogonal matrix $\mathbb{U}$ and a diagonal matrix $\bbD = \mathrm{diag}(D_{11}, \ldots, D_{dd})$. Then, for any scalar valued function $F:\bbR\to\bbR$, we define $F(\bbB)\coloneqq \mathbb{U}^\top F(\D) \mathbb{U}$, where $F(\D) \coloneqq \mathrm{diag}(F(D_{11}), \ldots, F(D_{dd}))$.

We recall the following result from \cite[Lem.~2.1]{barrett_boyaval_2009}.
\begin{lemma}
For all $\bbB,\bbG\in\bbR^{d\times d}_{\mathrm{S}}$ and for any $\delta \in (0,1)$, it holds
\begin{subequations}
	\begin{align}
		\label{eq:lemma_reg1b}
		\trace\big( \beta_\delta(\bbB) 
		+ \beta_\delta(\bbB)^{-1} - 2\I \big) 
		&\geq 0,
		\\
		\label{eq:lemma_reg1d}
		\big( \bbB - \beta_\delta(\bbB)\big) 
		: \big( \I - g_\delta'(\bbB) \big) 
		&\geq 0,
		\\
		\label{eq:lemma_reg1e}
		(\bbG-\bbB) : g_\delta'(\bbB) 
		&\geq \trace\big( g_\delta(\bbG) - g_\delta(\bbB) \big).
	\end{align}
	In addition, if $\delta\in(0,\frac{1}{2}]$, it holds
	\begin{align}
		\label{eq:lemma_reg1g}
		\trace\big( \bbB - g_\delta(\bbB) \big)  
		&\geq
		\begin{cases}
			\frac{1}{2} \abs{\bbB},
			\\
			\frac{1}{2\delta} \abs{ [\bbB]_{-} },
		\end{cases}
		\\
		\label{eq:lemma_reg1h}
		\bbB : \big( \I - g_\delta'(\bbB) \big) 
		&\geq \frac{1}{2}\abs{\bbB} - d,
	\end{align}
\end{subequations}
where $[\cdot]_{-}$ denotes the negative part function defined by $[s]_{-} \coloneqq \min\{s,0\}$ $\forall \,  s\in\bbR$.
\end{lemma}

Now, we introduce the following discrete problem with cut-offs in $\bbB$.
Let $\delta>0$.  
For any $n\in\{0,\ldots,N_T-1\}$, let $\Gamma^n$ be a closed polyhedral hypersurface and let $\bu^n \in C(\overline\Omega;\bbR^d)$ and $\bbB^n \in  C(\overline\Omega;\bbR^{d\times d}_\mathrm{S})$ be given. Then find $(\bu_\delta^{n+1}, \bbB_\delta^{n+1}, \bX_\delta^{n+1}, \kappa_\delta^{n+1}) \in \bbU_0^n \times \bbS^n \times \mathbf{V}(\Gamma^n) \times \mathrm{V}(\Gamma^n)$ such that 
\begin{subequations}
\label{eq:system_FE_delta}
\begin{align}
    \label{eq:u_FE_delta} \nonumber
    0 &= \frac{1}{2\Delta t} \skp{\rho^n \bu_\delta^{n+1} - (\mathrm{I}^n_0\rho^{n-1}) \bI^n_2 \bu^n }{\bw}
    + \frac{1}{2\Delta t} \skp{\mathrm{I}^n_0\rho^{n-1} (\bu_\delta^{n+1} - \bI^n_2 \bu^n) }{\bw}
    \\ \nonumber
    &\quad
    + \frac12 \skp{\rho^n }{ [(\bI^n_2\bu^n \cdot\nabla) \bu_\delta^{n+1} ] \cdot \bw}
    - \frac12 \skp{\rho^n }{\bu_\delta^{n+1} \cdot [(\bI^n_2\bu^n  \cdot\nabla) \bw]}
    \\ \nonumber
    &\quad 
    + \skp{2 \mu^n \D(\bu_\delta^{n+1})}{\D(\bw)}
    + G \skp{ \bbI_1^n \beta_\delta(\bbB_\delta^{n+1})  - \bbI }{\nabla\bw}
    \\
    &\quad
    - \gamma \dualp{\kappa_\delta^{n+1} \pmb\nu^n}{ \bw}_{\Gamma^n}
    - \skp{\rho^n \mathbf{f}_1^{n+1} + \mathbf{f}_2^{n+1}}{\bw} ,
    \\ 
    \label{eq:B_FE_delta} \nonumber
    0 &= 
    \frac{1}{\Delta t} \skp{\bbB_\delta^{n+1} - \bbB^n}{\bbG}_{\calT^n}^h
    - \sum_{i,j=1}^d \skp{[\bu_\delta^{n+1}]_i \,  \mathbf{\Lambda}_{i,j}^\delta(\bbB_\delta^{n+1}) }{\partial_{x_j} \bbG}
    \\
    &\quad 
     + \skp{\tfrac{1}{\lambda^n} (\bbB_\delta^{n+1} - \bbI) }{\bbG}_{\calT^n}^h 
    - 2 \skp{\nabla\bu_\delta^{n+1}}{ \bbI_1^n \big[ \bbG  \beta_\delta(\bbB_\delta^{n+1}) \big]}
    + \alpha \skp{\nabla\bbB_\delta^{n+1}}{\nabla\bbG},
    \\
    \label{eq:kappa_FE_delta}
    0 &= \dualp{\kappa_\delta^{n+1} \pmb\nu^n}{\pmb\zeta}_{\Gamma^n}^h
    + \dualp{\nabla_s \bX_\delta^{n+1}}{\nabla_s \pmb\zeta}_{\Gamma^n},
    \\
    \label{eq:X_FE_delta}
    0 &= \frac{1}{\Delta t} \dualp{(\bX_\delta^{n+1} - \mathbf{id}) \cdot \pmb\nu^n}{\chi}_{\Gamma^n}^h
    - \dualp{\bu_\delta^{n+1} \cdot \pmb\nu^n}{\chi}_{\Gamma^n},
\end{align}
\end{subequations}
for any $(\bw, \bbG, \pmb\zeta, \chi) \in \bbU_0^n \times \bbS^n \times \mathbf{V}(\Gamma^n) \times \mathrm{V}(\Gamma^n)$,
and set $\Gamma_\delta^{n+1} = \bX_\delta^{n+1}(\Gamma^n)$.

For $\bbB_\delta^{n+1} \in \bbS^n$ not necessarily positive definite, the nonlinear quantity $\mathbf{\Lambda}_{i,j}^\delta(\bbB_\delta^{n+1}) \in [\calS_0^n]^{d\times d}$ is constructed similarly to the unregularised version $\mathbf{\Lambda}_{i,j}$ such that the identity
\begin{align}
	\label{eq:Lambda_delta_kettenregel}
	\sum\limits_{j=1}^d 
	{\mathbf\Lambda} _{i,j}^\delta (\bbB_\delta^{n+1}) : \partial_{x_j} \mathbb{I}_1^n \big[\beta_\delta(\bbB_\delta^{n+1})^{-1} \big]
	= \partial_{x_i} \mathrm{I}_1^n \left[ 
    \trace \ln (\beta_\delta(\bbB_\delta^{n+1})^{-1}) \right] \, 
\end{align}
holds on any simplex $K\in\calT^n$ and for each $i\in\{1,\ldots,d\}$.
In particular, on the reference element $\hat K$, we define
\begin{align*}
	\hat{\mathbf\Lambda} _i^\delta(\hat\bbB) 
	\coloneqq 
	\begin{cases}
		\beta_\delta(\bbB_i^K) & \text{ if } \beta_\delta(\bbB_i^K) = \beta_\delta(\bbB_0^K),
		\\
		\beta_\delta(\bbB_i^K) + \lambda_i^\delta(\hat\bbB)
		\big( \beta_\delta(\bbB_0^K) - \beta_\delta(\bbB_i^K) \big)
		& \text{ if } \beta_\delta(\bbB_i^K) \not= \beta_\delta(\bbB_0^K),
	\end{cases}
\end{align*}
where $\lambda_i^\delta(\hat\bbB) \in \bbR$, $i\in\{1,\ldots,d\}$, are defined as
\begin{align*}
	\lambda_i^\delta(\hat\bbB) 
	&\coloneqq
	\frac{  \trace \ln (\beta_\delta(\bbB_i^K)^{-1})
		- \trace \ln (\beta_\delta(\bbB_0^K)^{-1}) }
	{ \big(\beta_\delta(\bbB_i^K) - \beta_\delta(\bbB_0^K) \big)  
		: \big( - \beta_\delta(\bbB_i^K)^{-1} 
		+ \beta_\delta(\bbB_0^K)^{-1}  \big)  }
	\\
	&\quad + 
	\frac{\beta_\delta(\bbB_i^K) : 
		\big( - \beta_\delta(\bbB_i^K)^{-1} 
		+ \beta_\delta(\bbB_0^K)^{-1} \big)}
	{ \big(\beta_\delta(\bbB_i^K) - \beta_\delta(\bbB_0^K) \big)  
		: \big( - \beta_\delta(\bbB_i^K)^{-1} 
		+ \beta_\delta(\bbB_0^K)^{-1} \big) } \,.
\end{align*}
Then, we define ${\mathbf\Lambda} _{i,j}^\delta(\bbB)_{|K} \in \bbR^{d\times d}_{\mathrm{S}}$, $i,j\in\{1,\ldots,d\}$, on $K\in\calT^n$ by
\begin{align*}
	{\mathbf\Lambda} _{i,j}^\delta(\bbB) _{|K}
	\coloneqq
	\sum\limits_{m=1}^d [(\calA_K^\top)^{-1}]_{i,m} \,
	\hat{\mathbf\Lambda} _{m}^\delta(\hat\bbB) \, [\calA_K^\top]_{m,j} \, \in\bbR^{d\times d}_{\mathrm{S}},
\end{align*}
on any $K \in \calT^n$ and for each $i,j\in\{1,\ldots,d\}$. 
Similarly to before, ${\mathbf\Lambda} _{i,j}^\delta(\bbB)_{|K} \in \bbR^{d\times d}_{\mathrm{S}}$ is constant on each element $K\in\calT^n$ and depends continuously on $\bbB_{|K}$ for any $K\in\calT^n$ and for each $i,j\in\{1,\ldots,d\}$. 
Moreover, it follows from the construction that \eqref{eq:Lambda_delta_kettenregel} is fulfilled and it holds $\lambda_i^\delta(\hat\bbB) \in [0,1]$. Again, the proof is similar to \cite[Lem.~3.6]{sieber_2020}, using that $-\ln(\cdot)$ is strictly convex as a function on the set of positive definite matrices. 

We recall the following result from \cite[Lem.~2.9]{GT_2023_DCDS}.
\begin{lemma}
Let $\delta\in(0,1)$ and let all simplices $K\in\calT^n$ be non-obtuse. Then, for all $\bbB\in\bbS^n$ and all $K\in\calT^n$, it holds 
\begin{align}
	\label{eq:nabla_ln} \nonumber
	- \int_{K} \nabla \bbB : \nabla\bbI_1^n [ \beta_\delta(\bbB)^{-1} ]  \dx
	&\geq 
	\frac1d \norm{\nabla\mathrm{I}_1^n \trace\ln \beta_\delta(\bbB)}_{L^2(K)}^2
	\\
	&= \frac1d  \norm{\nabla\mathrm{I}_1^n \trace\ln (\beta_\delta(\bbB)^{-1})}_{L^2(K)}^2.
\end{align}
\end{lemma}

The goal is to show an analogue of the discrete energy inequality \eqref{eq:stability_FE} for the regularised discrete scheme \eqref{eq:system_FE_delta}.

\begin{lemma} \label{lemma:stability_delta}
Let $\delta\in(0,1)$ and $n\in\{0,\ldots,N_T-1\}$. Let $\bu^n \in C(\overline\Omega;\bbR^d)$ and $\bbB^n \in  C(\overline\Omega;\bbR^{d\times d}_\mathrm{S})$ be given and suppose that Assumptions \ref{assumptions} hold. Then,
any solution $(\bu_\delta^{n+1}, \bbB_\delta^{n+1}, \bX_\delta^{n+1}, \kappa_\delta^{n+1}) \in \bbU_0^n \times \bbS^n \times \mathbf{V}(\Gamma^n) \times \mathrm{V}(\Gamma^n)$ to \eqref{eq:system_FE_delta}, if it exists, satisfies
\begin{align}
    \label{eq:stability_delta} \nonumber
    & \frac12 \norm{\sqrt{\rho^n} \bu_\delta^{n+1}}_{L^2}^2 
    + \skp{ W_\delta(\bbB_\delta^{n+1}) }{1}_{\calT^n}^h 
    + \gamma \calH^{d-1}(\Gamma_\delta^{n+1})
    \\ \nonumber
    &\quad
    + \frac12 \norm{\sqrt{\mathrm{I}_0^n \rho^{n-1}} (\bu_\delta^{n+1} - \bI_2^n \bu^n)}_{L^2}^2
    + 2 \Delta t \norm{\sqrt{\mu^n} \D(\bu_\delta^{n+1})}_{L^2}^2 
    \\ \nonumber
    &\quad
    + \Delta t \skp{\tfrac{G}{2\lambda^n}}{ \trace( \beta_\delta(\bbB_\delta^{n+1}) + \beta_\delta(\bbB_\delta^{n+1})^{-1} - 2 \bbI) }_{\calT^n}^h
    + \Delta t \frac{\alpha G}{2d} \norm{\nabla \mathrm{I}_1^n \trace\ln \beta_\delta(\bbB_\delta^{n+1})  }_{L^2}^2
    \\
    &\leq
    \frac12 \norm{\sqrt{\mathrm{I}_0^n \rho^{n-1}} \bI_2^n \bu^n}_{L^2}^2
    + \skp{W_\delta(\bbB^n)}{1}_{\calT^n}^h 
    + \gamma \calH^{d-1}(\Gamma^{n})
    + \Delta t \skp{\rho^n \mathbf{f}_1^{n+1} + \mathbf{f}_2^{n+1}}{\bu_\delta^{n+1}},
\end{align}
where $W_\delta(\bbB) = \tfrac{G}{2} \trace(\bbB - g_\delta(\bbB) - \bbI)$.
\end{lemma}

\begin{proof}
The main idea is to use discrete analogues of the energy estimates and geometric inequalities used in \eqref{eq:stability}. In particular, we set $\bw=\Delta t \bu_\delta^{n+1}$ in \eqref{eq:u_FE_delta} and note the elementary identity $2x(x-y) = x^2-y^2 + (x-y)^2$ for all $x,y\in\bbR$, which gives
\begin{align}
    \label{eq:stability_proof_1} \nonumber
     0&= \frac12 \norm{\sqrt{\rho^n} \bu_\delta^{n+1}}_{L^2}^2
     - \frac12 \norm{\sqrt{\mathrm{I}_0^n \rho^{n-1}} \bI_2^n \bu^n}_{L^2}^2
     \\ \nonumber
     &\quad
     + \frac12 \norm{\sqrt{\mathrm{I}_0^n \rho^{n-1}} (\bu_\delta^{n+1} - \bI_2^n \bu^n)}_{L^2}^2
    + 2 \Delta t \norm{\sqrt{\mu^n} \D(\bu_\delta^{n+1})}_{L^2}^2 
    \\
    &\quad 
    + \Delta t \skp{ \bbI_1^n \beta_\delta(\bbB_\delta^{n+1}) - \bbI }{\nabla\bu_\delta^{n+1}}
    - \Delta t \gamma \dualp{\kappa_\delta^{n+1} \pmb\nu^n}{ \bu_\delta^{n+1}}_{\Gamma^n}
    + \Delta t \skp{\rho^n \mathbf{f}_1^{n+1} + \mathbf{f}_2^{n+1}}{\bu_\delta^{n+1}}.
\end{align}
In addition, we obtain with $\bbG = \tfrac12  \Delta t  G(\bbI - \bbI_1^n [\beta_\delta(\bbB_\delta^{n+1})^{-1}])$ in \eqref{eq:B_FE_delta} that
\begin{align}
    \label{eq:stability_proof_2} \nonumber
    0 &= 
    \frac{G}{2} \skp{\bbB_\delta^{n+1} - \bbB^n}{\bbI -  \beta_\delta(\bbB_\delta^{n+1})^{-1}}_{\calT^n}^h
    \\ \nonumber
    &\quad
    + \Delta t G \sum_{i,j=1}^d \skp{[\bu_\delta^{n+1}]_i \, \mathbf{\Lambda}_{i,j}^\delta(\bbB_\delta^{n+1}) }{\partial_{x_j} \bbI_1^n [\beta_\delta(\bbB_\delta^{n+1})^{-1}] }
    \\ \nonumber
    &\quad 
    + \Delta t \skp{\tfrac{G}{2 \lambda^n} (\bbB_\delta^{n+1} - \bbI) }{\bbI - \beta_\delta(\bbB_\delta^{n+1})^{-1}}_{\calT^n}^h 
     \\
    &\quad
    - \Delta t G\skp{\nabla\bu_\delta^{n+1}}{ \bbI_1^n \beta_\delta(\bbB_\delta^{n+1}) - \bbI }
    - \Delta t \frac{\alpha G}{2} \skp{\nabla\bbB_\delta^{n+1}}{\nabla\bbI_1^n [\beta_\delta(\bbB_\delta^{n+1})^{-1}]}.
\end{align}
Using \eqref{eq:lemma_reg1e}, we get
\begin{align}
    \label{eq:stability_proof_3}
    \frac{G}{2} \skp{\bbB_\delta^{n+1} - \bbB^n}{\bbI -  \beta_\delta(\bbB_\delta^{n+1})^{-1}}_{\calT^n}^h
    \geq
    \skp{W_\delta(\bbB_\delta^{n+1})}{1}_{\calT^n}^h
    - \skp{W_\delta(\bbB^n)}{1}_{\calT^n}^h.
\end{align}
Moreover, we have with \eqref{eq:lemma_reg1b}, \eqref{eq:lemma_reg1d} and standard properties of the trace, that
\begin{align}
    \label{eq:stability_proof_3b}
    \skp{\tfrac{G}{2 \lambda^n} (\bbB_\delta^{n+1} - \bbI) }{\bbI - \beta_\delta(\bbB_\delta^{n+1})^{-1}}_{\calT^n}^h 
    \geq 
    \skp{\tfrac{G}{2\lambda^n}}{ \trace( \beta_\delta(\bbB_\delta^{n+1}) + \beta_\delta(\bbB_\delta^{n+1})^{-1} - 2 \bbI) }_{\calT^n}^h
    \geq 0.
\end{align}
Besides, it follows from \eqref{eq:Lambda_delta_kettenregel} and $\bu_\delta^{n+1}\in\bbU_0^n$ that
\begin{align}
    \label{eq:stability_proof_4} \nonumber
    &\sum_{i,j=1}^d \skp{[\bu^{n+1}_\delta]_i \, \mathbf{\Lambda}_{i,j}^\delta(\bbB_\delta^{n+1}) }{\partial_{x_j} \bbI_1^n [\beta_\delta(\bbB_\delta^{n+1})^{-1}] }
    \\
    &= \skp{\bu_\delta^{n+1}}{\nabla \mathrm{I}_1^n \trace\ln (\beta_\delta(\bbB_\delta^{n+1})^{-1})}
    = 0.
\end{align}
In the case $\alpha>0$, we assume that as all simplices $K\in\calT^n$ are non-obtuse, see Assumptions \ref{assumptions}. Therefore, it holds with \eqref{eq:nabla_ln} that
\begin{align}
    \label{eq:stability_proof_5}
    - \alpha \skp{\nabla\bbB_\delta^{n+1}}{\nabla\bbI_1^n [\beta_\delta(\bbB_\delta^{n+1})^{-1}]}
    \geq
    \frac{\alpha}{d} \norm{\nabla\mathrm{I}_1^n \trace\ln \beta_\delta(\bbB_\delta^{n+1})}_{L^2}^2,
\end{align}
while in the case $\alpha=0$ we can ignore all terms with $\alpha$.

Next, we choose $\pmb\zeta=\gamma (\bX_\delta^{n+1}-\mathbf{id}_{\Gamma^n})$ in \eqref{eq:kappa_FE_delta}, which gives together with \cite[Lem.~57]{BGN_2019_handbook} that
\begin{align}
    \label{eq:stability_proof_6} \nonumber
    0 &= \gamma \dualp{\kappa_\delta^{n+1} \pmb\nu^n}{\bX_\delta^{n+1}-\mathbf{id}}_{\Gamma^n}^h
    + \gamma \dualp{\nabla_s \bX_\delta^{n+1}}{\nabla_s (\bX_\delta^{n+1}-\mathbf{id})}_{\Gamma^n}
    \\
    &\geq 
    \gamma \dualp{\kappa_\delta^{n+1} \pmb\nu^n}{\bX_\delta^{n+1}-\mathbf{id}}_{\Gamma^n}^h
    + \gamma \calH^{d-1}(\Gamma_\delta^{n+1})
    - \gamma \calH^{d-1}(\Gamma^n),
\end{align}
where $\Gamma_\delta^{n+1} \coloneqq \bX_\delta^{n+1}(\Gamma^n)$. 
Moreover, setting $\chi=-\Delta t \gamma\kappa_\delta^{n+1}$ in \eqref{eq:X_FE_delta} leads to
\begin{align}
    \label{eq:stability_proof_7}
    0 &= -\gamma  \dualp{(\bX_\delta^{n+1} - \mathbf{id}) \cdot \pmb\nu^n}{\kappa_\delta^{n+1}}_{\Gamma^n}^h
    + \Delta t \gamma \dualp{\bu_\delta^{n+1} \cdot \pmb\nu^n}{\kappa_\delta^{n+1}}_{\Gamma^n}.
\end{align}

Finally, we can deduce the regularised energy inequality \eqref{eq:stability_delta} from \eqref{eq:stability_proof_1}--\eqref{eq:stability_proof_7}.
\end{proof}

Next, we show the existence of solutions to the regularised system \eqref{eq:system_FE_delta}. 
Here, we use the unique solvability of the subsystem \eqref{eq:kappa_FE_delta}--\eqref{eq:X_FE_delta} with continuous dependence on the velocity field to eliminate the parametrization $\bX^{n+1}_\delta$ and the mean curvature $\kappa^{n+1}_\delta$.
Then, we apply a fixed-point argument to show the existence of a velocity field $\bu^{n+1}_\delta$ and a tensor $\bbB^{n+1}_\delta$. After that, we reconstruct $\bX^{n+1}_\delta$ and $\kappa^{n+1}_\delta$. 

\begin{lemma} \label{lemma:eq:existence_FE_delta}
Let $\delta\in(0,\frac12]$ and $n\in\{0,\ldots,N_T-1\}$. Let $\bu^n \in C(\overline\Omega;\bbR^d)$ and $\bbB^n \in  C(\overline\Omega;\bbR^{d\times d}_\mathrm{S})$ be given and suppose that Assumptions \ref{assumptions} hold.
Then, there exists at least one solution $(\bu_\delta^{n+1}, \bbB_\delta^{n+1}, \bX_\delta^{n+1}, \kappa_\delta^{n+1}) \in \bbU_0^n \times \bbS^n \times \mathbf{V}(\Gamma^n) \times \mathrm{V}(\Gamma^n)$ to \eqref{eq:system_FE_delta} which satisfies the discrete energy inequality \eqref{eq:stability_delta}. 
\end{lemma}
\begin{proof}
First, we consider the subsystem \eqref{eq:kappa_FE_delta}--\eqref{eq:X_FE_delta}. Let $\bu\in\bbU^n$ be given. Using \eqref{eq:assumptions_vertex_normal} in Assumptions \ref{assumptions}, we conclude from the proof of \cite[Lem.~66]{BGN_2019_handbook} that there exists a unique solution tuple $(\bX, \kappa) \in \mathbf{V}(\Gamma^n)\times \mathrm{V}(\Gamma^n)$ to the linear discrete system
\begin{subequations}
\label{eq:kappa_X_FE_sub}
\begin{alignat}{2}
    \label{eq:kappa_FE_sub}
    \dualp{\kappa \,\pmb\nu^n}{\pmb\zeta}_{\Gamma^n}^h
    + \dualp{\nabla_s \bX}{\nabla_s \pmb\zeta}_{\Gamma^n} 
    &= 0 
    \qquad && \forall\,\pmb\zeta \in \mathbf{V}(\Gamma^n),
    \\
    \label{eq:X_FE_sub}
    \frac{1}{\Delta t} \dualp{\bX \cdot \pmb\nu^n}{\chi}_{\Gamma^n}^h
    &= 
    \frac{1}{\Delta t} \dualp{\mathbf{id} \cdot \pmb\nu^n}{\chi}_{\Gamma^n}^h
    +\dualp{\bu \cdot \pmb\nu^n}{\chi}_{\Gamma^n} 
    \qquad && \forall\,\chi \in \mathrm{V}(\Gamma^n),
\end{alignat}
\end{subequations}
which depends continuously on $\bu\in\bbU^n$ as we are in finite dimensions, i.e, $(\bX, \kappa) = (\bX(\bu), \kappa(\bu))$. Note that the continuous dependency may not be uniform with respect to the mesh size $h^n>0$. 

We define the following inner product
\begin{align*}
    \skpp{(\bu,\bbB)}{(\bw,\bbG)}
    &\coloneqq 
    \skp{\bu}{\bw}
    + \skp{\bbB}{\bbG}_{\calT^n}^h,
\end{align*}
for any $(\bu,\bbB), (\bw,\bbG) \in \bbU_0^n \times \bbS^n $.
For given functions $\bu^n\in C(\overline\Omega;\bbR^d)$ and $\bbB^n \in C(\overline\Omega;\bbR^{d\times d}_{\mathrm{S}})$, we define the continuous mapping
$\calH\colon\, \bbU_0^n \times \bbS^n \to \bbU_0^n \times \bbS^n$ such that for any $(\bu,\bbB) \in \bbU_0^n \times \bbS^n$
\begin{align*}
    &\skpp{\calH (\bu,\bbB)}{(\bw,\bbG)}
    \\
    &\coloneqq \frac{1}{2\Delta t} \skp{\rho^n \bu - (\mathrm{I}^n_0\rho^{n-1}) \bI^n_2 \bu^n }{\bw}
    + \frac{1}{2\Delta t} \skp{\mathrm{I}^n_0\rho^{n-1} (\bu - \bI^n_2 \bu^n) }{\bw}
    \\ 
    &\quad
    + \frac12 \skp{\rho^n }{ [(\bI^n_2\bu^n \cdot\nabla) \bu ] \cdot \bw}
    - \frac12 \skp{\rho^n }{\bu \cdot [(\bI^n_2\bu^n  \cdot\nabla) \bw]}
    \\
    &\quad 
    + \skp{2 \mu^n \D(\bu)}{\D(\bw)}
    + G \skp{ \bbI_1^n \beta_\delta(\bbB)  - \bbI }{\nabla\bw}
    - \gamma \dualp{\kappa(\bu) \, \pmb\nu^n}{ \bw}_{\Gamma^n}
    \\ 
    &\quad
    - \skp{\rho^n \mathbf{f}_1^{n+1} + \mathbf{f}_2^{n+1}}{\bw} 
    +\frac{1}{\Delta t} \skp{\bbB - \bbB^n}{\bbG}_{\calT^n}^h
    - \sum_{i,j=1}^d \skp{[\bu]_i \,  \mathbf{\Lambda}_{i,j}^\delta(\bbB) }{\partial_{x_j} \bbG}
    \\
    &\quad 
     + \skp{\tfrac{1}{ \lambda^n} (\bbB - \bbI) }{\bbG}_{\calT^n}^h 
    - 2 \skp{\nabla\bu}{ \bbI_1^n \big[ \bbG  \beta_\delta(\bbB) \big]}
    + \alpha \skp{\nabla\bbB}{\nabla\bbG},
\end{align*}
for all $(\bw,\bbG) \in \bbU_0^n \times \bbS^n $.
The function tuple $(\bu,\bbB,\bX(\bu),\kappa(\bu))$ with $(\bX(\bu),\kappa(\bu)) \in \mathbf{V}(\Gamma^n) \times \mathrm{V}(\Gamma^n)$ given by \eqref{eq:kappa_X_FE_sub}, if it exists, is a solution to \eqref{eq:system_FE_delta} if and only if $(\bu,\bbB) \in \bbU_0^n \times \bbS^n$ is a zero of $\calH$, i.e.,
\begin{align*}
    \skpp{\calH (\bu,\bbB)}{(\bw,\bbG)}=0 \qquad \forall\, (\bw,\bbG) \in \bbU_0^n \times \bbS^n .
\end{align*}
In the following, we use the testing procedure from the proof of \eqref{eq:stability_delta} in order to prove the following assumption wrong: Let $R>0$ and suppose that the continuous mapping $\calH$ has no zero in the closed ball 
\begin{align*}
    \calB_R \coloneqq \left\{ (\bw,\bbG) \in \bbU_0^n \times \bbS^n 
    \ \mid \
    ||| (\bw,\bbG) ||| \leq R \right\},
\end{align*}
where $||| \cdot |||$ is the norm induced by $\skpp{\cdot}{\cdot}$.
Then, for such $R>0$, we define a continuous mapping $\calG_R: \calB_R \to \partial \calB_R$ by
\begin{align*}
    \calG_R (\bw,\bbG)
    \coloneqq 
    -R \, \frac{\calH (\bw,\bbG)}{ 
    ||| \calH  (\bw,\bbG) ||| }
    \quad\quad \forall \,  (\bw,\bbG)\in \calB_R.
\end{align*}
It follows from Brouwer's fixed-point theorem \cite[Chap.~8.1.4, Thm.~3]{evans_2010} that there exists at least one fixed-point $(\bu^R,\bbB^R)\in \calB_R$ of the mapping $\calG_R$, i.e., it holds $(\bu^R,\bbB^R) = \calG_R(\bu^R,\bbB^R)$ and, as $\calG_R\colon \, \calB_R\to \partial\calB_R$, we also have
\begin{align}
\label{eq:existence_R}
    ||| (\bu^R,\bbB^R) |||
    = ||| \calG^R (\bu^R,\bbB^R) |||
    = R.
\end{align}
Next, the strategy is to show
\begin{align*}
    0 < \skpp{\calH(\bu^R,\bbB^R)}{ (\bw^R,\bbG^R) } < 0
\end{align*}
for one specific tuple of test functions $(\bw^R,\bbG^R) \in \bbU_0^n \times \bbS^n$, supposed that $R>0$ is large enough. This will disprove the assumption that $\calH$ has no zero in $\calB_R$. For the first inequality, we apply the strategy from the proof of the discrete energy inequality \eqref{eq:stability_delta}, and for the second inequality, we use that $(\bu^R,\bbB^R) \in \calB_R$ is a fixed-point of $\calG_R$. Let $(\bX(\bu^R),\kappa(\bu^R)) \in \mathbf{V}(\Gamma^n) \times \mathrm{V}(\Gamma^n)$ be the unique solution to \eqref{eq:kappa_X_FE_sub} with $\bu$ replaced by $\bu^R$.
We define 
\begin{align*}
    (\bw^R,\bbG^R) 
    = \Big(
    \Delta t \bu^R, \,
    \tfrac 12 \Delta t f (\bbI - \bbI_1^n[\beta_\delta(\bbB^R)^{-1}])
    \Big) \in \bbU_0^n \times \bbS^n.
\end{align*}
Then, choosing $\pmb\zeta=\gamma (\bX(\bu^R) - \mathbf{id}_{\Gamma^n})$ and $\chi = \Delta t \gamma \kappa(\bu^R)$ in \eqref{eq:kappa_X_FE_sub}, we obtain analogously to the proof of the regularised energy inequality \eqref{eq:stability_delta} with a straightforward computation that
\begin{align*}
    &\skpp{\calH(\bu^R,\bbB^R)}{(\bw^R,\bbG^R) }
    \\
    &\geq
    \frac12 \norm{\sqrt{\rho^n} \bu^R}_{L^2}^2 
    + \skp{ W_\delta(\bbB^R) }{1}_{\calT^n}^h 
    + \gamma \calH^{d-1}(\bX(\bu^R))
    \\ 
    &\quad
    + \frac12 \norm{\sqrt{\mathrm{I}_0^n \rho^{n-1}} (\bu^R - \bI_2^n \bu^n)}_{L^2}^2
    + 2 \Delta t \norm{\sqrt{\mu^n} \D(\bu^R)}_{L^2}^2 
    \\ 
    &\quad
    + \Delta t \skp{\tfrac{G}{2\lambda^n}}{ \trace( \beta_\delta(\bbB^R) + \beta_\delta(\bbB^R)^{-1} - 2 \bbI) }_{\calT^n}^h
    + \Delta t \frac{\alpha G}{2d} \norm{\nabla \mathrm{I}_1^n \trace\ln \beta_\delta(\bbB^R) }_{L^2}^2
    \\
    &\quad
    - \frac12 \norm{\sqrt{\mathrm{I}_0^n \rho^{n-1}} \bI_2^n \bu^n}_{L^2}^2
    - \skp{W_\delta(\bbB^n)}{1}_{\calT^n}^h 
    - \gamma \calH^{d-1}(\Gamma^n)
    - \Delta t \skp{\rho^n \mathbf{f}_1^{n+1} + \mathbf{f}_2^{n+1}}{\bu^R},
\end{align*}
where, by using Hölder's, Young's and Korn's inequalities, the force term is controlled by
\begin{align*}
    - \Delta t \skp{\rho^n \mathbf{f}_1^{n+1} + \mathbf{f}_2^{n+1}}{\bu^R} 
    \geq - C \Delta t \norm{\rho^n \mathbf{f}_1^{n+1} + \mathbf{f}_2^{n+1}}_{L^2}^2 - \Delta t \min\{\mu_+, \mu_-\} \, \norm{\bbD(\bu^R)}_{L^2}^2,
\end{align*}
where the constant $C$ only depends on $\mu_\pm$ and the domain.
With the help of \eqref{eq:lemma_reg1g} and an inverse inequality, it follows analogously to \cite{barrett_boyaval_2009} that there exist positive constants $C_1(h^n)$, $C_2(h^n)$ which depend on the mesh size $h^n$ such that
\begin{align*}
    \skp{ W_\delta(\bbB^R) }{1}_{\calT^n}^h 
    \geq 
    C_1(h^n) R^{-1} \skp{\bbB^R}{\bbB^R}_{\calT^n}^h
    - C_2(h^n),
\end{align*}
which implies that there exists another positive constant $C_3(h^n)$ such that
\begin{align*}
\skpp{\calH(\bu^R,\bbB^R)}{(\bw^R,\bbG^R) }
\geq 
C \norm{\bu^R}_{L^2}^2
+ C_1(h^n) R^{-1} \skp{\bbB^R}{\bbB^R}_{\calT^n}^h
- C_3(h^n),
\end{align*}
which is greater than zero if $R>0$ is large enough.

On the other hand, as $(\bu^R,\bbB^R) \in \calB_R$ is a fixed-point of $\calG_R$, we have
\begin{align*}
    \skpp{\calH (\bu^R,\bbB^R)}{(\bw^R,\bbG^R)}
    &= - \frac{||| \calH^h (\bu^R,\bbB^R)  ||| }{R} 
    \skpp{\calG_R (\bu^R,\bbB^R)}{(\bw^R,\bbG^R)}
    \\
    &= - \frac{||| \calH^h (\bu^R,\bbB^R)  ||| }{R} 
    \skpp{(\bu^R,\bbB^R)}{(\bw^R,\bbG^R)}.
\end{align*}
Again, it follows analogously to \cite{barrett_boyaval_2009} by using \eqref{eq:lemma_reg1h} and an inverse inequality that there exist positive constants $C_3(h^n)$, $C_4(h^n)$ depending on $h^n$ such that
\begin{align*}
    \skpp{(\bu^R,\bbB^R)}{(\bw^R,\bbG^R)}
    \geq 
    C \norm{\bu^R}_{L^2}^2 
    + C_3(h^n) R^{-1} \skp{\bbB^R}{\bbB^R}_{\calT^n}^h - C_4(h^n) ,
\end{align*}
which is greater than zero if $R>0$ is large enough, which implies $\skpp{\calH (\bu^R,\bbB^R)}{(\bw^R,\bbG^R)}<0$ for $R>0$ large enough.
In summary, we have shown 
\begin{align*}
    0 < \skpp{\calH (\bu^R,\bbB^R)}{(\bw^R,\bbG^R)} < 0
\end{align*}
for $R>0$ large enough, which yields a contradition. Hence, assuming that $R>0$ is large enough, the assumption that $\calH$ does not have a zero in $\calB_R$ has been proven to be wrong. Thus, there exists a zero $(\bu^R,\bbB^R) \in \calB_R$ of $\calH$ and therefore we found a solution $(\bu^R,\bbB^R, \bX(\bu^R), \kappa(\bu^R))$ to \eqref{eq:system_FE_delta}.
The discrete energy inequality \eqref{eq:stability_delta} follows from Lemma \ref{lemma:stability_delta} as all assumptions are fulfilled.
\end{proof}

We now finalise the proof Theorem \ref{thm:stability_FE}.
In particular, we identify subsequences of solutions to the regularised system \eqref{eq:system_FE_delta} and show that they converge to a solution to the original discrete system \eqref{eq:system_FE2}, as $\delta\to 0$. The proof relies mainly on the regularised energy inequality \eqref{eq:stability_delta} and the fact that we are in finite dimensions. The precise statement is formulated in the following lemma.

\begin{lemma}
Let $\delta\in(0,\frac12]$ and $n\in\{0,\ldots,N_T-1\}$. Let $\Gamma^n$ be a closed polyhedral hypersurface, and let $\bu^n \in C(\overline\Omega;\bbR^d)$ and $\bbB^n \in  C(\overline\Omega;\bbR^{d\times d}_\mathrm{S})$ be given with $\bbB^n$ being positive definite on $\overline\Omega$. Further, suppose that Assumptions \ref{assumptions} hold. Then, there exists a non-relabeled subsequence of $\{ (\bu^{n+1}_\delta, \bbB^{n+1}_\delta, \bX^{n+1}_\delta, \kappa^{n+1}_\delta) \}_{\delta>0}$, where $(\bu^{n+1}_\delta, \bbB^{n+1}_\delta, \bX^{n+1}_\delta, \kappa^{n+1}_\delta) \in \bbU_0^n \times \bbS^n \times \mathbf{V}(\Gamma^n) \times \mathrm{V}(\Gamma^n)$ solves the regularised system \eqref{eq:system_FE_delta}, and there exist functions $(\bu^{n+1}, \bbB^{n+1}, \bX^{n+1}, \kappa^{n+1}) \in \bbU_0^n \times \bbS^n \times \mathbf{V}(\Gamma^n) \times \mathrm{V}(\Gamma^n)$ such that for the subsequence
\begin{align*}
    \bu_\delta^{n+1} \to \bu^{n+1},
    \quad 
    \bbB_\delta^{n+1} \to \bbB^{n+1},
    \quad 
    \bX_\delta^{n+1} \to \bX^{n+1},
    \quad 
    \kappa_\delta^{n+1} \to \kappa^{n+1},
\end{align*}
as $\delta\to0$. In addition, $\bbB^{n+1} \in \bbS^n$ is positive definite and $(\bu^{n+1}, \bbB^{n+1}, \bX^{n+1}, \kappa^{n+1})$ is a solution to \eqref{eq:system_FE2}. Moreover, all solutions to \eqref{eq:system_FE2} with $\bbB^{n+1}$ positive definite satisfy \eqref{eq:stability_FE}.
\end{lemma}

\begin{proof}
As $\beta_\delta(\bbB_\delta^{n+1})$ is positive definite, we deduce from \eqref{eq:stability_delta}, \eqref{eq:lemma_reg1g}, Hölder's, Young's and Korn's inequalities that there exists a positive constant $C$ which is independent of $\delta$ such that
\begin{align}
    \label{eq:stability_delta_proof} \nonumber
    & \norm{\bu_\delta^{n+1}}_{L^2}^2 
    + \nnorm{\mathrm{I}_1^n \abs{\bbB_\delta^{n+1}} }_{L^1}
    + \frac{1}{\delta} \nnorm{\mathrm{I}_1^n \abs{ [\bbB_\delta^{n+1}]_- } }_{L^1}
    + \Delta t  \nnorm{\mathrm{I}_1^n \abs{\beta_\delta(\bbB_\delta^{n+1})^{-1}} }_{L^1} 
    \\
    &\leq C \big( \norm{\mathbf{I}_2^n \bu^n}_{L^2}^2 
    + \skp{W_\delta(\bbB^n)}{1}_{\calT^n}^h
    + \calH^{d-1}(\Gamma^n) 
    + \Delta t \norm{\rho^n \mathbf{f}_1^{n+1} + \mathbf{f}_2^{n+1}}_{L^2}^2 \big).
\end{align}
As $\bbB^n$ is positive definite, we have that $\skp{W_\delta(\bbB^n)}{1}_{\calT^n}^h \to \skp{W(\bbB^n)}{1}_{\calT^n}^h$ in the limit $\delta\to 0$. Hence, the right-hand side of \eqref{eq:stability_delta_proof} is uniformly bounded in $\delta>0$.
Therefore, since $\bu^{n+1}_\delta$ and $\bbB^{n+1}_\delta$ are finite dimensional, we directly deduce from \eqref{eq:stability_delta_proof} that there exists a non-relabeled subsequence such that $\bu^{n+1}_\delta \to \bu^{n+1} \in \bbU_0^n$ and $\bbB^{n+1}_\delta \to \bbB^{n+1} \in \bbS^n$ pointwise on $\overline\Omega$, as $\delta\to 0$. 

Moreover, it follows from \eqref{eq:stability_delta_proof} that $\abs{ [\bbB_\delta^{n+1}]_{-} }$ vanishes in the limit $\delta\to 0$. Therefore, $\bbB^{n+1}$ is at least positive semi-definite. 
As $\beta_\delta$ is Lipschitz continuous with Lipschitz constant equal to $1$, and as $\bbB^{n+1}$ is positive semi-definite, we obtain $\beta_\delta(\bbB_\delta^{n+1}) \to \bbB^{n+1}$, as $\delta\to 0$.
Moreover, \eqref{eq:stability_delta_proof} implies that $\mathrm{I}_1^n \abs{\beta_\delta(\bbB_\delta^{n+1})^{-1}}$ converges pointwise to a function in $\calS_1^n$, which is identified with $\mathrm{I}_1^n \abs{[\bbB^{n+1}]^{-1}} \in \calS_1^n$ using the pointwise convergence of $\beta_\delta(\bbB_\delta^{n+1})$ to $\bbB^{n+1}$ and the fact that $\beta_\delta(\bbB_\delta^{n+1}) \beta_\delta(\bbB_\delta^{n+1})^{-1} = \bbI$. Therefore, $\beta_\delta(\bbB_\delta^{n+1})^{-1} \to [\bbB^{n+1}]^{-1}$, as $\delta\to 0$, and hence $\bbB^{n+1}$ is positive definite.
From this we also have $\trace g_\delta(\bbB_\delta^{n+1}) \to \trace \ln(\bbB^{n+1})$ and $\trace\ln \beta_\delta(\bbB_\delta^{n+1}) \to \trace\ln(\bbB^{n+1})$, as $\delta\to 0$.
Similarly, $\mathbf\Lambda^\delta_{i,j}(\bbB_\delta^{n+1}) \to \mathbf\Lambda_{i,j}(\bbB^{n+1})$, as $\delta\to 0$, where $\mathbf\Lambda_{i,j}(\bbB^{n+1})$ is defined similarly to $\mathbf\Lambda^\delta_{i,j}(\bbB_\delta^{n+1})$, with $\hat{\mathbf\Lambda}_i^\delta$ replaced by $\hat{\mathbf\Lambda}_i$, which is defined similarly to $\hat{\mathbf\Lambda}_i^\delta$ with $\lambda_i^\delta$ and $\beta_\delta$ replaced by $\lambda_i$ and the identity function, respectively, where $\lambda_i$ is defined similarly to $\lambda_i^\delta$ with $\beta_\delta$ replaced by the identity function.

We recall that for any given $\bu \in \bbU^n$, we find a unique solution tuple $(\bX, \kappa) \in \mathbf{V}(\Gamma^n)\times \mathrm{V}(\Gamma^n)$ to the linear discrete system \eqref{eq:kappa_X_FE_sub} which depends continuously on $\bu \in\bbU^n$, as we are in finite dimensions. 
Therefore, taking the limit $\delta\to 0$ gives $\bX^{n+1}_\delta\to \bX^{n+1} \in \mathbf{V}(\Gamma^n)$ and $\kappa^{n+1}_\delta \to \kappa^{n+1} \in \mathrm{V}(\Gamma^n)$, where $(\bX^{n+1},\kappa^{n+1}) \in \mathbf{V}(\Gamma^n)\times \mathrm{V}(\Gamma^n)$ is defined as the unique solution tuple to \eqref{eq:kappa_X_FE_sub} with $\bu$ replaced by $\bu^{n+1}$.
In summary, we can now pass to the limit in \eqref{eq:system_FE_delta} and obtain the existence of a solution to the unregularised system \eqref{eq:system_FE2}.
Fatou's lemma and the positive definiteness of $\bbB^n$ allow us to pass to the limit in \eqref{eq:stability_delta} to justify \eqref{eq:stability_FE} for the limit functions $(\bu^{n+1}, \bbB^{n+1}, \bX^{n+1}, \kappa^{n+1})$.

To show that an arbitrary solution $(\bu^{n+1}, \bbB^{n+1}, \bX^{n+1}, \kappa^{n+1})$ to \eqref{eq:system_FE2} with $\bbB^{n+1}$ positive definite satisfies the discrete energy inequality \eqref{eq:stability_FE}, a similar testing procedure as in Lemma \ref{lemma:stability_delta} has to be performed. This involves choosing the test functions $\bw = \Delta t \bu^{n+1}$ in \eqref{eq:u_FE2}, $\bbG = \frac12 \Delta t G (\bbI - \bbI_1^n[ \bbB^{n+1}]^{-1} )$ in \eqref{eq:B_FE2}, $\zeta = \gamma (\bX^{n+1} - \mathbf{id}_{\Gamma^n})$ in \eqref{eq:kappa_FE2} and $\chi = - \Delta t \gamma \kappa^{n+1}$ in \eqref{eq:X_FE2}.
We point out the inequality
\begin{align*}
    \frac{G}{2} \skp{\bbB^{n+1} - \bbB^n}{\bbI -  (\bbB^{n+1})^{-1}}_{\calT^n}^h
    &\geq
    \skp{W(\bbB^{n+1})}{1}_{\calT^n}^h
    - \skp{W(\bbB^n)}{1}_{\calT^n}^h,
\end{align*}
which follows directly from the convexity of $W(\cdot)$, and
\begin{align*}
    - \alpha \skp{\nabla\bbB^{n+1}}{\nabla\bbI_1^n [(\bbB^{n+1})^{-1}]}
    &\geq
    \frac{\alpha}{d} \norm{\nabla\mathrm{I}_1^n \trace\ln \bbB^{n+1}}_{L^2}^2,
\end{align*}
where the proof is similar to that of \eqref{eq:nabla_ln} and is based on the convexity of $-\ln(\cdot)$, requiring only minor adjustments. Besides these two inequalities, the remaining steps are identical to the proof of Lemma \ref{lemma:stability_delta}.
\end{proof}

\subsection{Semi-discrete scheme and equidistribution of nodes} \label{sec:equidistribution}

Similarly to other works, we have good mesh properties for our numerical method in practice. To see this, we first state a semi-discrete analogue of \eqref{eq:system_FE}. 

Let $\calT^h$ be an arbitrarily fixed regular partitioning of $\Omega$ into disjoint open simplices such that $\overline\Omega = \bigcup_{K\in\calT^h} \overline{K}$ and such that $\partial_\mathrm{D}\Omega$ is exactly matched by the sides in $\calT^h$. We define the finite element spaces $\calS_k^h$, $\bbS^h$, $\bbU^h$, $\bbP^h$ similarly to $\calS_k^n$, $\bbS^n$, $\bbU^n$, $\bbP^n$, where $k\geq0$. Moreover, we denote the corresponding interpolation operators by $\mathrm{I}_k^h$, $\bI_k^h$, $\bbI_k^h$ for scalar, vector or matrix valued functions, respectively. 
Given $\Gamma^h(t)$, we denote the exterior of $\Gamma^h(t)$ by $\Omega_+^h(t)$ and the interior of $\Gamma^h(t)$ by $\Omega_-^h(t)$, so that $\Gamma^h(t) = \partial\Omega_-^h(t)$. Similarly to before, we define the piecewise constant unit normal $\pmb\nu^h(t)$ to $\Gamma^h(t)$ such that $\pmb\nu^h(t)$ points into $\Omega_+^h(t)$.
Moreover, we define $\rho^h(t) \in \calS_0^h$, $\mu^h(t) \in \calS_0^h$ and $\lambda^h(t) \in \calS_0^h$ as semi-discrete analogues of \eqref{eq:def_density}.
In the following, we always integrate over the current interface $\Gamma^h(t)$ using $\dualp{\cdot}{\cdot}_{\Gamma^h(t)}$, where the interface $\Gamma^h(t)$ is described by the identity function $\bX^h(t) = \mathbf{id}_{|\Gamma^h(t)} \in \mathbf{V}(\Gamma^h(t))$. 

We then formulate the semi-discrete analogue of \eqref{eq:system_FE} as follows. Let $\bu^h(0) \in \bbU^h(0)$ and $\bbB^h(0)\in \bbS^h(0)$ be given with $\bbB^h(0)$ positive definite. Moreover, let $\Gamma^h(0)$ be a closed polyhedral hypersurface. Then, for all $t\in(0,T]$, find $\bu^h(t) \in \bbU^h$, $p^h(t) \in \bbP^h$, $\bbB^h(t) \in \bbS^h$, $\bX^h(t) \in \mathbf{V}(\Gamma^h(t))$ and $\kappa^h(t) \in \mathrm{V}(\Gamma^h(t))$, 
such that $\bbB^h(t)$ is positive definite and 
\begin{subequations}
\label{eq:system_semi}
\begin{align}
    \label{eq:u_semi} \nonumber
    0 &= \ddt \frac{1}{2} \skp{\rho^h \bu^h }{\bw}
    + \frac{1}{2} \skp{\rho^h \partial_t \bu^h }{\bw}
    -  \frac{1}{2} \skp{\rho^h \bu^h}{\partial_t \bw}
    \\ \nonumber
    &\quad
    + \frac12 \skp{\rho^h }{ [(\bu^h \cdot\nabla) \bu^h ] \cdot \bw}
    - \frac12 \skp{\rho^n }{\bu^{n+1} \cdot [(\bI^n_2\bu^n  \cdot\nabla) \bw]}
    \\ \nonumber
    &\quad 
    + \skp{2 \mu^h \D(\bu^h)}{\D(\bw)}
    - \skp{p^h}{\Div \bw}
    + G \skp{\bbB^h - \bbI}{\nabla\bw} 
    \\ 
    &\quad
    - \gamma \dualp{\kappa^h \,\pmb\nu^h}{ \bw}_{\Gamma^h(t)}
    - \skp{\rho^h \mathbf{f}_1^h + \mathbf{f}_2^h}{\bw} ,
    \\
    \label{eq:div_semi}
    0 &= \skp{\Div \bu^h}{q},
    \\
    \label{eq:B_semi} \nonumber
    0 &= 
    \ddt \skp{\bbB^h}{\bbG}_{\calT^h}^h
	- \skp{\bbB^h}{\partial_t \bbG}_{\calT^h}^h
    - \sum_{i,j=1}^d \skp{[\bu^h]_i \, \mathbf{\Lambda}_{i,j}(\bbB^h)}{\partial_{x_j} \bbG}
    \\
    &\quad 
     + \skp{\tfrac{1}{\lambda^n} (\bbB^h - \bbI) }{\bbG}_{\calT^h}^h 
    - 2 \skp{\nabla\bu^h}{ \bbI_1^h \left[ \bbG \bbB^h \right]}
    + \alpha \skp{\nabla\bbB^h}{\nabla\bbG},
    \\
    \label{eq:kappa_semi}
    0 &= \dualp{\kappa^h \,\pmb\nu^h}{\pmb\zeta}_{\Gamma^h(t)}^h
    + \dualp{\nabla_s \bX^h}{\nabla_s \pmb\zeta}_{\Gamma^h(t)},
    \\
    \label{eq:X_semi}
    0 &=  \dualp{\partial_t \bX^h \cdot \pmb\nu^h}{\chi}_{\Gamma^h(t)}^h
    - \dualp{\bu^h \cdot \pmb\nu^h}{\chi}_{\Gamma^h(t)},
\end{align}
\end{subequations}
for all $t\in(0,T]$ and any $(\bw, q, \bbG, \pmb\zeta, \chi) \in \bbU^h \times \bbP^h \times \bbS^h \times \mathbf{V}(\Gamma^h(t)) \times \mathrm{V}(\Gamma^h(t))$.
Here we define $\mathbf{f}_i^h(t) \coloneqq \bI_2^h \, \mathbf{f}_i(t) \in \bbU^h$, $i\in\{1,2\}$.

For the semi-discrete scheme \eqref{eq:system_semi}, we have the following result.
\begin{theorem}\label{thm:stability_semi}
For any $t\in(0,T]$, let $\bu^h(t) \in \bbU^h(t)$, $p^h(t) \in \bbP^h(t)$, $\bbB^h(t) \in \bbS^h(t)$ with $\bbB^h(t)$ positive definite, $\bX^h(t) \in \mathbf{V}(\Gamma^h(t))$ and $\kappa^h(t) \in \mathrm{V}(\Gamma^h(t))$ be a solution to the semi-discrete system \eqref{eq:system_semi}. Moreover, let Assumptions \ref{assumptions} hold with $\Gamma^n$ replaced by $\Gamma^h(t)$.
\begin{itemize}
\item[(i)] It holds 
\begin{align}
    \label{eq:stability_semi} \nonumber
    &\ddt \Big( \frac12 \norm{\sqrt{\rho^h} \bu^h}_{L^2}^2
    + \skp{W(\bbB^h)}{1}_{\calT^h}^h 
    + \gamma \calH^{d-1}(\Gamma^h(t)) \Big)
    + \norm{ \sqrt{2\mu^h} \D(\bu^h)}_{L^2}^2
    \\ \nonumber
    &\quad 
    + \skp{\tfrac{G}{2\lambda^h}}{ \trace( \bbB^h + [\bbB^h]^{-1} - 2 \bbI) }_{\calT^h}^h
    + \frac{\alpha G}{2d} \norm{\nabla\mathrm{I}_1^h \ln \det \bbB^h}_{L^2}^2
    \\
    &\leq \skp{\rho^h \mathbf{f}_1^h + \mathbf{f}_2^h}{\bu^h},
\end{align}
where $W(\bbB) = \tfrac{G}{2} \trace(\bbB - \ln\bbB - \bbI)$.


\item[(ii)] For any $t\in(0,T]$, it holds that $\Gamma^h(t)$ is a conformal polyhedral hypersurface.
Moreover, if $d=2$, then $\Gamma^h(t)$ is weakly equidistributed for all $t\in(0,T]$, i.e., any two neighbouring elements of the curve $\Gamma^h(t)$ have equal lengths or are parallel.

\end{itemize}
\end{theorem}

\begin{proof}
\begin{itemize}
\item[(i)] This follows similarly to \eqref{eq:stability} and \eqref{eq:stability_FE} by choosing the test functions $\bw=\bu^h$, $q=p^h$, $\bbG = \tfrac12 G (\bbI - \bbI_1^h[(\bbB^h)^{-1}])$, $\chi=-\gamma\kappa^h$ and $\pmb\zeta = \gamma \partial_t\bX^h$ in the semi-discrete scheme \eqref{eq:system_semi}.


\item[(ii)] This follows directly from Equation \eqref{eq:kappa_semi}, see also \cite[Def.~60 and Thm.~62]{BGN_2019_handbook}.
\end{itemize}
\end{proof}

The second statement in Theorem \ref{thm:stability_semi} motivates why the vertices of $\Gamma^n$ in our fully discrete scheme \eqref{eq:system_FE} are well distributed in practice. 
In \cite[Sect.~4.1]{BGN_2008_pfem_three_dimensions} it is discussed that for $d=3$ the geometric property of (ii) means that the interface mesh has a good mesh quality.

There are also other methods in the literature that lead to a good mesh quality. One possibility is to use an implicit interface description that, on the fully discrete level, leads to a weakly equidistributed mesh for $d=2$ and to a conformal polyhedral hypersurface for $d=3$. However, in general it is not clear whether solutions exist then, see also \cite[Thm.~82]{BGN_2019_handbook}. For other approaches, we refer to \cite[Sect.~4.6]{BGN_2019_handbook}.

\subsection{XFEM for the conservation of the phase volumes} \label{sec:XFEM}

We have seen in \eqref{eq:volume_conservation} that the volume of the two phases $\Omega_\pm(t)$ is conserved in time for sufficiently smooth solutions to \eqref{eq:system_weak}. The key feature is to use the characteristic function $\chi_{|\Omega_-(t)}$ as a test function in \eqref{eq:div_weak}. 
By this idea, we have volume conservation for the semi-discrete scheme.

In the following, we allow for a time-dependent discrete pressure space $\bbP^h(t)$, while $\bbU^h$ and $\bbS^h$ are fixed.

\begin{lemma}
For any $t\in(0,T]$, let $\bu^h(t) \in \bbU^h$, $p^h(t) \in \bbP^h(t)$, $\bbB^h(t) \in \bbS^h$ with $\bbB^h(t)$ positive definite, $\bX^h(t) \in \mathbf{V}(\Gamma^h(t))$ and $\kappa^h(t) \in \mathrm{V}(\Gamma^h(t))$ form a solution to the semi-discrete scheme \eqref{eq:system_semi}. 
If $\chi_{|\Omega_-^h(t)} \in \bbP^h(t)$, then the volume of $\Omega_-^h(t)$ is conserved, i.e., it holds
\begin{align}
\ddt \mathcal{L}^d(\Omega_-^h(t)) = 0,
\end{align}
where $\mathcal{L}^d(\cdot)$ denotes the $d$-dimensional Lebesgue measure.
\end{lemma}

\begin{proof}
This follows similarly to \eqref{eq:volume_conservation} by using the transport identity \eqref{eq:transport_subdomain}, 
\eqref{eq:X_semi} with $\chi=1$, the Gauss theorem and \eqref{eq:div_semi} with $q=\chi_{|\Omega_-^h(t)}$, that
\begin{align*}
    \ddt \int_{\Omega^h_-(t)} 1 \dx
    &= \dualp{\partial_t \bX^h}{\pmb\nu^h}_{\Gamma^h(t)}^h
    = \dualp{\bu^h}{\pmb\nu^h}_{\Gamma^h(t)}
    = \skp{\Div \bu^h}{\chi_{|\Omega^h_{-}(t)}}
    = 0.
\end{align*}
\end{proof}

This result for the semi-discrete scheme motivated the authors of \cite{BGN_2013_stokes} to extend the pressure space $\bbP^n$ by a single basis function, that is, $\chi_{|\Omega_-^n}$, which leads to good volume conservation properties in practice. This is also referred to as the $\text{XFEM}_\Gamma$ approach, as the extra contributions in \eqref{eq:div_FE} due to $\chi_{|\Omega_-^n}$ can be written in terms of integrals over the interface $\Gamma^n$ by using the Gauss theorem, i.e., 
\begin{align*}
	0 = \skp{\Div \bu^{n+1}}{\chi_{|\Omega_-^n}}  = \dualp{\bu^{n+1}}{\pmb\nu^n}_{\Gamma^n}.
\end{align*}
Unfortunately, in general we do not have exact volume conservation with the $\text{XFEM}_\Gamma$ approach for the fully discrete scheme \eqref{eq:system_FE}. In fact, choosing the test function $\chi=1$ in \eqref{eq:X_FE}, using the Gauss theorem and noting that $\chi_{|\Omega_{-}^n} \in \bbP^n$, yields
\begin{align*}
    \frac{1}{\Delta t}
    \dualp{ \bX^{n+1} - \bX^n}{\pmb\nu^n}_{\Gamma^n}^h
    = \dualp{\bu^{n+1}}{\pmb\nu^n}_{\Gamma^n}
    = \skp{\Div \bu^{n+1}}{\chi_{|\Omega^n_{-}}}
    = 0,
\end{align*}
but, in general, it holds
\begin{align*}
\frac{1}{\Delta t} (  \mathcal{L}^d(\Omega_{-}^{n+1}) - \mathcal{L}^d({\Omega_{-}^{n}})  )
    & \neq \frac{1}{\Delta t}
    \dualp{ \bX^{n+1} - \bX^n}{\pmb\nu^n}_{\Gamma^n}^h.
\end{align*}
However, we observe very good volume conservation properties for the fully discrete scheme, as we show in Section \ref{sec:numerics}.
Similarly to other XFEM approaches, one cannot expect the inf-sup condition \eqref{eq:LBB} for $(\bbU^n,\hat\bbP^n)$, where \cblue{$\hat\bbP^n = L^2_0(\Omega) \cap \bbP^n$}
with
\begin{align*}
    \bbP^n = \calS_1^n + \mathrm{span}\, \{\chi_{|\Omega_{-}^n}\}
    \quad \text{or} \quad
    \bbP^n = \calS_1^n + \calS_0^n +\mathrm{span}\, \{\chi_{|\Omega_{-}^n}\}.
\end{align*}
Therefore, (iii) in Theorem \ref{thm:stability_FE} in general does not hold, while (i), (ii) and (iv) in Theorem \ref{thm:stability_FE} and also Theorem \ref{thm:stability_FE_sum} remain valid.
For more details, we refer to \cite{BGN_2013_stokes}.

As a natural generalisation of \cite[Lem.~5]{BGN_2013_stokes}, we can show with the $\text{XFEM}_\Gamma$ approach and in the absence of external forces that the reduced discrete system \eqref{eq:system_FE2} is capable of exactly recovering discrete stationary solutions. Specifically, consider the initial conditions $\bu^0 = \mathbf{0}$, $\bbB^0 = \bbI$ and let $\Gamma^0$ be a polyhedral surface with a constant discrete mean curvature $\overline{\kappa}\in\bbR$, i.e., $\Gamma^0$ is an equidistributed polygonal curve for $d=2$ and a conformal approximation of a sphere for $d=3$, see \cite{BGN_2008_pfem_three_dimensions}.
If in addition $\mathbf{f}_1^{1}=\mathbf{f}_2^{1}=\mathbf{0}$ and $\chi_{|\Omega_-^0} \in \bbP^0$, then 
\begin{align*}
    (\bu^1, \bbB^1, \bX^1, \kappa^1) = (\mathbf{0}, \bbI, \bX^0, \overline{\kappa}) \in \bbU^0_0 \times \bbS^0 \times \mathbf{V}(\Gamma^0) \times \mathrm{V}(\Gamma^0) 
\end{align*}
is a solution to \eqref{eq:system_FE2} for $n=0$.
However, as in Theorem \ref{thm:stability_FE}, the issue of uniqueness remains unsolved. 

The authors of \cite{GNZ_2023} have recently introduced a stable fully discrete scheme for two-phase Navier--Stokes flow with exact volume conservation. They use a time-weighted approximation $\pmb\nu^{n+\frac12}$ of the unit normal, which, unlike the approach in \cite{BGN_2015_navierstokes}, results in a nonlinear scheme. \cblue{We refer to \cite[(3.12)]{GNZ_2023} for the precise definition.}
The existence of solutions to the scheme from \cite{GNZ_2023} has not yet been proven due to the nonlinear definition of the time-weighted approximation of the unit normal. However, numerical tests have shown that their scheme can be easily solved with a fixed-point iteration in practice. 
Let us consider our approximation for viscoelastic two-phase flow combined with the time-weighted approximation of the unit normal from \cite{GNZ_2023}, i.e., \eqref{eq:u_FE}--\eqref{eq:B_FE} combined with 
\begin{subequations}
\label{eq:kappa_X_FE3}
\begin{align}
\label{eq:kappa_FE3}
    0 &= \dualp{\kappa^{n+1} \pmb\nu^{n+\frac12}}{\pmb\zeta}_{\Gamma^n}^h
    + \dualp{\nabla_s \bX^{n+1}}{\nabla_s \pmb\zeta}_{\Gamma^n},
    \\
    \label{eq:X_FE3}
    0 &= \frac{1}{\Delta t} \dualp{(\bX^{n+1} - \mathbf{id}) \cdot \pmb\nu^{n+\frac12}}{\chi}_{\Gamma^n}^h
    - \dualp{\bu^{n+1} \cdot \pmb\nu^n}{\chi}_{\Gamma^n}.
\end{align}
\end{subequations}
In this case, the identity
\begin{align*}
\frac{1}{\Delta t} (  \mathcal{L}^d(\Omega_{-}^{n+1}) - \mathcal{L}^d({\Omega_{-}^{n}})  )
    & = \frac{1}{\Delta t}
    \dualp{ \bX^{n+1} - \bX^n }{\pmb\nu^{n+\frac12}}_{\Gamma^n}^h 
\end{align*}
is actually valid, see \cite[Lem.~3.1]{GNZ_2023}, and it leads to exact volume conservation, i.e.,
\begin{align*}
    \frac{1}{\Delta t} (  \mathcal{L}^d(\Omega_{-}^{n+1}) - \mathcal{L}^d({\Omega_{-}^{n}})  )
    & = \frac{1}{\Delta t}
    \dualp{ \bX^{n+1} - \bX^n }{\pmb\nu^{n+\frac12}}_{\Gamma^n}^h
    = \dualp{\bu^{n+1}}{\pmb\nu^n}_{\Gamma^n}
    = \skp{\Div \bu^{n+1}}{\chi_{|\Omega^n_{-}}}
    = 0.
\end{align*}
Therefore, if we consider \eqref{eq:u_FE}--\eqref{eq:B_FE} together with \eqref{eq:kappa_X_FE3}, we can show the same stability result as in Theorem \ref{thm:stability_FE}, i.e., all solutions to \eqref{eq:u_FE}--\eqref{eq:B_FE} together with \eqref{eq:kappa_X_FE3} with $\bbB^{n+1}$ positive definite, if they exist, satisfy the discrete energy inequality \eqref{eq:stability_FE}. 
Although now the existence of discrete solutions is open, we have the desirable property of exact volume conservation in the case $\chi_{|\Omega_-^n}\in\bbP^n$.



\subsection{Extension to a variable shear modulus} \label{sec:variable_shear_modulus}
We now discuss an approach to deal with a variable shear modulus
\begin{align*}
    G(\cdot,t) = 
    \begin{cases}
        G_+ & \text{ in } \Omega_{+}(t),
        \\
        G_- & \text{ in } \Omega_{-}(t),
    \end{cases}
\end{align*}
which can be important in applications.
For the natural analogue of the fully discrete system \eqref{eq:system_FE}, we cannot derive stability estimates with a variable shear modulus. To derive stability estimates similar to Theorem \ref{thm:stability_FE}, we would have to perform a testing procedure that involves setting $\bbG = \frac12 \bbI_1^n[G^n (\bbI - (\bbB^{n+1})^{-1} ) ]$ in \eqref{eq:B_FE}, where $G^n \in \calS_0^n$ is defined as in \eqref{eq:def_density}. 
However, this approach encounters the problem that we would obtain multiple terms containing $\bbI_1^n [ G^n (\bbB^{n+1})^{-1} ]$ that we cannot control. 
Thus, we need to consider a different approach.

An idea to avoid these problems is to include the variable shear modulus $G$ in a variational formulation and to find a possible discretization for it.
Therefore, we take the Frobenius inner product of \eqref{eq:B} with $G \bbC$, where $\bbC$ is a time-dependent test function, integrate over $\Omega$ and perform integration by parts over the subdomains $\Omega_\pm(t)$, and use the fact that $G$ is constant in the subdomains $\Omega_\pm(t)$. Applying \eqref{eq:X} and \eqref{eq:jump}, and in the case $\alpha>0$ also using \eqref{eq:jump2} and \eqref{eq:bc_B}, we obtain
\begin{align}
    \label{eq:B_weak_G} \nonumber
    0 &= \ddt \skp{G \bbB}{\bbC}
    - \skp{\bbB}{\partial_t (G \bbC)}
    - \skp{G\bbB}{(\bu\cdot\nabla)\bbC}
    + \skp{\tfrac{G}{ \lambda} (\bbB - \bbI)}{\bbC}
    \\
    &\quad
    - 2 \skp{\nabla\bu}{G \bbC\bbB}
    + \alpha \skp{ G \nabla\bbB}{\nabla\bbC},
\end{align}
where $\bbC(\cdot,t) \in H^1(\Omega; \bbR^{d\times d}_\mathrm{S})$ for a.e.~$t\in(0,T)$.

Let us note that a natural choice for the discretization of $\ddt \skp{G \bbB}{\bbC} - \skp{\bbB}{\partial_t (G \bbC)}$ is
\begin{align*}
    &\frac{1}{\Delta t} \skp{G^n \bbB^{n+1} - G^{n-1} \bbB^n}{\bbC}_{\calT^n}^h
    - \frac{1}{\Delta t} \skp{\bbB^{n}}{(G^n - G^{n-1}) \bbC}_{\calT^n}^h
    \\
    &= \frac{1}{\Delta t} \skp{G^n (\bbB^{n+1} - \bbB^n)}{\bbC}_{\calT^n}^h.
\end{align*}

Then, the fully discrete system with variable shear modulus reads as follows. 
Let $\Gamma^0$ be an approximation of $\Gamma(0)$ and let $\bu^0 \in \bbU^0$ and $\bbB^0 \in \bbS^0$ positive definite be approximations of $\bu_0$ and $\bbB_0$, respectively.
Then, for any $n\in\{0,\ldots,N_T-1\}$, find $(\bu^{n+1}, p^{n+1}, \bbB^{n+1}, \bX^{n+1}, \kappa^{n+1}) \in \bbU^n \times \bbP^n \times \bbS^n \times \mathbf{V}(\Gamma^n) \times \mathrm{V}(\Gamma^n)$ with $\bbB^{n+1}$ positive definite such that 
\begin{subequations}
\label{eq:system_FE_G}
\begin{align}
    \label{eq:u_FE_G} \nonumber
    0 &= \frac{1}{2\Delta t} \skp{\rho^n \bu^{n+1} - (\mathrm{I}^n_0\rho^{n-1}) \bI^n_2 \bu^n }{\bw}
    + \frac{1}{2\Delta t} \skp{\mathrm{I}^n_0\rho^{n-1} (\bu^{n+1} - \bI^n_2 \bu^n) }{\bw}
    \\ \nonumber
    &\quad
    + \frac12 \skp{\rho^n }{ [(\bI^n_2\bu^n \cdot\nabla) \bu^{n+1} ] \cdot \bw}
    - \frac12 \skp{\rho^n }{\bu^{n+1} \cdot [(\bI^n_2\bu^n  \cdot\nabla) \bw]}
    \\ \nonumber
    &\quad 
    + \skp{2 \mu^n \D(\bu^{n+1})}{\D(\bw)}
    - \skp{p^{n+1}}{\Div \bw}
    + \skp{G^n (\bbB^{n+1} - \bbI)}{\nabla\bw}
    \\
    &\quad
    - \gamma \dualp{\kappa^{n+1} \pmb\nu^n}{ \bw}_{\Gamma^n}
    - \skp{\rho^n \mathbf{f}_1^{n+1} + \mathbf{f}_2^{n+1}}{\bw},
    \\
    \label{eq:div_FE_G}
    0 &= \skp{\Div \bu^{n+1}}{q},
    \\
    \label{eq:B_FE_G} \nonumber
    0 &= 
    \frac{1}{\Delta t} \skp{G^n  (\bbB^{n+1} - \bbB^n)}{\bbC}_{\calT^n}^h
    - \sum_{i,j=1}^d \skp{G^n [\bu^{n+1}]_i \, \mathbf{\Lambda}_{i,j}(\bbB^{n+1})}{\partial_{x_j} \bbC}
    \\
    &\quad 
     + \skp{\tfrac{G^n }{\lambda^n} (\bbB^{n+1} - \bbI) }{\bbC}_{\calT^n}^h 
    - 2 \skp{\nabla\bu^{n+1}}{ G^n \bbI_1^n \left[ \bbC \bbB^{n+1} \right]}
    + \alpha \skp{G^n \nabla\bbB^{n+1}}{\nabla\bbC},
    \\
    \label{eq:kappa_FE_G}
    0 &= \dualp{\kappa^{n+1} \pmb\nu^n}{\pmb\zeta}_{\Gamma^n}^h
    + \dualp{\nabla_s \bX^{n+1}}{\nabla_s \pmb\zeta}_{\Gamma^n},
    \\
    \label{eq:X_FE_G}
    0 &= \frac{1}{\Delta t} \dualp{(\bX^{n+1} - \mathbf{id}) \cdot \pmb\nu^n}{\chi}_{\Gamma^n}^h
    - \dualp{\bu^{n+1} \cdot \pmb\nu^n}{\chi}_{\Gamma^n},
\end{align}
\end{subequations}
for any $(\bw, q, \bbC, \pmb\zeta, \chi) \in \bbU^n \times \bbP^n \times \bbS^n \times \mathbf{V}(\Gamma^n) \times \mathrm{V}(\Gamma^n)$,
and set $\Gamma^{n+1} = \bX^{n+1}(\Gamma^n)$.

Here we have the following result.

\begin{lemma} Let Assumptions \ref{assumptions} hold. Then, all solutions to \eqref{eq:system_FE_G} with $\bbB^{n+1}$ positive definite, if they exist, satisfy
\begin{align}
    \label{eq:stability_FE_G} \nonumber
    & \frac12 \norm{\sqrt{\rho^n} \bu^{n+1}}_{L^2}^2 
    + \frac12 \skp{G^n}{\trace(\bbB^{n+1} - \ln\bbB^{n+1} - \bbI)}_{\calT^n}^h 
    + \gamma \calH^{d-1}(\Gamma^{n+1})
    \\ \nonumber
    &\quad
    + \frac12 \norm{\sqrt{\mathrm{I}_0^n \rho^{n-1}} (\bu^{n+1} - \bI_2^n \bu^n)}_{L^2}^2
    + 2 \Delta t \norm{\sqrt{\mu^n} \D(\bu^{n+1})}_{L^2}^2 
    \\ \nonumber
    &\quad
    + \Delta t \skp{\tfrac{ G^n}{2\lambda^n} }{\trace(\bbB^{n+1} + [\bbB^{n+1}]^{-1} - 2\bbI)}_{\calT^n}^h
    + \Delta t \frac{\alpha}{2d} 
    \skp{G^n}{\abs{\nabla \mathrm{I}_1^n \trace\ln\bbB^{n+1}}^2}
    \\ \nonumber
    &\leq
    \frac12 \norm{\sqrt{\mathrm{I}_0^n \rho^{n-1}} \bI_2^n \bu^n}_{L^2}^2
    + \frac12 \skp{G^n}{\trace(\bbB^{n} - \ln\bbB^{n} - \bbI)}_{\calT^n}^h 
    + \gamma \calH^{d-1}(\Gamma^n)
    \\
    &\quad
    + \Delta t \skp{\rho^n \mathbf{f}_1^{n+1} + \mathbf{f}_2^{n+1}}{\bu^{n+1}} 
    + \frac12 \Delta t \skp{G^n \bu^{n+1}}{\nabla\mathrm{I}_1^n \trace\ln\bbB^{n+1} }.
\end{align}
\end{lemma}
\begin{proof}
This follows similarly to Lemma \ref{lemma:stability_delta} by choosing the test functions $\bw=\Delta t \bu^{n+1}$, $q=p^{n+1}$, $\bbG = \tfrac12 \Delta t (\bbI - \bbI_1^h[(\bbB^{n+1})^{-1}])$, $\chi=-\Delta t \gamma \kappa^{n+1}$ and $\pmb\zeta = \gamma (\bX^{n+1} - \mathbf{id}_{\Gamma^n}) $ in the discrete system \eqref{eq:system_FE_G}. 
Here we note that the last term in \eqref{eq:stability_FE_G} follows from \eqref{eq:Lambda_kettenregel}.
\end{proof}

We note that in the case when $G(\cdot,t) = G_+ = G_-$, the system \eqref{eq:system_FE_G} is equivalent to the original discrete system \eqref{eq:system_FE}. However, in the case where the shear moduli are unmatched, \eqref{eq:system_FE_G} allows for the inequality \eqref{eq:stability_FE_G} while the system \eqref{eq:system_FE} does not. This motivates why the fully discrete system \eqref{eq:system_FE_G}, where $G^n \in \calS_0^n$ is directly incorporated into the scheme, can actually be used in practice. This is particularly important when trying to study the interaction between two viscoelastic fluids with different shear moduli $G_\pm$.

\section{Numerical results} \label{sec:numerics}
In this section, we present numerical examples in two space dimensions to illustrate the practicability of the numerical method analysed in this work. 
For the implementation, we use the finite element toolbox FEniCS \cite{fenics_book_2012}, which includes a framework for multiple intersecting meshes \cite{fenics_multimesh_2019} and offers access to the linear algebra package PETSc \cite{petsc-user-ref_2021}.
\cmagenta{The simulations were performed on a workstation with an Intel i7-12700 (2.1 GHz) processor with 32 GB of main memory, and none took longer than 96 hours to complete.}

\subsection{Solution strategy}
Now, we describe the solution techniques for the numerical scheme \eqref{eq:system_FE}. To effectively implement the computational procedure, we follow these steps sequentially:
\begin{enumerate}
\item Define the computational domain and the parameters, and set up the initial configurations for both bulk and parametric meshes. Then, for $n\in\{0,\ldots,N_T-1\}$, repeat the following steps.
\item Perform local adaptation of the bulk mesh to improve accuracy in regions close to the interface.
\item Apply local refinement to the parametric mesh to ensure that the mesh resolution is sufficient to capture detailed features of the solution.
\item Compute a solution to the fully discrete system \eqref{eq:system_FE} with a fixed-point iteration.
\item Calculate benchmark quantities to evaluate the performance of the model and save the results for visualisation or further analysis.
\end{enumerate}

\cteal{In the following, we describe the local adaption procedures from steps 2 and 3, and the fixed-point method from step 4. The domain $\Omega$, the final time $T$, the model parameters and the initial meshes will be defined later in Sections \ref{sec:numerics_retraction} and \ref{sec:numerics_rising_bubble}, respectively. 
In all cases, we use a rectangular domain $\Omega \subset\bbR^2$ in two space dimensions. 
We always start with a polygonal curve $\Gamma^0$ consisting of $400$ elements, and with a uniform initial triangulation of the domain $\Omega$ with a mesh size $h_c=\frac{2}{N_c}$, where $N_c=20$. We note that this initial mesh will be refined close to the interface using the routine from step 2, as described below, before we compute a solution to the discrete system.}

In all calculations, we use the P2--P1 element for the velocity and pressure. 
In addition to that, we employ $\text{XFEM}_\Gamma$ to ensure good volume conservation properties in practice, i.e., we enrich the basis functions of the pressure space by the characteristic function of the bubble $\chi_{|\Omega_-^n}$, as described in Section \ref{sec:XFEM}. The finite element spaces $\bbS^n$ for the tensor $\bbB$ and the parametric finite element spaces $\mathrm{V}(\Gamma^n)$ and $\mathbf{V}(\Gamma^n)$ are defined as in Section \ref{sec:fem}.

\begin{figure}
\centering
\includegraphics[width=0.23\textwidth,trim={1.22cm 1.25cm 1.1cm 1.31cm},clip]{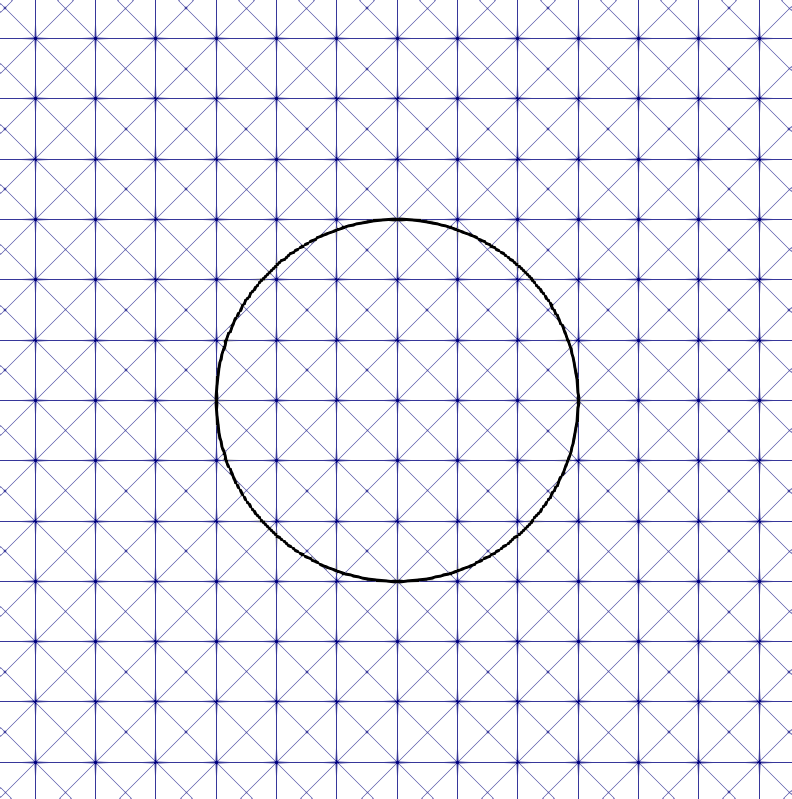} \,
\includegraphics[width=0.23\textwidth,trim={1.22cm 1.25cm 1.1cm 1.31cm},clip]{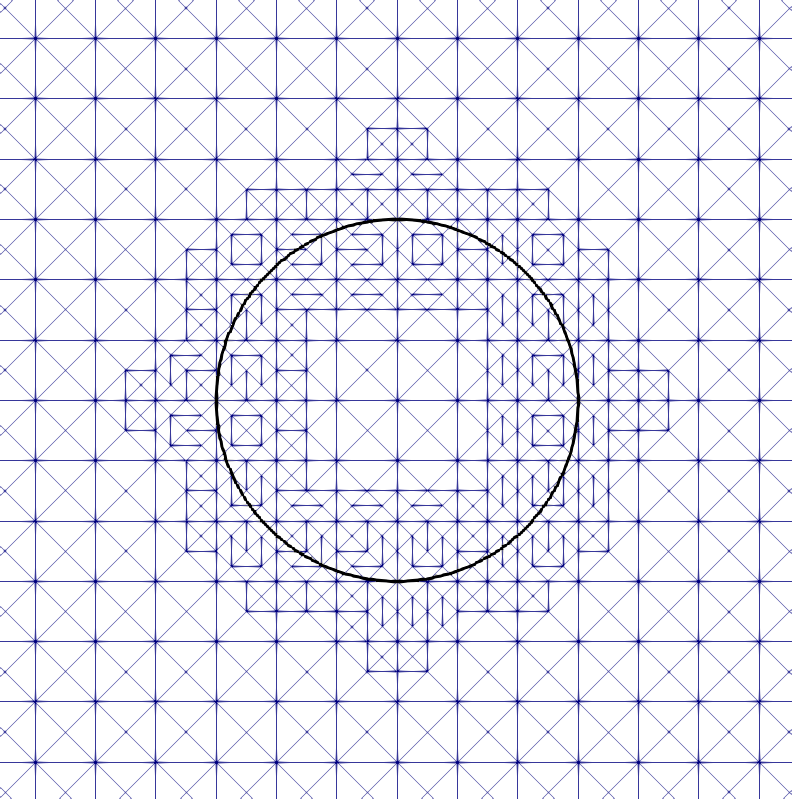} \,
\includegraphics[width=0.23\textwidth,trim={1.22cm 1.25cm 1.1cm 1.31cm},clip]{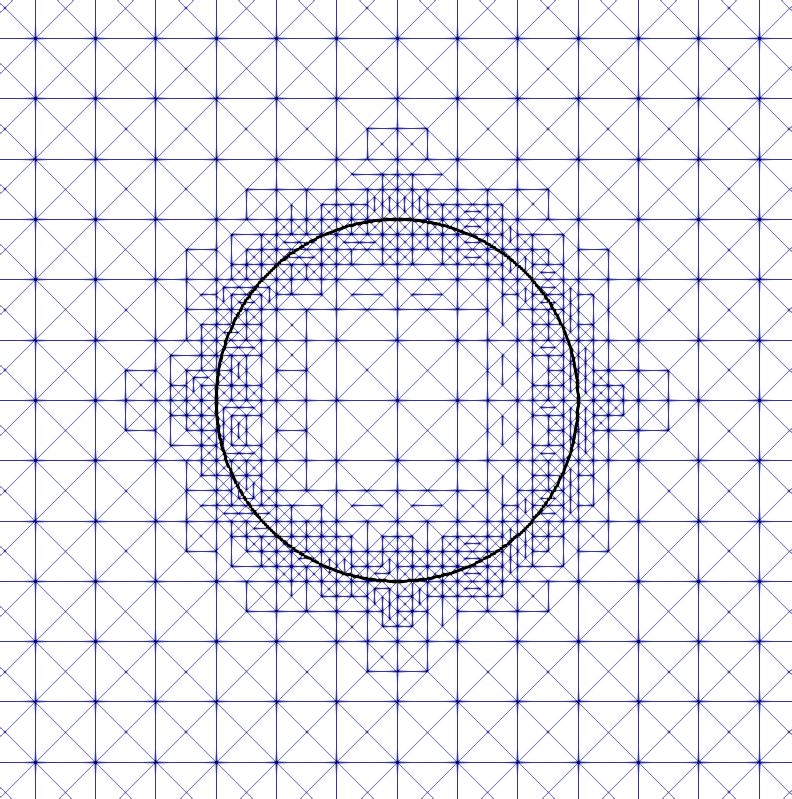} \,
\includegraphics[width=0.23\textwidth,trim={1.22cm 1.25cm 1.1cm 1.31cm},clip]{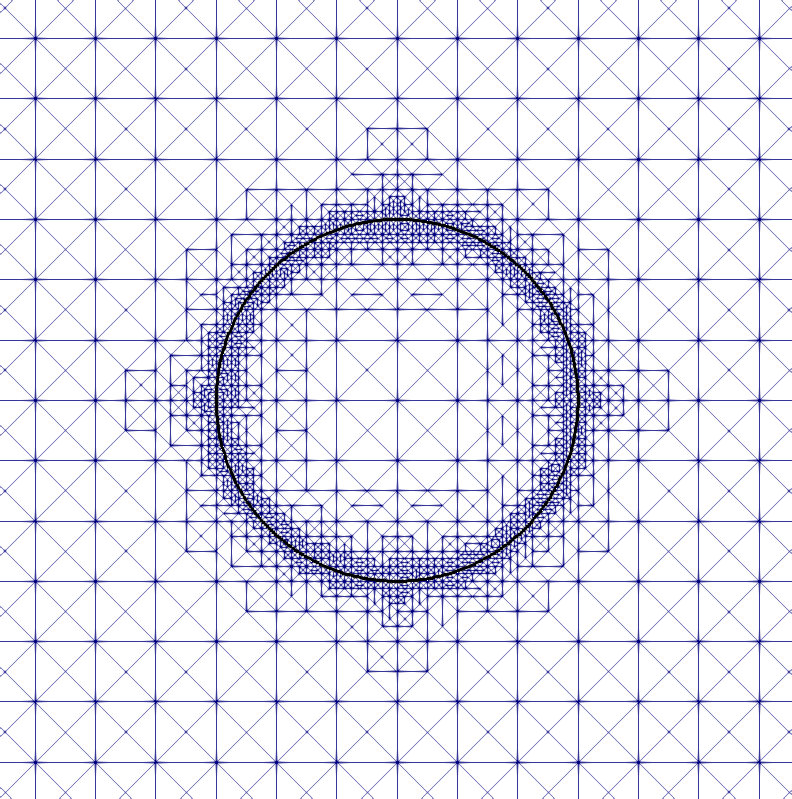}
\caption{A uniform bulk mesh, successively refined close to the interface.}
\label{fig:refinement}
\end{figure}

We employ a mesh refinement strategy inspired by similar concepts presented in \cite{BGN_2015_navierstokes}. Our goal is to create a mesh with a finer resolution close to the interface and a coarser resolution further away, optimising the mesh for computational accuracy and efficiency.
\cblue{At each time step,} we start with \cblue{the generic coarse triangulation $\mathcal{T}^{-1}$ with uniform mesh size $h_c = \frac{2}{N_c}$.}
We selectively mark elements \cblue{$K \in \mathcal{T}^{-1}$} for refinement based on the following criteria: An element \cblue{$K \in \mathcal{T}^{-1}$} is marked if $|K|>2 \frac{(h_f)^d}{d!}$ and if it intersects the \cblue{current} interface $\Gamma^n$, i.e., $K\cap \Gamma^{n} \not= \emptyset$. Here, we take $h_f = \frac{1}{8} h_c$. 
In addition, neighbouring cells are also marked for refinement to ensure spatial continuity and accuracy.
Refinement is performed using the bisection method \cblue{from \cite{rivara_1984} which is included in FEniCS \cite{fenics_book_2012}}, where we divide all marked elements into smaller, more finely resolved subelements. We note that, in general, more elements are refined than those initially selected for refinement to avoid hanging nodes. We repeat this refinement strategy iteratively until no further element satisfies the refinement criteria. We visualise three iterations of the refinement strategy in Figure \ref{fig:refinement}.

\begin{figure}
\centering
\includegraphics[width=0.32\textwidth]{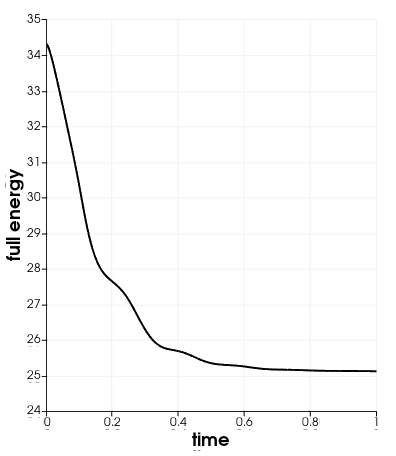}
\includegraphics[width=0.32\textwidth]{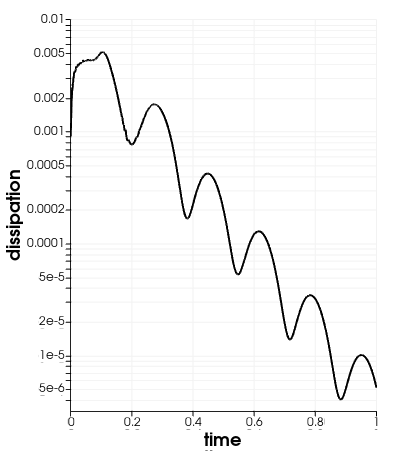}
\caption{\cmagenta{Time evolution of the full energy and the dissipations $\mathcal{D}^{n+1}$ defined in \eqref{eq:dissipation}.}}
\label{fig:retraction_energy_full}
\end{figure}

Similarly to \cite{BGN_2015_navierstokes}, we apply a local refinement strategy also to the interface mesh itself. In practice, the mesh vertices that represent the evolving surface $\Gamma^n$ are naturally well distributed. Thus, there is no need for additional mesh smoothing or redistribution procedures.
However, if the interface length grows significantly over time, it might be beneficial to refine parts of the mesh where elements have expanded too much.
To address this, we define the parameter $\mathrm{vol}_{max} \coloneqq \max \{ \calH^{d-1}(\sigma) \mid \sigma\in\calT(\Gamma^0) \}$,
which represents the length of the largest element of the initial interface $\Gamma^0$.
Then we identify all elements $\sigma \in \calT(\Gamma^n)$ of the current interface $\Gamma^n$ that have grown too large, that is, $\calH^{d-1}(\sigma) \geq \frac32 \mathrm{vol}_{max}$ for $\sigma\in\calT(\Gamma^n)$, and mark them for refinement.
The refinement process involves bisection of the identified elements, dividing them \cred{into two smaller segments} without changing the overall polyhedral surface $\Gamma^{n}$. 
\cblue{We stress that the parametric mesh is only refined, not coarsened, and that the simple refinement strategy for it is completely independent from the bulk mesh. In fact, apart from local refinement, no further modifications such as smoothing of the parametric mesh were required in our simulations. 
As mentioned before, we use a fine bulk mesh close to the moving interface to improve computational efficiency. Moreover, we recall that Theorem \ref{thm:stability_FE} is valid under bulk mesh adaption (i.e., local refinement and coarsening), and that the energy-dissipation inequality in Theorem \ref{thm:stability_FE_sum} is valid if no bulk mesh coarsening is performed. In practice we do employ coarsening, but our energy plots in Figures \ref{fig:retraction_energy_full}, \ref{fig:test94_energy}, and \ref{fig:test95_energy} suggest that we observe a discrete energy stability in practice nonetheless.}

As initial data for the velocity and the elastic Cauchy--Green tensor, we take $\bu^0$ and $\bbB^0$ to be the zero vector and the identity matrix, respectively, that is, $\bu^0 \coloneqq \mathbf{0} \in \bbU^0$ and $\bbB^0 \coloneqq \bbI \in \bbS^0$. Then, for each time step $n\in\{0,\ldots,N_T-1\}$, we compute a solution $(\bu^{n+1}, p^{n+1}, \bbB^{n+1}, \bX^{n+1}, \kappa^{n+1})$ to \eqref{eq:system_FE} with the help of a fixed-point iteration. We set $\bu^n_{0} \coloneqq \bu^n$ and $\bbB^n_{0} \coloneqq \bbB^n$. Then, for any $\ell\in\{0,1,\ldots\}$, compute $(\bu^{n}_{\ell+1}, p^{n}_{\ell+1}, \bbB^{n}_{\ell+1}, \bX^{n}_{\ell+1}, \kappa^{n}_{\ell+1}) \in \bbU^n \times \bbP^n \times \bbS^n \times \mathbf{V}(\Gamma^n) \times \mathrm{V}(\Gamma^n)$ such that 
\begin{subequations}
\label{eq:system_fixedpoint}
\begin{align}
    \label{eq:u_fixedpoint} \nonumber
    0 &= \frac{1}{2\Delta t} \skp{\rho^n \bu^{n}_{\ell+1} - (\mathrm{I}^n_0\rho^{n-1}) \bI^n_2 \bu^n }{\bw}
    + \frac{1}{2\Delta t} \skp{\mathrm{I}^n_0\rho^{n-1} (\bu^{n}_{\ell+1} - \bI^n_2 \bu^n) }{\bw}
    \\ \nonumber
    &\quad
    + \frac12 \skp{\rho^n }{ [(\bI^n_2\bu^n \cdot\nabla) \bu^{n}_{\ell+1} ] \cdot \bw}
    - \frac12 \skp{\rho^n }{\bu^{n}_{\ell+1} \cdot [(\bI^n_2\bu^n  \cdot\nabla) \bw]}
    \\ \nonumber
    &\quad 
    + \skp{2 \mu^n \D(\bu^{n}_{\ell+1})}{\D(\bw)}
    - \skp{p^{n}_{\ell+1}}{\Div \bw}
    + G \skp{\bbB^{n}_{\ell} - \bbI}{\nabla\bw}
    \\
    &\quad
    - \gamma \dualp{\kappa^{n}_{\ell+1} \pmb\nu^n}{ \bw}_{\Gamma^n}
    - \skp{\rho^n \mathbf{f}_1^{n+1} + \mathbf{f}_2^{n+1}}{\bw} ,
    \\
    \label{eq:div_fixedpoint}
    0 &= \skp{\Div \bu^{n}_{\ell+1}}{q},
    \\
    \label{eq:B_fixedpoint} \nonumber
    0 &= 
    \frac{1}{\Delta t} \skp{\bbB^{n}_{\ell+1} - \bbB^n}{\bbG}_{\calT^n}^h
    - \sum_{i,j=1}^d \skp{[\bu^{n}_{\ell}]_i \, \mathbf{\Lambda}_{i,j}(\bbB^{n}_{\ell})}{\partial_{x_j} \bbG}
    \\
    &\quad 
     + \skp{\tfrac{1}{\lambda^n} (\bbB^{n}_{\ell+1} - \bbI) }{\bbG}_{\calT^n}^h 
    - 2 \skp{\nabla \bu^{n}_{\ell}}{ \bbI_1^n \left[ \bbG \bbB^{n}_{\ell} \right]}
    + \alpha \skp{\nabla\bbB^{n}_{\ell+1}}{\nabla\bbG},
    \\
    \label{eq:kappa_fixedpoint}
    0 &= \dualp{\kappa^{n}_{\ell+1} \pmb\nu^n}{\pmb\zeta}_{\Gamma^n}^h
    + \dualp{\nabla_s \bX^{n}_{\ell+1}}{\nabla_s \pmb\zeta}_{\Gamma^n},
    \\
    \label{eq:X_fixedpoint}
    0 &= \frac{1}{\Delta t} \dualp{(\bX^{n}_{\ell+1} - \mathbf{id}) \cdot \pmb\nu^n}{\chi}_{\Gamma^n}^h
    - \dualp{\bu^{n}_{\ell+1} \cdot \pmb\nu^n}{\chi}_{\Gamma^n},
\end{align}
\end{subequations}
for any $(\bw, q, \bbG, \pmb\zeta, \chi) \in \bbU^n \times \bbP^n \times \bbS^n \times \mathbf{V}(\Gamma^n) \times \mathrm{V}(\Gamma^n)$. 


We observe that \eqref{eq:system_fixedpoint} is a linear scheme where the two--phase Navier--Stokes subsystem \eqref{eq:u_fixedpoint}--\eqref{eq:div_fixedpoint}, \eqref{eq:kappa_fixedpoint}--\eqref{eq:X_fixedpoint}, and the viscoelastic equation \eqref{eq:B_fixedpoint} can be solved independently. 
For a reduced version of \eqref{eq:system_fixedpoint}, where the discrete pressure is eliminated by using weakly divergence-free test functions $\bw \in \bbU_0^n$, the existence and uniqueness of solutions $\bu^n_{\ell+1}\in \bbU_0^n$, $\bX^n_{\ell+1} \in \mathbf{V}(\Gamma^n)$ and $\kappa^n_{\ell+1} \in \mathrm{V}(\Gamma^n)$ for \eqref{eq:u_fixedpoint} and \eqref{eq:kappa_fixedpoint}--\eqref{eq:X_fixedpoint} can be shown following the approach of \cite{BGN_2015_navierstokes}. 
If the inf-sup stability condition \eqref{eq:LBB} is satisfied, e.g., without $\text{XFEM}_\Gamma$, a corresponding unique pressure $p^{n}_{\ell+1} \in \hat{\bbP}^n$ can be reconstructed.
The existence and uniqueness of $\bbB^{n}_{\ell+1} \in \bbS^n$ solving \eqref{eq:B_fixedpoint} are guaranteed, since \eqref{eq:B_fixedpoint} clearly defines an invertible linear system of equations. 
To determine when to stop the iteration, we compute the $\ell^\infty$ norm of the residual of \eqref{eq:system_FE} with $(\bu^{n+1}, p^{n+1}, \bbB^{n+1}, \bX^{n+1}, \kappa^{n+1})$ replaced by $(\bu^{n}_{\ell+1}, p^{n}_{\ell+1}, \bbB^{n}_{\ell+1}, \bX^{n}_{\ell+1}, \kappa^{n}_{\ell+1})$, and stop if it is below a given small tolerance, e.g., $10^{-12}$. Then, we set 
\begin{align*}
    (\bu^{n+1}, p^{n+1}, \bbB^{n+1}, \bX^{n+1}, \kappa^{n+1}) \coloneqq
(\bu^{n}_{\ell+1}, p^{n}_{\ell+1}, \bbB^{n}_{\ell+1}, \bX^{n}_{\ell+1}, \kappa^{n}_{\ell+1})
\end{align*} 
and proceed with the next time step. We note that in practice, in all our numerical experiments, the converged matrix $\bbB^{n+1}$ is always positive definite. 
\cblue{Moreover, in all our numerical experiments, the fixed point iteration \eqref{eq:system_fixedpoint} required between 3 and 10 steps to converge. We observed a general trend that the number of sub-iterations increased when the elastic properties of the fluid became more dominant, i.e., for larger values of $G$ and $\lambda$.}

In practice, the variables $\bX^n_{\ell+1}$ and $\kappa^n_{\ell+1}$ are not explicitly computed during the fixed-point iteration \eqref{eq:system_fixedpoint}. Instead, $\bX^n_{\ell+1}$ and $\kappa^n_{\ell+1}$ are implicitly defined as unique solutions to the linear subsystem \eqref{eq:kappa_fixedpoint}--\eqref{eq:X_fixedpoint}, using the Schur complement method from \cite{BGN_2015_navierstokes}. 
This allows us to focus on solving for $\bu^{n+1}$, $p^{n+1}$, and $\bbB^{n+1}$ within the fixed-point iteration \eqref{eq:u_fixedpoint}--\eqref{eq:B_fixedpoint}, where $\bX^n_{\ell+1}$ and $\kappa^n_{\ell+1}$ are implicitly defined by \eqref{eq:kappa_fixedpoint}--\eqref{eq:X_fixedpoint}.
\cblue{The Navier--Stokes subsystem \eqref{eq:u_fixedpoint}--\eqref{eq:div_fixedpoint} is solved using the BiCGSTAB iterative solver from PETSc \cite{petsc-user-ref_2021} with right-oriented ``BFBt'' preconditioning for the Navier--Stokes component. This preconditioning method is explained in detail in \cite[Sect.~5]{BGN_2015_navierstokes} and is based on the techniques outlined in \cite[Sect.~9.2]{elman_silvester_wathen_2014}.}
The viscoelastic equation \eqref{eq:B_fixedpoint} corresponds to a linear system of equations with a positive definite system matrix. Therefore, well known solution techniques, such as a preconditioned conjugate gradient solver, can be applied.
After completion of the fixed-point iteration, $\bX^{n+1}$ and $\kappa^{n+1}$ are then calculated as unique solutions to the linear subsystem \eqref{eq:kappa_FE}--\eqref{eq:X_FE}. The detailed strategy for the elimination process of $\bX^n_{\ell+1}$ and $\kappa^n_{\ell+1}$ using the Schur complement approach is elaborated in \cite{BGN_2015_navierstokes}.

\begin{figure}
\centering
\includegraphics[width=0.32\textwidth]{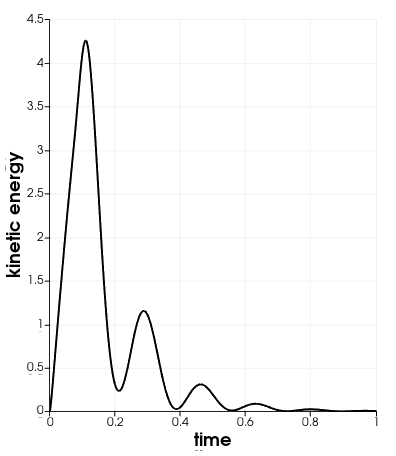}
\includegraphics[width=0.32\textwidth]{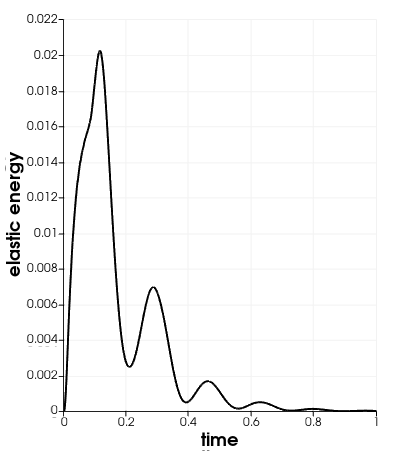}
\includegraphics[width=0.32\textwidth]{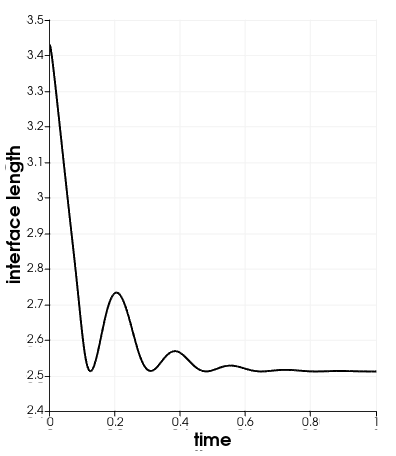}
\caption{\cmagenta{Time evolution of the energy contributions for the retraction experiment. From left to right: kinetic energy, elastic energy and interface length.}}
\label{fig:retraction_energy}
\end{figure}

\cmagenta{\subsection{Retraction of an elongated bubble}
\label{sec:numerics_retraction}

To illustrate the energy stability of the scheme \eqref{eq:system_FE}, we present a retraction experiment of an elongated bubble without external forces (i.e., $\mathbf{f}_1 = \mathbf{f}_2 = \mathbf{0}$).

We define the rectangular domain $\Omega = (0,2)^2 \subset\mathbb{R}^2$ for the numerical computation and we fix the no-slip boundary conditions for the velocity on the whole boundary of $\Omega$, i.e., $\partial_\mathrm{D}\Omega = \partial\Omega$. Moreover, we define the final time $T=1$ and use the time step size $\Delta t = 10^{-4}$.

The initial parametric surface $\Gamma^0$ consists of $400$ elements and is defined as a polygonal approximation of an ellipse with center $(1,1)^\top$ and horizontal and vertical axis of length $0.8$ and $0.2$, respectively. 
In particular, the vertices $\mathbf{q}_k^0$, $k\in\{1,\ldots,400\}$, of $\Gamma^0$ are set to
\begin{align*}
    \mathbf{q}_k^0 = 
    \begin{pmatrix}
        1 \\ 1
    \end{pmatrix}
    +
    \begin{pmatrix}
        0.8 \cos(\phi_k) \\ 0.2 \sin(\phi_k)
    \end{pmatrix},
    \quad 
    \phi_k = \frac{2\pi k}{400}, \quad k \in \{1,\ldots,400\}.
\end{align*}
For the initial triangulation of the domain $\Omega$, we consider a uniform partitioning with a mesh size $h_c = \frac{2}{N_c}$, where $N_c = 20$. This initial mesh will be refined close to the interface using the method from step 2 from above before we compute a solution to the discrete system.

The model parameters are fixed as
\begin{gather*}
    \rho_+ = \rho_- = 1, \quad
    \mu_+ = \mu_- = 0.1, \quad
    \gamma = 10, \quad
    \lambda_+ = \lambda_- = 0.01, \\
    G = 1, \quad
    \alpha = 10^{-2}, \quad
    \mathbf{f}_1 = \mathbf{f}_2 = \mathbf{0}.
\end{gather*}

\begin{figure}
\centering
\includegraphics[width=0.65\textwidth,trim={1cm 3cm 1cm 3cm},clip]{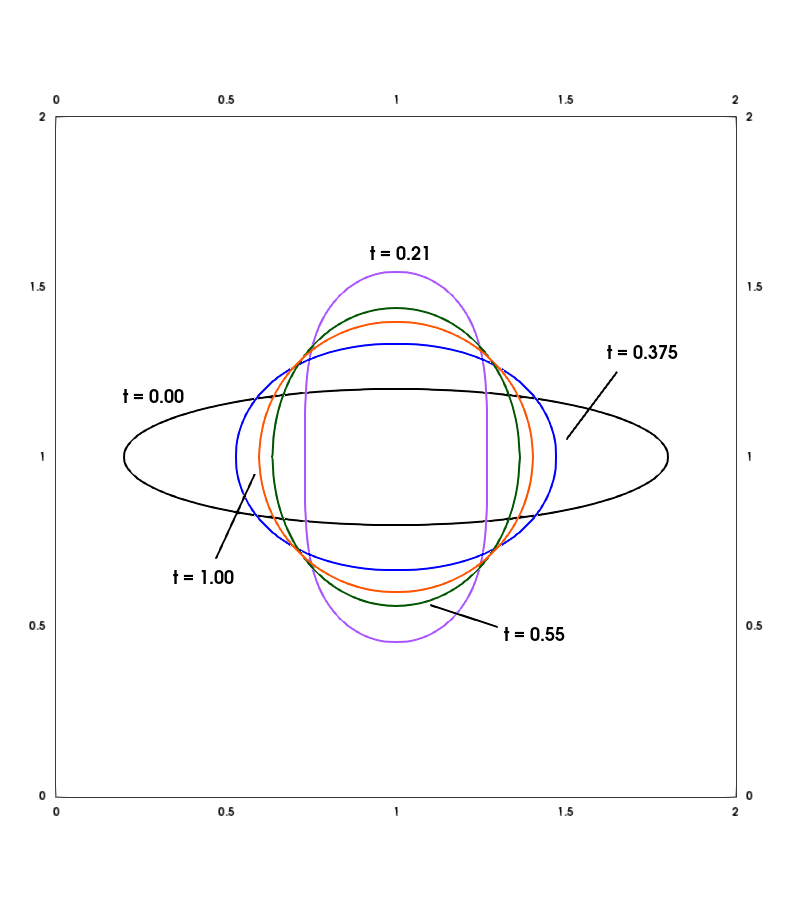}
\caption{\cmagenta{Time evolution of the interface in the retraction experiment at times $t \in \{0, 0.21, 0.375, 0.55, 1\}$.}}
\label{fig:retraction}
\end{figure}

In Figure \ref{fig:retraction}, we visualize the interface at different times $t \in \{0, 0.21, 0.375, 0.55, 1\}$. The bubble has an oscillatory behaviour and the amplitude becomes smaller over time until the oscillations are not visible any more. In Figure \ref{fig:retraction_energy}, we show the time evolution of the energy contributions from \eqref{eq:stability_FE}, i.e., the kinetic energy, the elastic energy and the length of the interface. 
As predicted by the discrete energy inequality \eqref{eq:stability_FE}, we observe a decay over time of the discrete energy which is shown in Figure \ref{fig:retraction_energy_full}.
To better illustrate the magnitude of the dissipative terms in the discrete energy inequality \eqref{eq:stability_FE}, we define the change in the discrete energy (excluding the contribution from body forces) by  
\begin{align}
    \label{eq:dissipation} \nonumber
    - \mathcal{D}^{n+1} &:= \frac12 \|\sqrt{\rho^n} \mathbf{u}^{n+1}\|_{L^2}^2 
    + \left( W(\mathbb{B}^{n+1}) , 1 \right)_{\mathcal{T}^n}^h 
    + \gamma \mathcal{H}^{d-1}(\Gamma^{n+1})
    \\ \nonumber
    &- \frac12 \| \sqrt{\mathrm{I}_0^n \rho^{n-1}} \mathbf{I}_2^n \mathbf{u}^n \|_{L^2}^2
    - \left(W(\mathbb{bbB}^n), 1\right)_{\mathcal{T}^n}^h 
    - \gamma \mathcal{H}^{d-1}(\Gamma^n)
    \\
    &- \Delta t \left( \rho^n \mathbf{f}_1^{n+1} + \mathbf{f}_2^{n+1}, \mathbf{u}^{n+1}  \right),
\end{align}
where, for this specific experiment, we recall that $\mathbf{f}_1^{n+1} = \mathbf{f}_2^{n+1} = \mathbf{0}$.
The dissipations $\mathcal{D}^{n+1}$ are expected to be non-negative due to \eqref{eq:stability_FE} and this is also what we observe in Figure \ref{fig:retraction_energy_full}.}

\subsection{The rising bubble experiment}
\label{sec:numerics_rising_bubble}

\cteal{The following numerical calculations are motivated by \cite{pillapakkam_2007_viscoelastic}. 
We define the rectangular domain $\Omega = (0,2) \times (0,4) \subset\bbR^2$ for our computations and we fix the no-slip boundary conditions for the velocity on the whole boundary of $\Omega$, i.e., $\partial_\mathrm{D}\Omega = \partial\Omega$. Moreover, we define the final time $T=1$ and the time step size $\Delta t = 10^{-4}$.
The initial parametric surface $\Gamma^0$ is defined as an equidistributed polygonal curve with a constant discrete mean curvature, consisting of $400$ elements. 
In particular, the vertices of $\Gamma^0$ are positioned on a circle of radius $0.3$ centred at the point $(1, 0.8)^\top$.
For the initial triangulation of the domain $\Omega$, we consider a uniform partitioning with a mesh size $h_c = \frac{2}{N_c}$, where $N_c = 20$. We note that this initial mesh will be refined close to the interface using the routine from step 2, as described above, before we compute a solution to the discrete system.


We choose the model parameters
\begin{align*}
    \rho_+ &= 1, \quad
    \rho_- = 0.1, \quad
    \mu_+ = \frac{10.25}{1+c_0}, \quad 
    \mu_- = 1.025, \quad
    \gamma = 10, \\
    \lambda_+ &\in \{0.05, 0.075\}, \quad
    \lambda_- = 10^{-3}, \quad
    G = \frac{c_0 \mu_+}{\lambda_+}, \quad
    \alpha = 10^{-2}.
\end{align*}
In this context, the constant $c_0 \geq 0$ represents a parameter that characterises the elastic contribution to the viscoelastic behaviour of the fluid. More specifically, adjusting the value of the constant $c_0 \geq 0$ leads to variations in the viscosity ratio
\begin{align*}
     \frac{\mu_+}{\mu_+ + G \lambda_+ } 
    = \frac{1}{1+c_0} \in (0, 1]
\end{align*}
for the outer phase $\Omega_+(t)$, while the total shear viscosity $\mu_+ + G \lambda_+ = 10.25$ for the outer phase remains fixed across all experiments. 
For the inner phase $\Omega_-(t)$, we select a very small relaxation parameter $\lambda_-$ to approximate the rheology of a Newtonian fluid.
This ensures that any observed differences are attributed solely to the viscoelasticity of the outer phase. 
The numerical examples presented here are computed with $c_0 \in \{0, 1, 19.5\}$. 
\cblue{In the case $c_0=0$, there is no contribution of elastic stresses in the momentum equation, thus we do not compute $\mathbb{B}$ and use the linear scheme of \cite{BGN_2015_navierstokes} for the two-phase Navier--Stokes problem without viscoelasticity. This corresponds to performing exactly one iteration of \eqref{eq:u_fixedpoint}--\eqref{eq:div_fixedpoint} without the term $G \left(  {\mathbb{B}^{n}_{\ell} - \mathbb{I}}, {\nabla\mathbf{w}} \right)$, combined with \eqref{eq:kappa_fixedpoint}--\eqref{eq:X_fixedpoint}.}

In addition to the model parameters, we account for gravitational forces within our simulations. Specifically, we consider the body acceleration
\begin{align*}
    \cred{\mathbf{f}_1} = -980 \,\mathbf{e}_2,
\end{align*}
with $\mathbf{e}_2 = (0, 1)^\top$ being the unit vector in positive $x_2$-direction, and we set $\cred{\mathbf{f}_2} = \mathbf{0}$.
}

\begin{figure}
\centering
\includegraphics[width=0.32\textwidth,trim={0 2cm 0 2cm},clip]{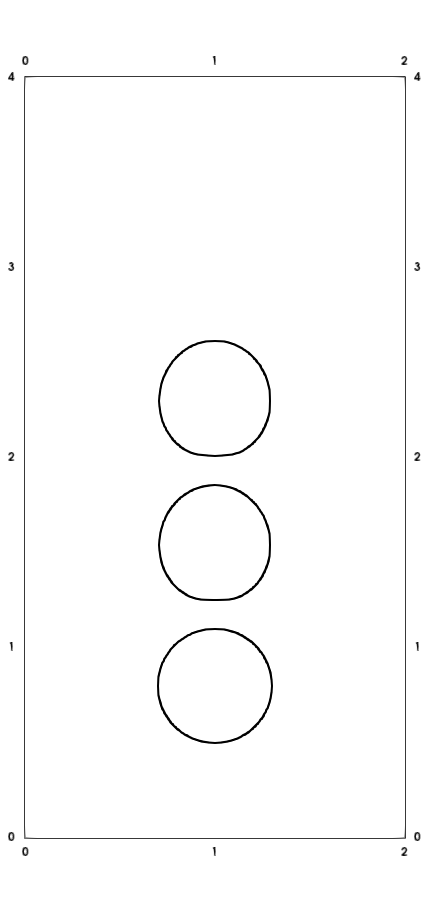}
\includegraphics[width=0.32\textwidth,trim={0 2cm 0 2cm},clip]{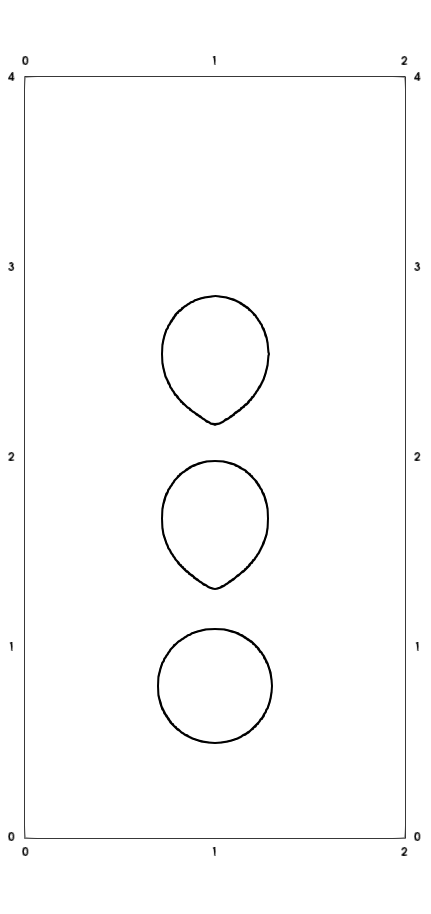}
\includegraphics[width=0.32\textwidth,trim={0 2cm 0 2cm},clip]{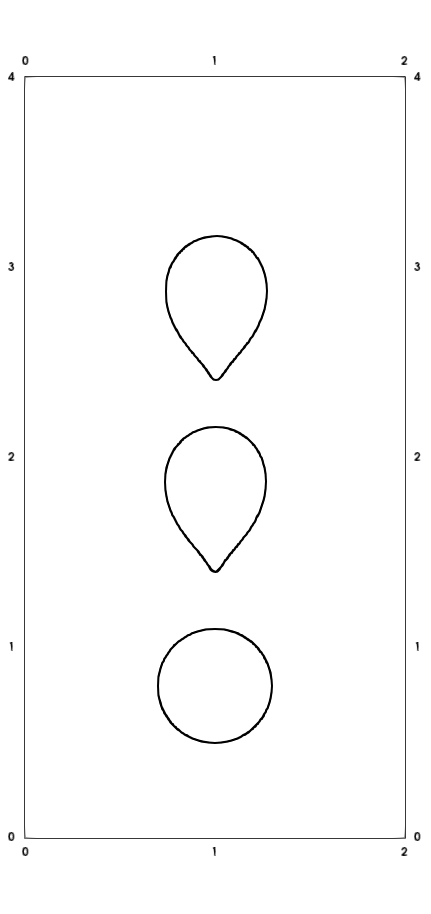}
\caption{Time evolution of the rising bubble for the cases $c_0=0$ (left), $c_0=1$ (center) and $c_0=19.5$ (right) with the relaxation parameter $\lambda_+ = 0.05$ at times $t\in\{0, 0.5, 1\}$.}
\label{fig:test94_simulation}
\end{figure}

In the following, we present the numerical results for two different relaxation parameters $\lambda_+\in\{0.05, 0.075\}$, where we discuss the influence of viscoelasticity by varying the constant $c_0\in\{0, 1, 19.5\}$.

\begin{figure}
\centering
\includegraphics[width=0.24\textwidth,trim={0 2cm 0 2cm},clip]{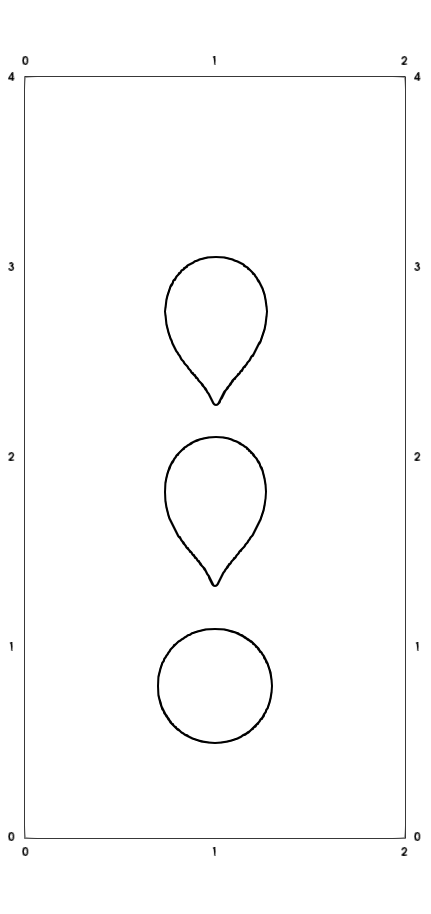}
\includegraphics[width=0.24\textwidth,trim={0 2cm 0 2cm},clip]{figures/test94_evolution.png}
\includegraphics[width=0.24\textwidth,trim={0 2cm 0 2cm},clip]{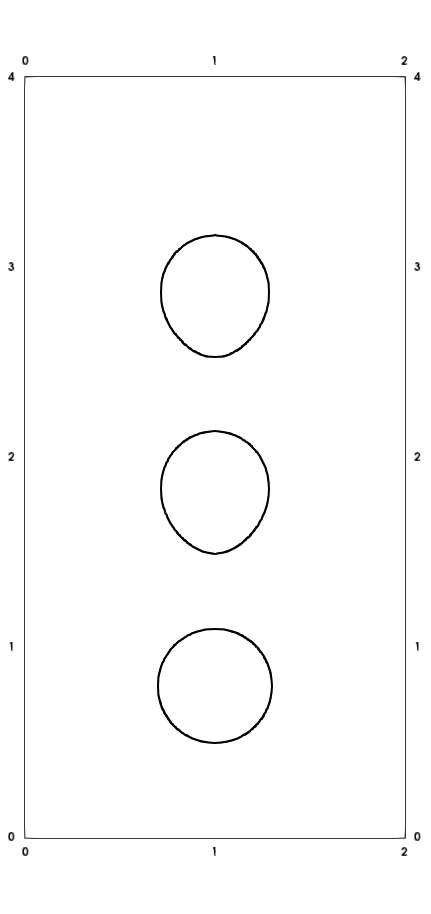}
\includegraphics[width=0.24\textwidth,trim={0 2cm 0 2cm},clip]{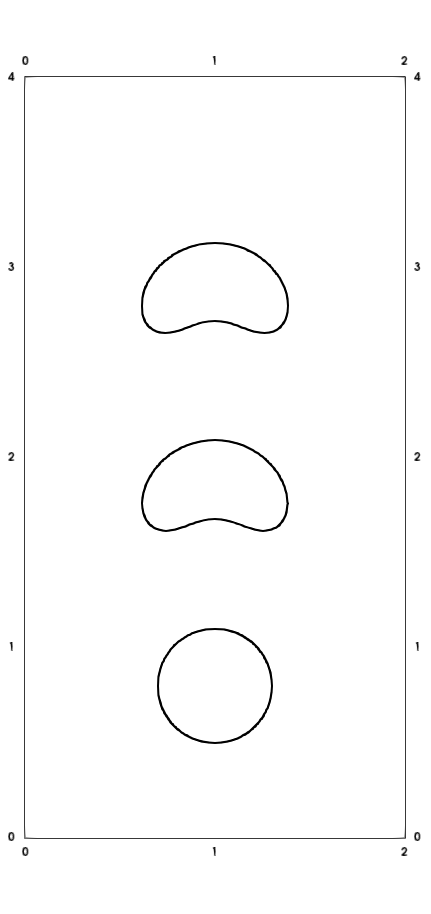}
\caption{\cblue{Time evolution of the rising bubble for the case $c_0=19.5$ with relaxation parameter $\lambda_+=0.05$ for different values of the viscoelastic diffusion parameter $\alpha\geq0$. From left to right: interface for $\alpha=0$ at times $t\in\{0,0.5,1\}$, for $\alpha=0.01$ at $t\in\{0,0.5,1\}$, for $\alpha=1$ at $t\in\{0,0.25,0.5\}$, and for $\alpha=10$ at $t\in\{0,0.15,0.3\}$.}}
\label{fig:test94_simulation_alpha}
\end{figure}

In Figure \ref{fig:test94_simulation}, we show the time evolution of the rising bubble at times $t\in\{0, 0.5, 1\}$ for the cases $c_0 \in \{0, 1, 19.5\}$ and with the relaxation parameter $\lambda_+ = 0.05$. 
In the first scenario in Figure \ref{fig:test94_simulation}, corresponding to the case $c_0 = 0$, both the bubble and the surrounding fluid are Newtonian with a viscosity ratio of $\frac{\mu_+}{\mu_+ + G \lambda_+} = 1$. Due to the high surface tension ($\gamma = 10$), the bubble maintains a shape very close to its initial spherical form.
In the second scenario in Figure \ref{fig:test94_simulation}, representing the case $c_0 = 1$, the surrounding fluid is viscoelastic, characterised by a viscosity ratio of $\frac{\mu_+}{\mu_+ + G \lambda_+} = 0.5$. In this scenario, notable differences from the fully Newtonian case are observed. Specifically, a small tail forms at the bottom of the bubble, indicating a departure from the spherical shape because of the viscoelastic nature of the surrounding fluid and its interaction with the surface tension.
In the third scenario in Figure \ref{fig:test94_simulation}, corresponding to the case $c_0 = 19.5$, the surrounding fluid is highly viscoelastic with a viscosity ratio of $\frac{\mu_+}{\mu_+ + G \lambda_+} \approx 0.05$. Here we observe the most significant difference from the Newtonian case. The tail at the bottom of the bubble is larger than in the case with $c_0 = 1$. This demonstrates the influence of the highly viscoelastic properties of the surrounding fluid, which are capable of overcoming the high surface tension and thereby changing the shape of the bubble more significantly.
\cblue{We note that significant changes in the interface evolution can be observed when the viscoelastic diffusion parameter $\alpha\geq 0$ increases, as shown in Figure \ref{fig:test94_simulation_alpha} for different values of $\alpha\in\{0, 0.01, 1, 10\}$. For larger $\alpha$, the elastic stresses around the bubble are less localized which is caused by the increased diffusion of $\mathbb{B}$. As a result, the elastic effects of the surrounding fluid are reduced, leading to a noticeably faster rise of the bubble as $\alpha$ increases. For this reason, we focussed on the smaller
value $\alpha=0.01$ in our previous simulations.}

\begin{figure}
\centering
\includegraphics[width=0.24\textwidth]{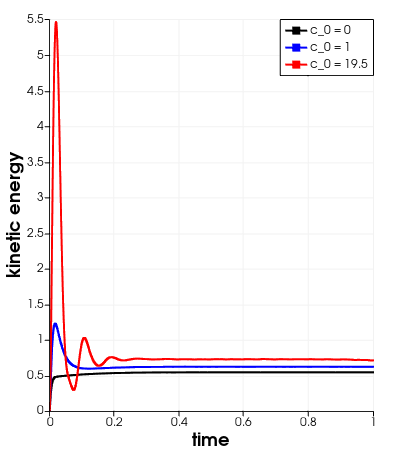}
\includegraphics[width=0.24\textwidth]{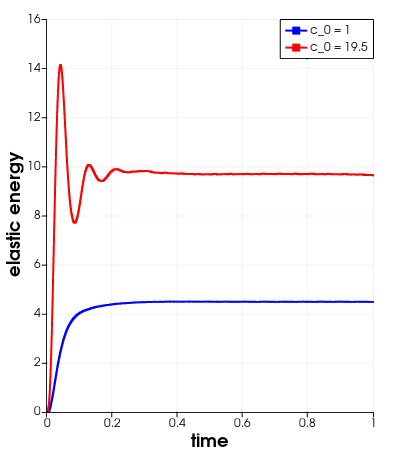}
\includegraphics[width=0.24\textwidth]{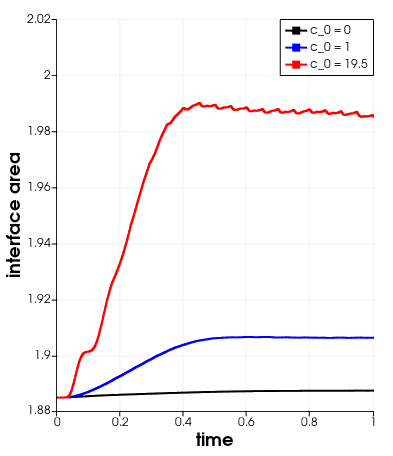}
\includegraphics[width=0.24\textwidth]{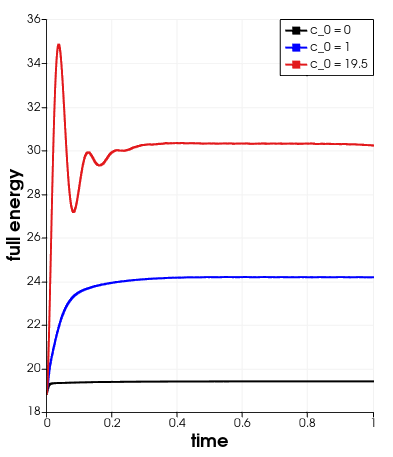}
\caption{Time evolution of the energy contributions for the cases $c_0 = 0$ (black), $c_0 = 1$ (blue) and $c_0 = 19.5$ (red) with the relaxation parameter $\lambda_+ = 0.05$. From left to right: kinetic energy, elastic energy, interfacial length and full energy. Note that the elastic energy is zero for the case $c_0 = 0$ and therefore not plotted. \cmagenta{We also point out that the full energy is not monotonically decreasing due to the presence of a nonzero body acceleration.}}
\label{fig:test94_energy}
\end{figure}

\begin{figure}
\centering
\includegraphics[width=0.24\textwidth]{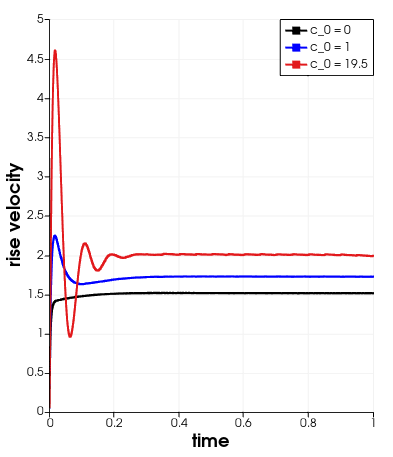}
\includegraphics[width=0.24\textwidth]{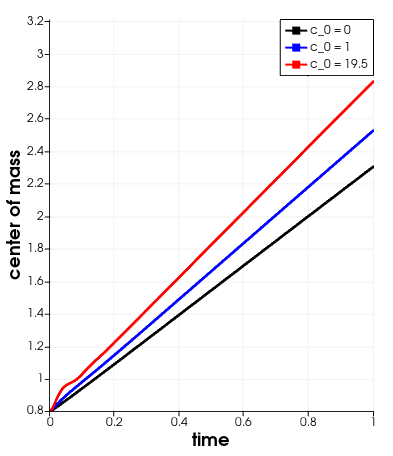}
\includegraphics[width=0.24\textwidth]{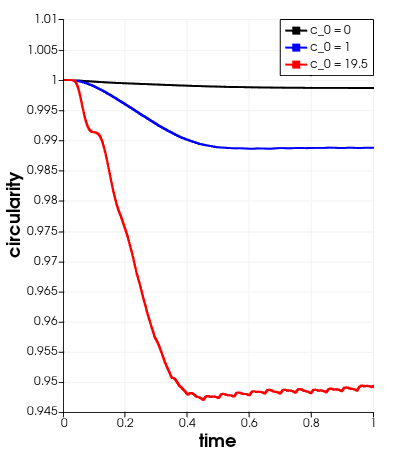}
\caption{Time evolution of some benchmark quantities for the cases $c_0 = 0$ (black), $c_0 = 1$ (blue) and $c_0 = 19.5$ (red) with the relaxation parameter $\lambda_+ = 0.05$. From left to right: rise velocity $V_c^n$, the $x_2$ component of the center of mass $y_c^n$ and the circularity $\mathrm{circ}^n$.}
\label{fig:test94_benchmark}
\end{figure}

In Figure \ref{fig:test94_energy}, we visualise the individual contributions to the energy (kinetic energy, elastic energy, interfacial length) and the total energy for the three cases $c_0 \in \{0, 1, 19.5\}$ with the relaxation parameter $\lambda_+ = 0.05$. \cmagenta{We note that the full energy is not monotonically decreasing due to the presence of a nonzero body acceleration.}
For the Newtonian case ($c_0 = 0$), all energy contributions stabilise to constant values shortly after the initial phase. We note that the elastic energy is zero for this case and is therefore not shown in the plot. This behaviour agrees with the constant rise velocity of the bubble and the stationary shape, as indicated by the circularity being very close to 1 (see Figure \ref{fig:test94_benchmark}) and the interfacial length remaining constant over time (Figure \ref{fig:test94_energy}).
Here, the rise velocity of the bubble at time step $n$ is computed as
\begin{align*}
    V_c^n = \frac{\skp{\rho_-^n \bu^n}{\mathbf{e}_2}}{\skp{\rho_-^n}{1}},
\end{align*}
where $\rho_-^n \in \calS_0^n$ is defined as in \eqref{eq:def_density} but with $\rho_+$ replaced by zero, and the circularity of the interface $\Gamma^n$ is defined as the ratio of the perimeter of an area-equivalent circle and $\calH^1(\Gamma^n)$, i.e.,
\begin{align*}
    \mathrm{circ}^n =  \frac{2 [ \pi \mathcal{L}^2(\Omega_-^n)]^{1/2}}{\calH^1(\Gamma^n)}.
\end{align*}

\begin{figure}
\centering
\includegraphics[width=0.4\textwidth]{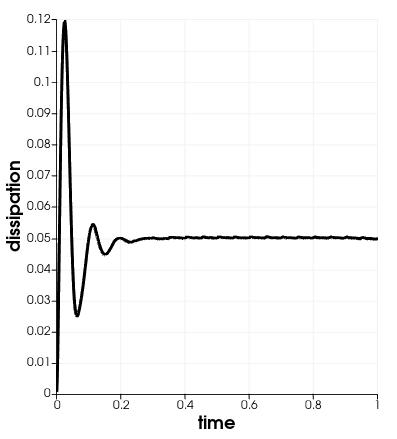}
\caption{\cblue{Time evolution of the quantity $\mathcal{D}^{n+1}$ defined in \eqref{eq:dissipation} for the case $c_0 = 19.5$, illustrating the dissipative effects in \eqref{eq:stability_FE}.}}
\label{fig:test94_dissipation}
\end{figure}

In the case $c_0 = 1$, the rise velocity of the bubble reaches its peak shortly after the initial phase, then decreases, and finally stabilises at a constant rise velocity (Figure \ref{fig:test94_benchmark}). This is reflected in the kinetic energy, which follows a similar pattern. During the same time, the elastic energy increases as the bubble decelerates. As observed in Figure \ref{fig:test94_simulation}, the bubble develops a small tail over time, which influences the interfacial length (Figure \ref{fig:test94_energy}) and the circularity (Figure \ref{fig:test94_benchmark}). After some time, the bubble reaches a stationary shape, it rises with constant velocity, and all energy contributions become constant.
\cblue{In Figure \ref{fig:test94_dissipation}, we visualize the quantity $\mathcal{D}^{n+1}$ defined in \eqref{eq:dissipation} in order to examine the magnitude of the dissipative terms in \eqref{eq:stability_FE} for the numerical experiment with $c_0=19.5$. As expected, the term remains nonnegative throughout.}

In the case $c_0 = 19.5$, we observe the most significant deviation from the Newtonian case. In the beginning, there is an oscillatory behaviour in the rise velocity, as shown in Figure \ref{fig:test94_benchmark}. The rise velocity initially peaks at a value much higher than in the other cases, followed by a few oscillations before it reaches a constant rise velocity. 
The kinetic energy plot in Figure \ref{fig:test94_energy} mirrors this trend. The elastic energy also shows oscillations, but with a slight shift in time. As the bubble develops a tail, the interfacial length increases until it stabilises at a nearly constant value with small oscillations in the bottom tip. These oscillations, which are visible in the plots of the interface length and circularity, persist even with finer discretization parameters and are influenced by the relaxation parameters $\lambda_\pm$.
Compared to the other cases, the final constant rise velocity in the case $c_0=19.5$ is higher and also the plots for the bubble's center of mass for the three cases differ, see Figure \ref{fig:test94_benchmark}, where (the $x_2$-component of) the bubble's center of mass at time step $n$ is defined by 
\begin{align*}
    y_c^n = \frac{1}{\mathcal{L}^2(\Omega_-^n)} \int_{\Omega_-^n} x_2 \dx.
\end{align*}
However, this behaviour is in agreement with other works which investigate the effect of viscoelasticity on the rise velocity of the bubble \cite{pillapakkam_2007_viscoelastic}.

\begin{figure}
\centering
\includegraphics[height=5cm]{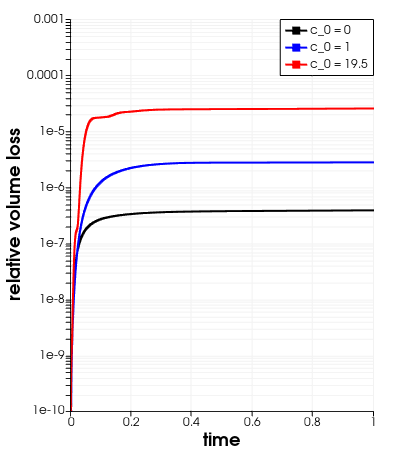}
\includegraphics[height=5cm]{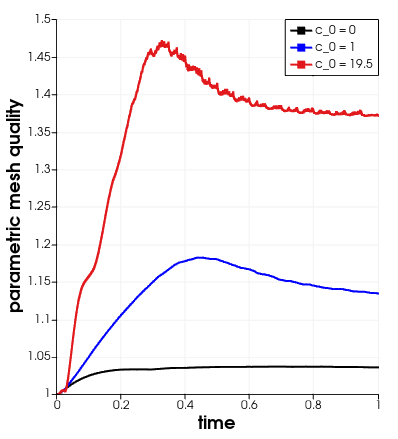}
\caption{Time evolution of the relative volume loss $\mathcal{L}_\mathrm{loss}^n$ in a logarithmic plot, and the element ratios $r^n$ as a measure for the parametric mesh quality for the three cases $c_0\in\{0, 1, 19.5\}$ with the relaxation parameter $\lambda_+ = 0.05$.}
\label{fig:test94_benchmark2}
\end{figure}

We note that in all three cases $c_0\in\{0, 1, 19.5\}$, we observe excellent volume conservation properties. 
In Figure \ref{fig:test94_benchmark2}, we visualize the relative volume loss from the initial configuration to the time step $n$, which is defined as
\begin{align*}
    \mathcal{L}_\mathrm{loss}^n = \frac{\mathcal{L}^2(\Omega_-^0) - \mathcal{L}^2(\Omega_-^n)}{\mathcal{L}^2(\Omega_-^0)},
\end{align*}
measuring the decrease in area of the inner phase $\Omega_-^n$ relative to its initial size. Remarkably, the relative volume loss was much lower than $0.01\%$ in all computations. In addition, the mesh quality of the interface mesh was also very high. 
Here we use the element ratios as a measure,
\begin{align*}
    r^n = \frac{\max_{\sigma \in \calT(\Gamma^n)} \calH^{1}(\sigma)}{\min_{\sigma \in \calT(\Gamma^n)} \calH^{1}(\sigma)},
\end{align*}
which consistently remained much smaller than $2$, indicating a good distribution of the vertices on the interface mesh.

\begin{figure}
\centering
\includegraphics[width=0.32\textwidth,trim={0 2cm 0 2cm},clip]{figures/test92_evolution.png}
\includegraphics[width=0.32\textwidth,trim={0 2cm 0 2cm},clip]{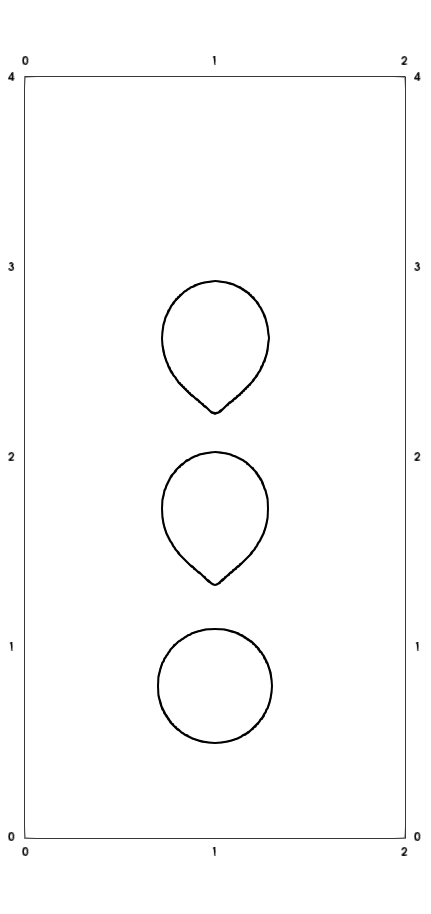}
\includegraphics[width=0.32\textwidth,trim={0 2cm 0 2cm},clip]{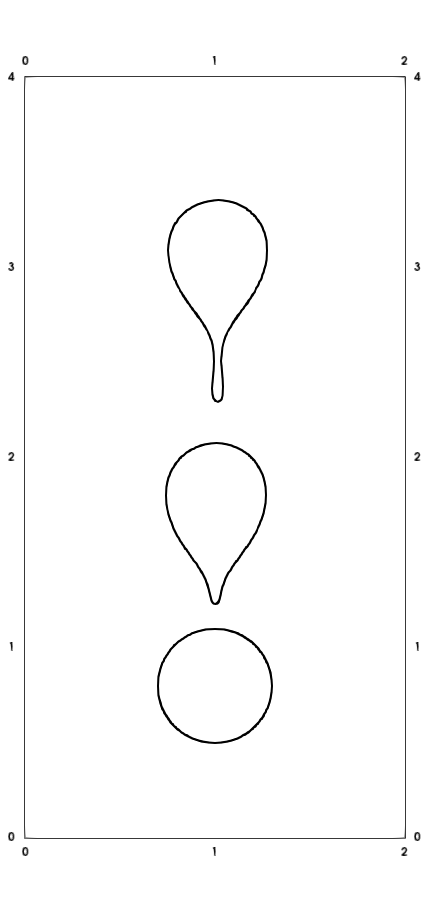}
\caption{Time evolution of the rising bubble with $\lambda_+=0.075$ for the cases $c_0=0$ (left) and $c_0=1$ (center) at times $t\in\{0, 0.5, 1\}$, and for $c_0=19.5$ (right) at times $t\in\{0, 0.4, 1\}$.}
\label{fig:test95_simulation}
\end{figure}

\begin{figure}
\centering
\includegraphics[width=0.24\textwidth]{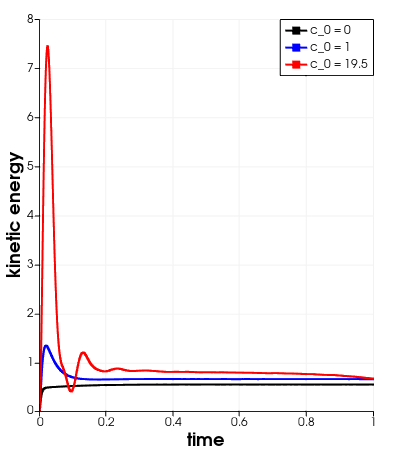}
\includegraphics[width=0.24\textwidth]{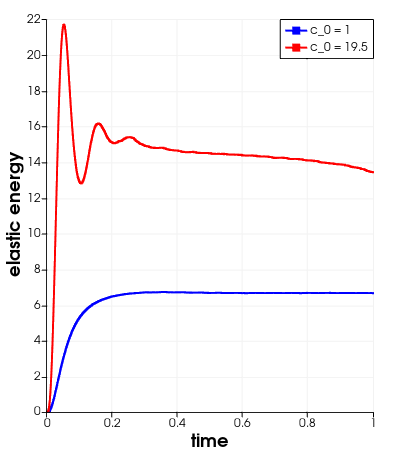}
\includegraphics[width=0.24\textwidth]{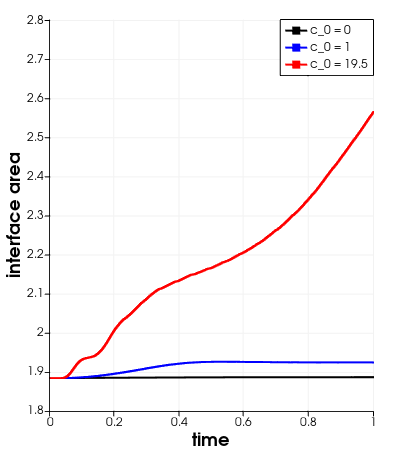}
\includegraphics[width=0.24\textwidth]{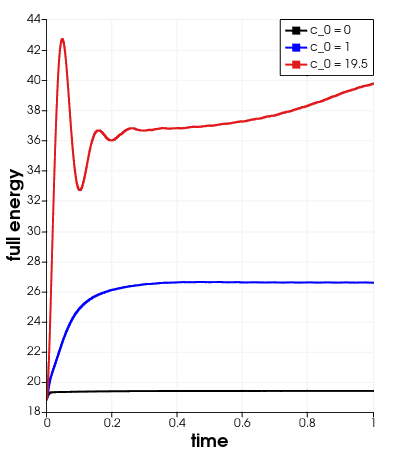}
\caption{Time evolution of the energy contributions for the cases $c_0 = 0$ (black), $c_0 = 1$ (blue) and $c_0 = 19.5$ (red) with the relaxation parameter $\lambda_+ = 0.075$. From left to right: kinetic energy, elastic energy, interfacial length and full energy. Note that the elastic energy is zero for the case $c_0 = 0$ and therefore not plotted. \cmagenta{We also point out that the full energy is not monotonically decreasing due to the presence of a nonzero body acceleration.}}
\label{fig:test95_energy}
\end{figure}

\begin{figure}
\centering
\includegraphics[width=0.24\textwidth]{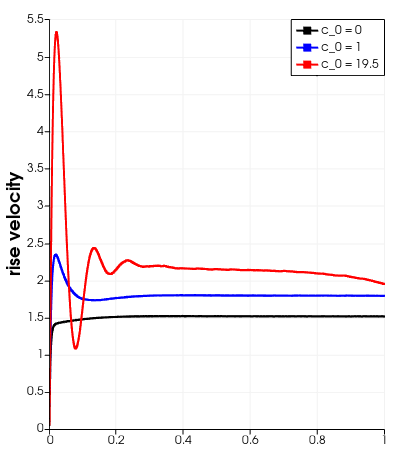}
\includegraphics[width=0.24\textwidth]{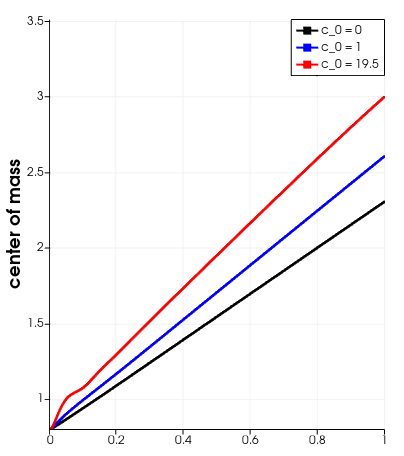}
\includegraphics[width=0.24\textwidth]{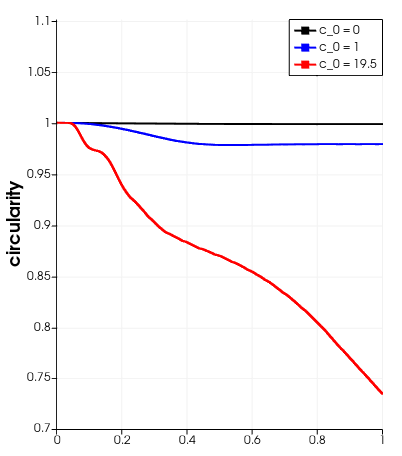}
\caption{Time evolution of some benchmark quantities for the cases $c_0 = 0$ (black), $c_0 = 1$ (blue) and $c_0 = 19.5$ (red) with the relaxation parameter $\lambda_+ = 0.075$. From left to right: rise velocity $V_c^n$, the $x_2$ component of the center of mass $y_c^n$ and the circularity $\mathrm{circ}^n$.}
\label{fig:test95_benchmark}
\end{figure}

\begin{figure}
\centering
\includegraphics[height=5cm]{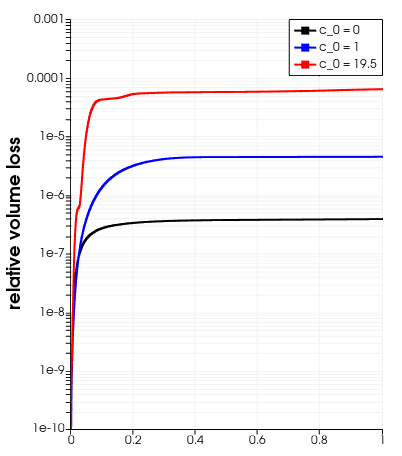}
\includegraphics[height=5cm]{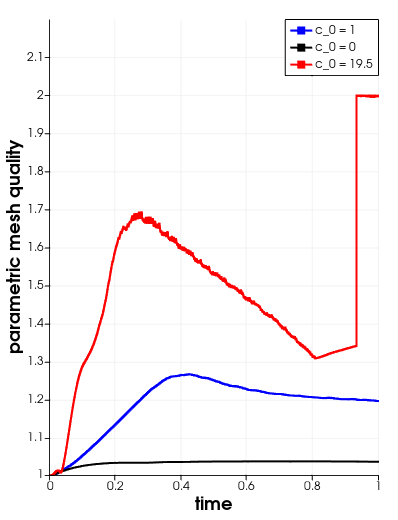}
\caption{Time evolution of the relative volume loss $\mathcal{L}_\mathrm{loss}^n$ in a logarithmic plot, and the element ratios $r^n$ as a measure for the parametric mesh quality for the three cases $c_0\in\{0, 1, 19.5\}$ with the relaxation parameter $\lambda_+ = 0.075$. }
\label{fig:test95_benchmark2}
\end{figure}


In the following, we present the results for the computations with the relaxation parameter $\lambda_+=0.075$. In Figure \ref{fig:test95_simulation}, the time evolution of the interface is shown for the values $c_0\in\{0,1\}$ at times $t\in\{0,0.5,1\}$ and for $c_0=19.5$ at times $t\in\{0,0.4,1\}$. We note that the Newtonian case with $c_0=0$ is identical in both experiments.
For the case $c_0=1$, the shapes of the bubble in Figure \ref{fig:test95_simulation} are comparable to those in Figure \ref{fig:test94_simulation}, where a smaller relaxation parameter was used. This similarity is due to the influence of viscoelasticity being quite small, as indicated by the viscosity ratio $\frac{\mu_+}{\mu_+ + G \lambda_+} = 0.5$. In the case $c_0=19.5$, where the viscosity ratio is $\frac{\mu_+}{\mu_+ + G \lambda_+} \approx 0.05$, the influence of viscoelasticity is much stronger, and we observe more differences compared to before. In particular, the final shape of the bubble in Figure \ref{fig:test95_simulation} for $c_0=19.5$ is more elongated, and the tail at the bottom tip is longer and thinner than for the same value of $c_0$ in Figure \ref{fig:test94_simulation}.

The influence of the larger relaxation parameter is also observed in the energy plots in Figure \ref{fig:test95_energy}. For $c_0=1$, the energy plots resemble those from Figure \ref{fig:test94_energy}, with larger peaks and oscillations in the beginning for the kinetic and elastic energies for the case $c_0=19.5$. Moreover, in the case $c_0=19.5$, as the tail of the bubble evolves, the interface length also increases over time and does not attain a constant value. 
\cmagenta{We note that the full energy is not monotonically decreasing due to the presence of a nonzero body acceleration.}
Similarly, in Figure \ref{fig:test95_benchmark}, the circularity for $c_0=19.5$ mirrors this trend by decreasing its values.
As before, the rise velocity is larger in the case $c_0=19.5$ than for the values $c_0\in\{0, 1\}$, which is reflected in the plot of (the $x_2$-component of) the bubble's center of mass.

As before, we observe excellent volume conservation properties with a relative volume loss much smaller than $0.01\%$ and also a good quality of the parametric interface mesh, as seen in Figure \ref{fig:test95_benchmark2}. We note that the jump in the plot of the element ratios $r^n$ for the case $c_0=19.5$ in Figure \ref{fig:test95_benchmark2} is due to the local refinement of the interface mesh, as some elements have grown too large as the bottom end of the bubble has elongated.

\subsection{Interpretation of the viscoelastic behaviour}

\begin{figure}
\centering
\includegraphics[width=0.32\textwidth]{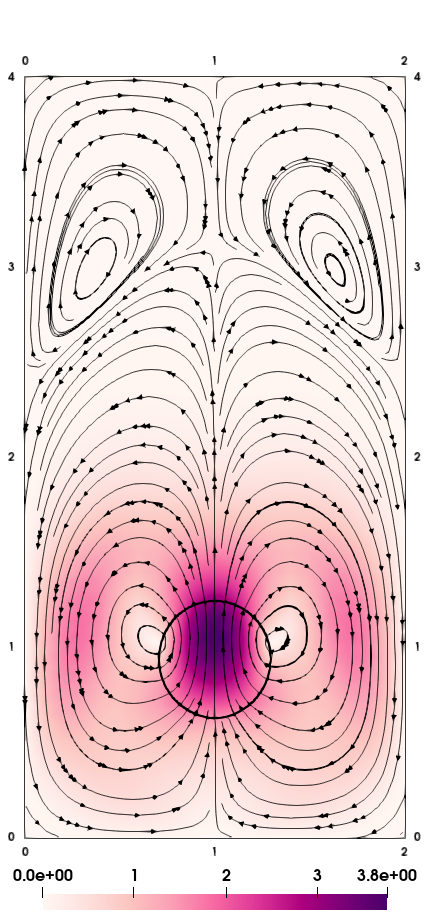}
\includegraphics[width=0.32\textwidth]{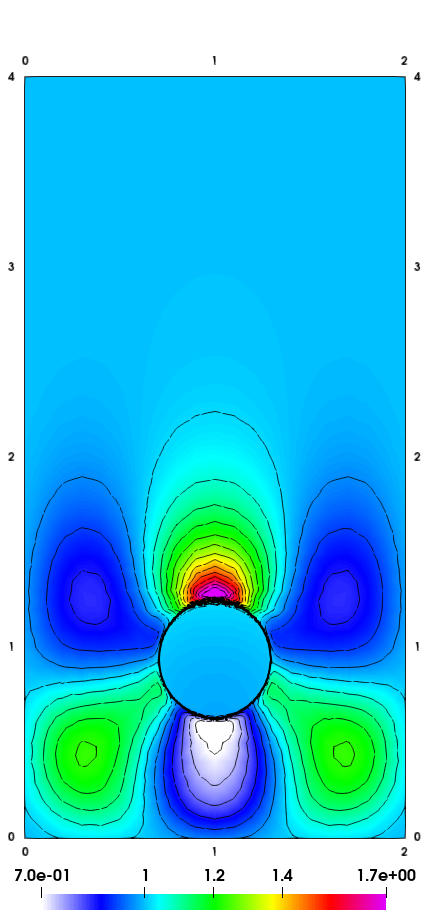}
\includegraphics[width=0.32\textwidth]{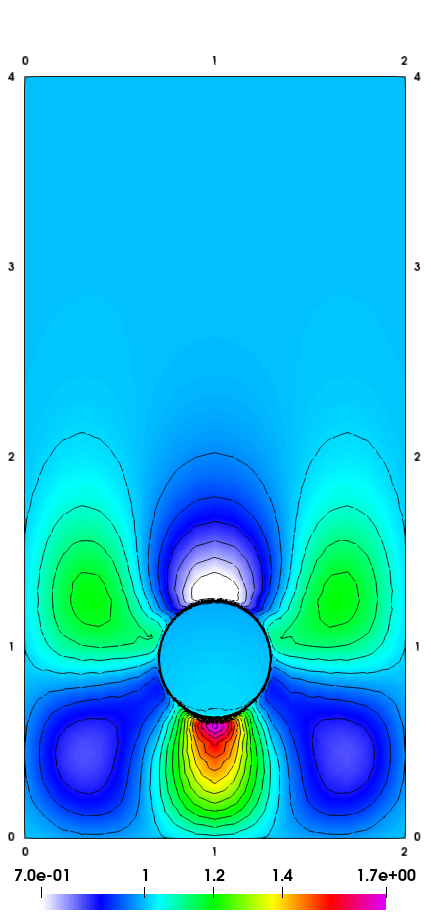}
\caption{Solution at time $t=0.04$ for the case $c_0 = 19.5$ with the relaxation parameter $\lambda_+ = 0.05$. Left: magnitude and streamlines of the velocity $\bu^n$. Center and right: diagonal components $\bbB^n_{11}$ and $\bbB^n_{22}$, respectively, of the tensor $\bbB^n$.}
\label{fig:test94_steamlines_isovalues_t=0.04}
\end{figure}

\begin{figure}
\centering
\includegraphics[width=0.32\textwidth]{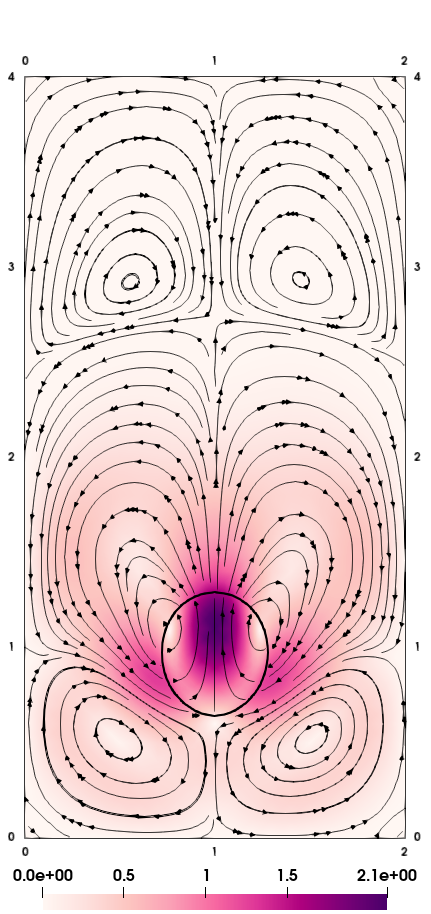}
\includegraphics[width=0.32\textwidth]{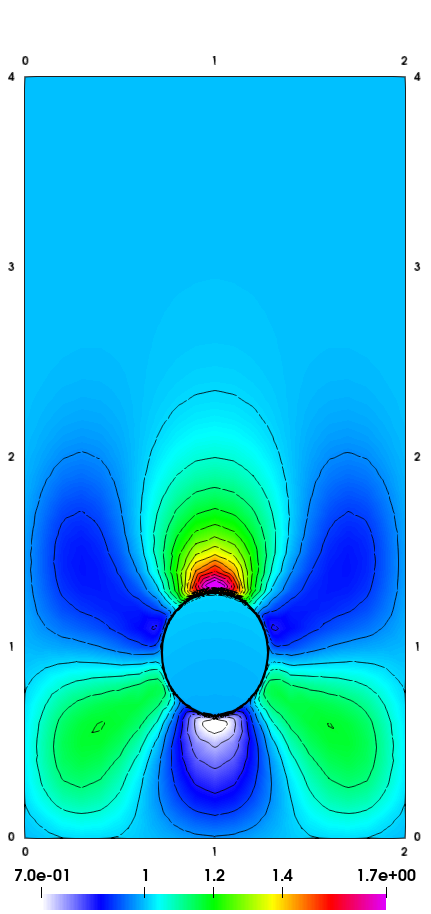}
\includegraphics[width=0.32\textwidth]{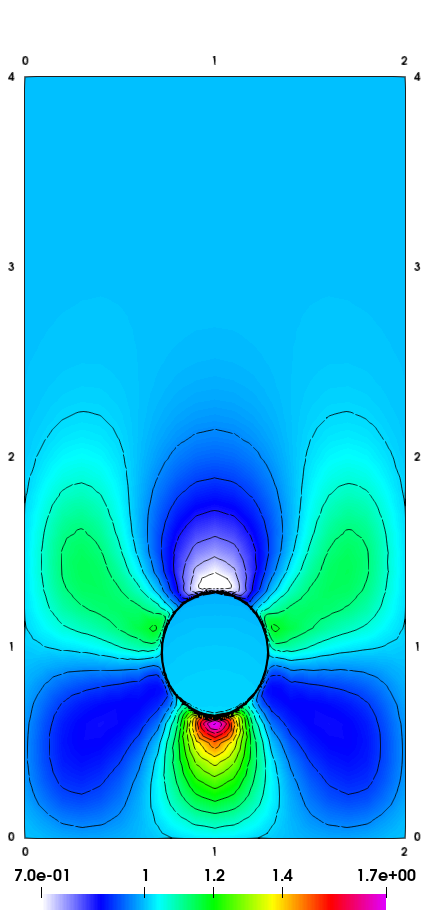}
\caption{Solution at time $t=0.06$ for the case $c_0 = 19.5$ with the relaxation parameter $\lambda_+ = 0.05$. Left: magnitude and streamlines of the velocity $\bu^n$. Center and right: diagonal components $\bbB^n_{11}$ and $\bbB^n_{22}$, respectively, of the tensor $\bbB^n$.
}
\label{fig:test94_steamlines_isovalues_t=0.06}
\end{figure}

\begin{figure}
\centering
\includegraphics[width=0.32\textwidth]{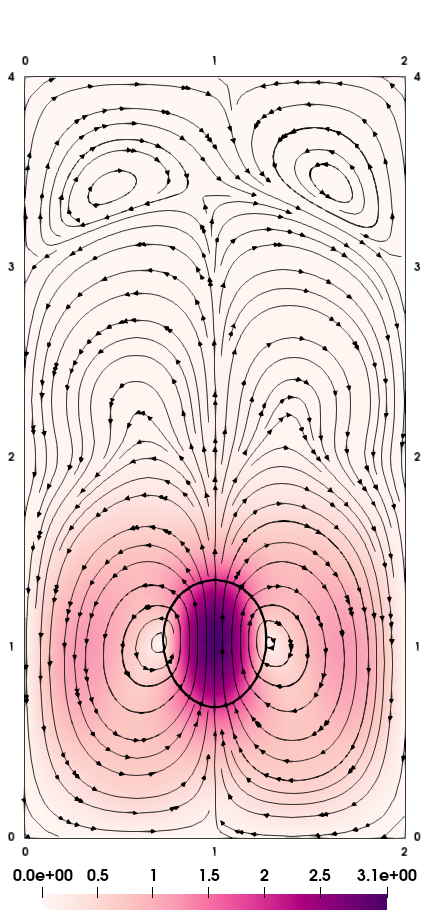}
\includegraphics[width=0.32\textwidth]{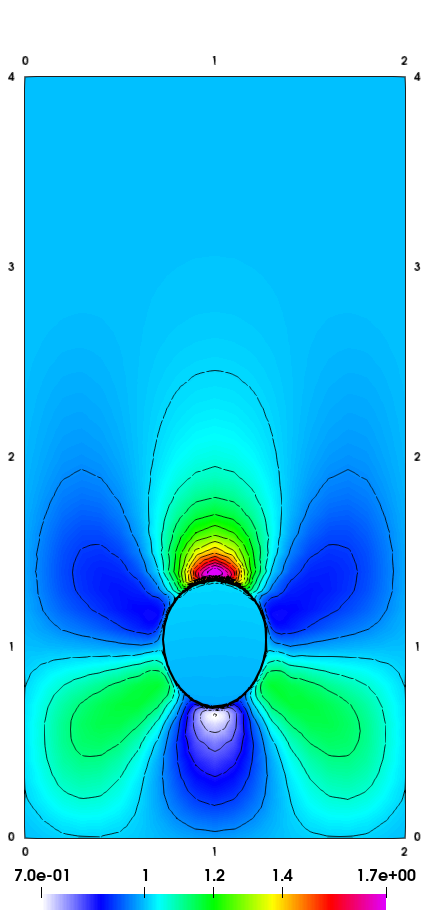}
\includegraphics[width=0.32\textwidth]{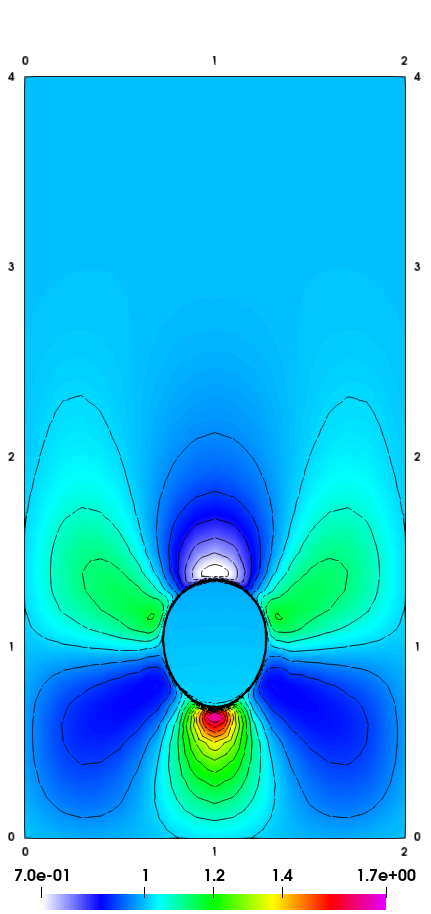}
\caption{Solution at time $t=0.1$ for the case $c_0 = 19.5$ with the relaxation parameter $\lambda_+ = 0.05$. Left: magnitude and streamlines of the velocity $\bu^n$. Center and right: diagonal components $\bbB^n_{11}$ and $\bbB^n_{22}$, respectively, of the tensor $\bbB^n$.}
\label{fig:test94_steamlines_isovalues_t=0.10}
\end{figure}

\begin{figure}
\centering
\includegraphics[width=0.32\textwidth]{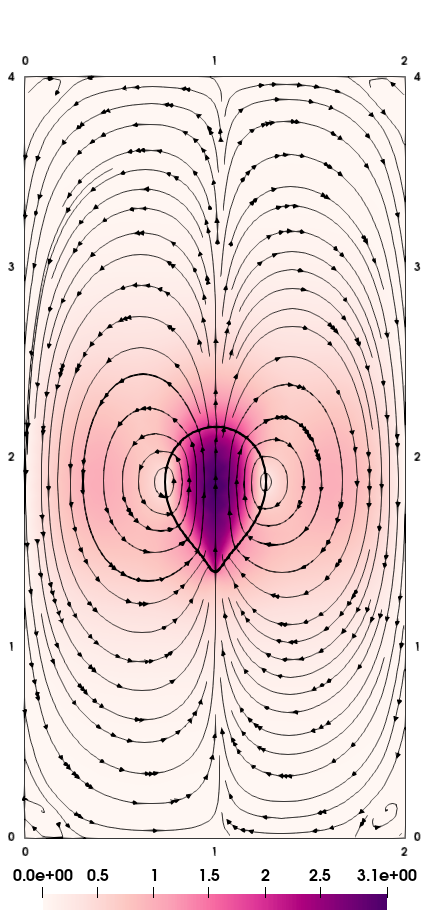}
\includegraphics[width=0.32\textwidth]{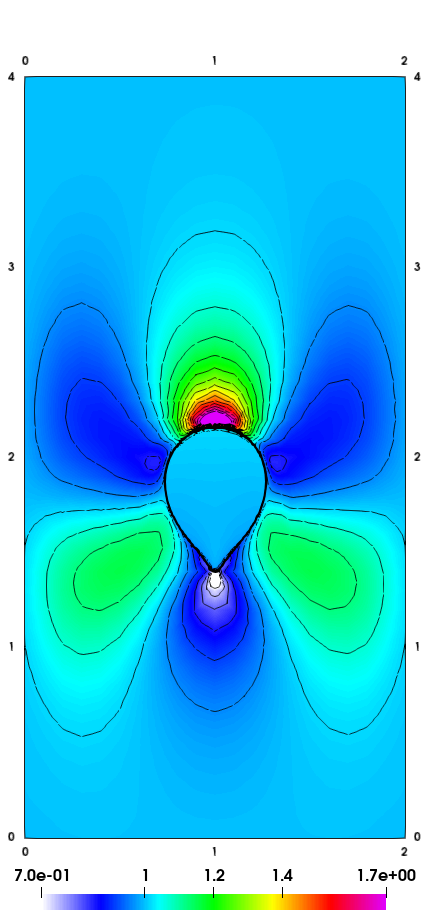}
\includegraphics[width=0.32\textwidth]{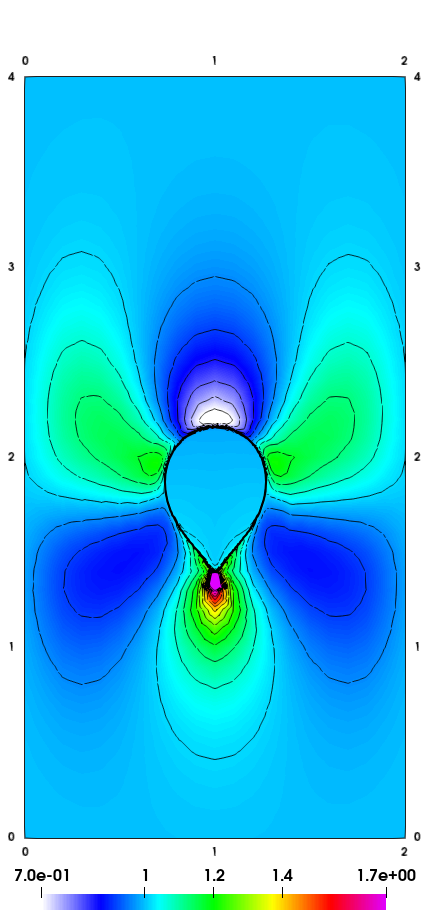}
\caption{Solution at time $t=0.5$ for the case $c_0 = 19.5$ with the relaxation parameter $\lambda_+ = 0.05$. Left: magnitude and streamlines of the velocity $\bu^n$. Center and right: diagonal components $\bbB^n_{11}$ and $\bbB^n_{22}$, respectively, of the tensor $\bbB^n$.}
\label{fig:test94_steamlines_isovalues_t=0.50}
\end{figure}

\begin{figure}
\centering
\includegraphics[width=0.32\textwidth]{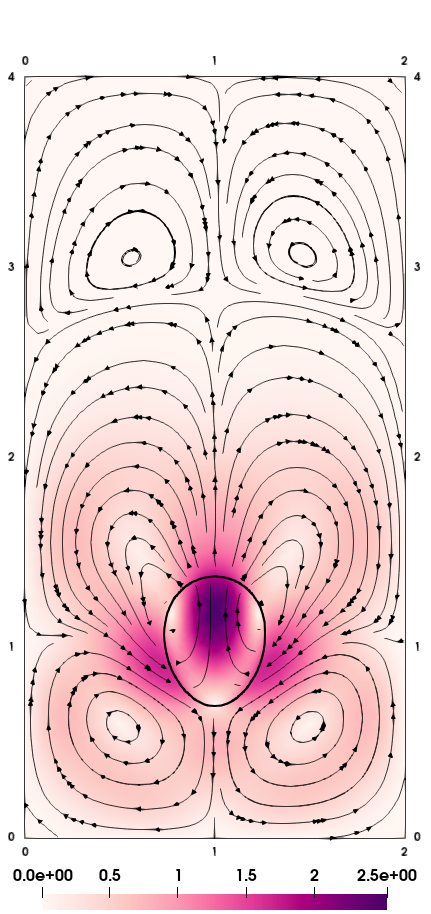}
\includegraphics[width=0.32\textwidth]{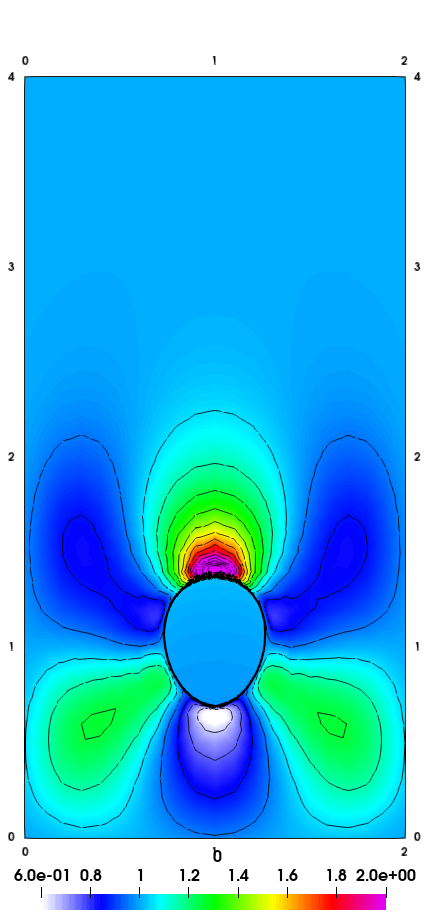}
\includegraphics[width=0.32\textwidth]{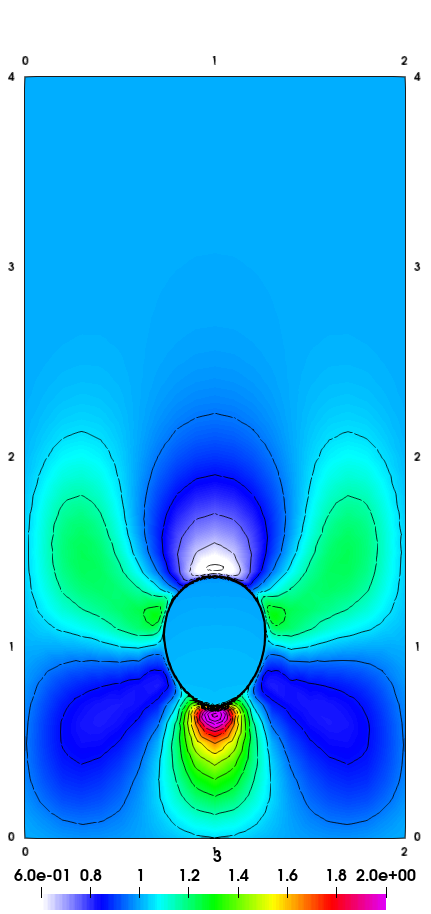}
\caption{Solution at time $t=0.075$ for the case $c_0 = 19.5$ with the relaxation parameter $\lambda_+ = 0.075$. Left: magnitude and streamlines of the velocity $\bu^n$. Center and right: diagonal components $\bbB^n_{11}$ and $\bbB^n_{22}$, respectively, of the tensor $\bbB^n$.
}
\label{fig:test95_steamlines_isovalues_t=0.075}
\end{figure}

\begin{figure}
\centering
\includegraphics[width=0.32\textwidth]{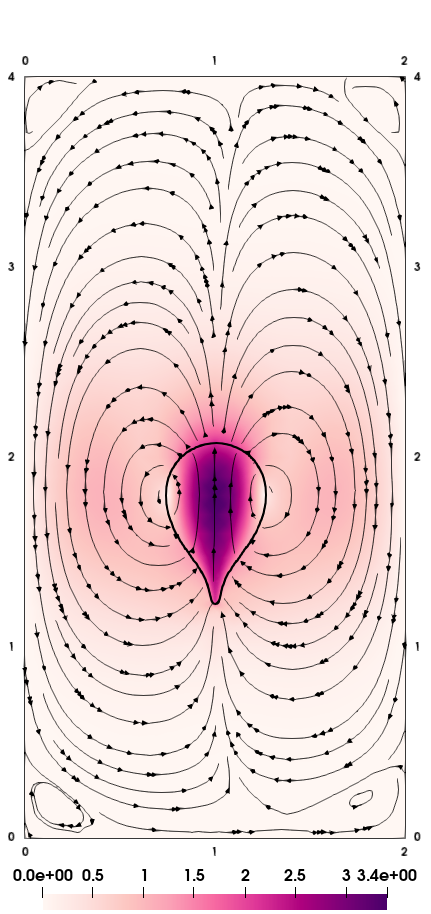}
\includegraphics[width=0.32\textwidth]{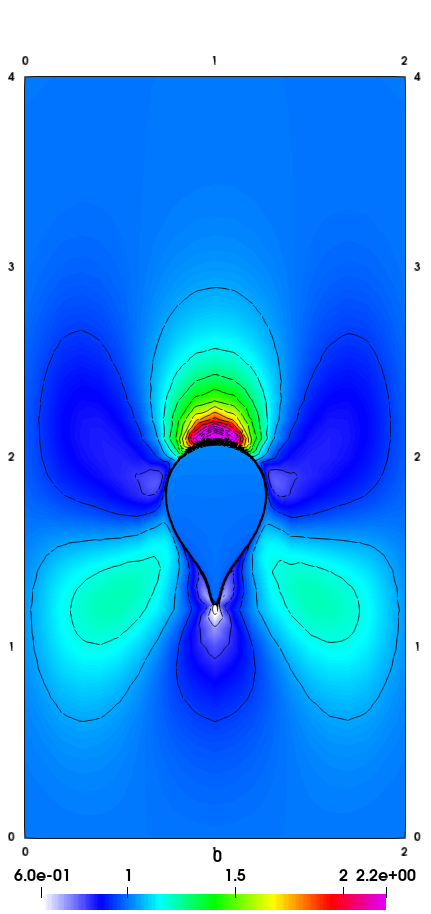}
\includegraphics[width=0.32\textwidth]{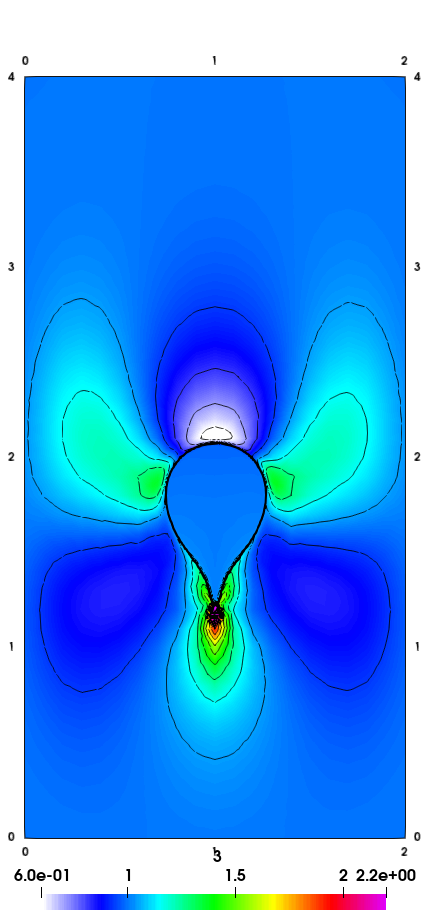}
\caption{Solution at time $t=0.4$ for the case $c_0 = 19.5$ with the relaxation parameter $\lambda_+ = 0.075$. Left: magnitude and streamlines of the velocity $\bu^n$. Center and right: diagonal components $\bbB^n_{11}$ and $\bbB^n_{22}$, respectively, of the tensor $\bbB^n$.
}
\label{fig:test95_steamlines_isovalues_t=0.4}
\end{figure}

Now, we discuss the insights for the rising bubble experiments \cred{from Section \ref{sec:numerics_rising_bubble}}, which demonstrate the complex interactions between the viscoelastic properties of the surrounding fluid and the movement of the bubble.
The overall behaviour of the bubble in the highly viscoelastic case ($c_0 = 19.5$) is typical for viscoelastic materials, which exhibit elastic (oscillatory) behaviour on short time scales and viscous (damped) behaviour on longer time scales. This phenomenon for rising bubbles in a viscoelastic fluid was first observed in \cite{hassager_1979_viscoelastic} and has been verified in several numerical studies, including \cite{pillapakkam_2007_viscoelastic, phanthien_2020_rising_bubble} and references therein. 
To improve our understanding, we visualise the velocity $\bu^n$ and the diagonal components of the tensor $\bbB^n$ at four different times in Figures \ref{fig:test94_steamlines_isovalues_t=0.04}--\ref{fig:test94_steamlines_isovalues_t=0.50} for the case $\lambda_+=0.05$, and at two different times in Figures \ref{fig:test95_steamlines_isovalues_t=0.075} and \ref{fig:test95_steamlines_isovalues_t=0.4} for the case $\lambda_+=0.075$. In Figure \ref{fig:test94_steamlines_isovalues_t=0.04}, we present the solution at $t=0.04$ with the relaxation parameter $\lambda_+=0.05$, when the rise velocity and kinetic energy of the bubble are at their maximum (see Figures \ref{fig:test94_energy} and \ref{fig:test94_benchmark}). The plots of the diagonal components of $\bbB^n$ reveal significant deviations from the identity matrix above and below the bubble, indicating high local elastic stresses in these regions. Recall that the elastic stress tensor for the Oldroyd-B model is defined as $\bbT_e = G (\bbB - \bbI)$.
In the next plot, Figure \ref{fig:test94_steamlines_isovalues_t=0.06}, we show the solution at $t=0.06$, when the bubble starts to decelerate and the elastic energy increases. At this moment, two additional vortices appear in the velocity plot, and a negative wake forms right below the bubble, causing it to slow down and stretch vertically. This is due to the high elastic stresses that have built up below the bubble.
Next, in Figure \ref{fig:test94_steamlines_isovalues_t=0.10}, we show the bubble at $t=0.1$, when it begins to accelerate again, and the elastic energy reaches a local minimum. Here, the two additional vortices have disappeared. In the following moments, the rise velocity of the bubble continues to oscillate slightly, but with a much smaller amplitude until it reaches a constant value.
Finally, we present the solution at $t=0.5$ in Figure \ref{fig:test94_steamlines_isovalues_t=0.50}, where the elastic stresses are concentrated at the top and bottom of the bubble, acting more strongly than the surface tension forces. As a result, a tail is formed behind the rising bubble due to the high local elastic stresses.

The effects that we observed for the relaxation parameter $\lambda_+=0.05$ are even stronger for $\lambda_+=0.075$. To highlight the enhanced effects of viscoelasticity, we present the solution at two different times. As before, the bubble exhibits an oscillatory behaviour in the beginning, as previously mentioned. In Figure \ref{fig:test95_steamlines_isovalues_t=0.075}, we show the solution at time $t=0.075$, when the rise velocity in Figure \ref{fig:test95_benchmark} is at its lowest value. As before, there are two additional vortices below the bubble, causing a negative wake and stretching the bubble in the vertical direction. Compared to Figure \ref{fig:test94_steamlines_isovalues_t=0.06}, this negative wake is stronger, resulting in an even more elongated bubble. This can be explained with the larger deviations of $\bbB^n$ from the identity matrix below the bubble, as shown in Figure \ref{fig:test95_steamlines_isovalues_t=0.075}.
Next, we show the solution at time $t=0.4$ in Figure \ref{fig:test95_steamlines_isovalues_t=0.4}. At this moment, the rise velocity of the bubble (see Figure \ref{fig:test95_benchmark}) is almost constant, but the bubble does not attain a stationary shape as the tail continues to grow. This has already been observed in the plot of the interface length in Figure \ref{fig:test95_energy}. The reason for this is the very high local elastic stresses near the bottom tip of the bubble, which result in the elongated bottom tip.

\begin{figure}
\centering
\includegraphics[width=0.32\textwidth]{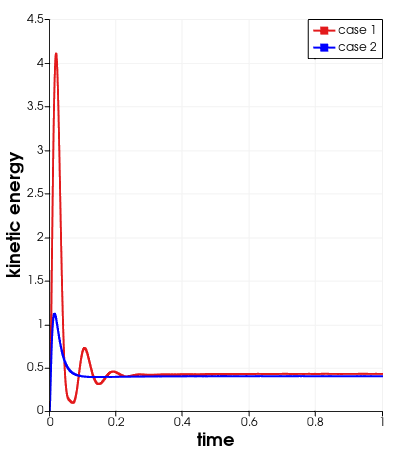}
\includegraphics[width=0.32\textwidth]{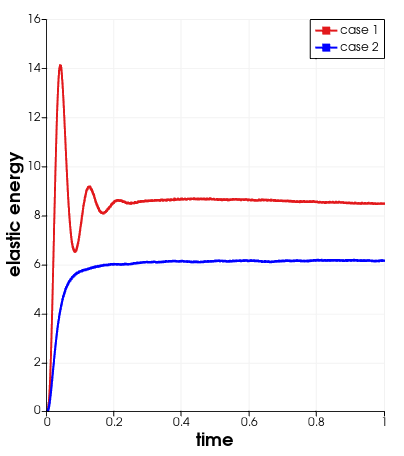}
\includegraphics[width=0.32\textwidth]{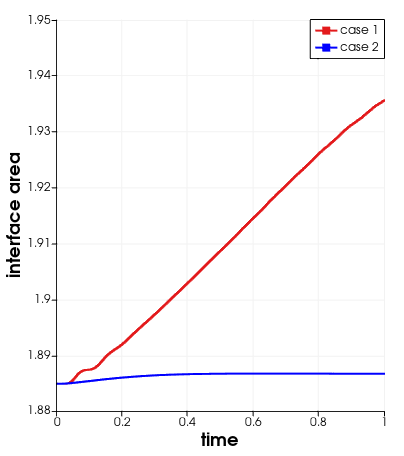}
\caption{\cblue{Time evolution of the energy contributions for the rising bubble experiment with variable shear modulus. From left to right: kinetic energy, elastic energy and interface length for the first case $c_+=19.5$, $c_-=1$ (red) and the second case $c_+=1$, $c_-=19.5$ (blue).}}
\label{fig:test94_Gvar_energy}
\end{figure}

\begin{figure}
\centering
\includegraphics[width=0.32\textwidth]{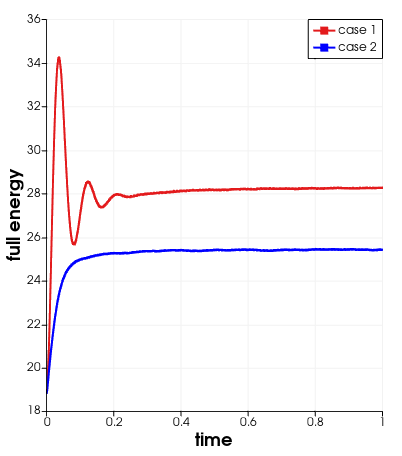}
\includegraphics[width=0.32\textwidth]{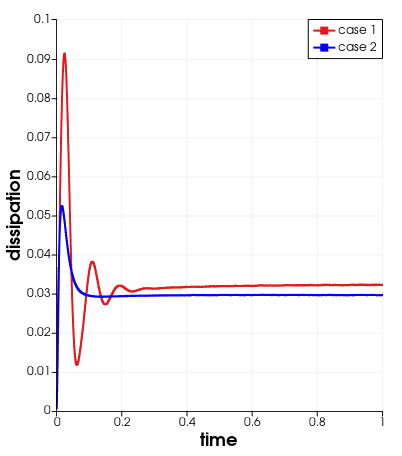}
\caption{\cblue{Time evolution of the full energy and the dissipations $\mathcal{D}^{n+1}$ defined in \eqref{eq:dissipation_Gvar} for the rising bubble experiment with variable shear modulus. The first case with $c_+=19.5$, $c_-=1$ is shown in red, and the second case $c_+=1$, $c_-=19.5$ is shown in blue.}}
\label{fig:test94_Gvar_energy_full}
\end{figure}

\cblue{\subsection{The rising bubble experiment with variable shear modulus}
\label{sec:numerics_Gvar}

\begin{figure}
\centering
\includegraphics[width=0.32\textwidth]{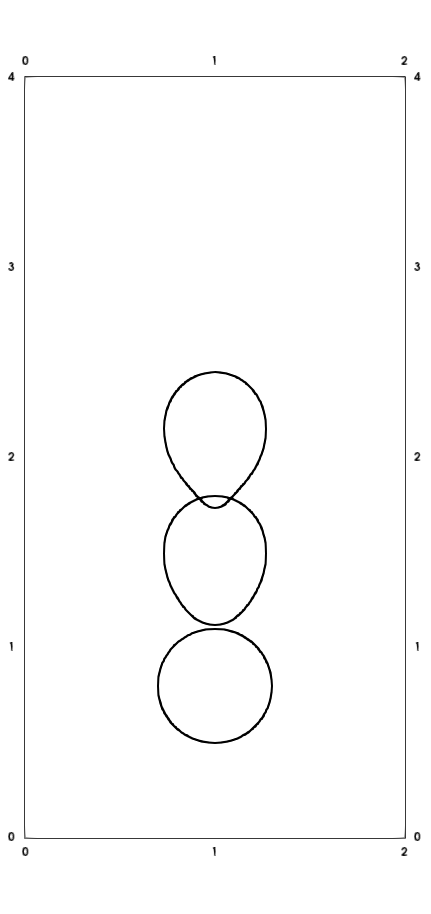}
\includegraphics[width=0.32\textwidth]{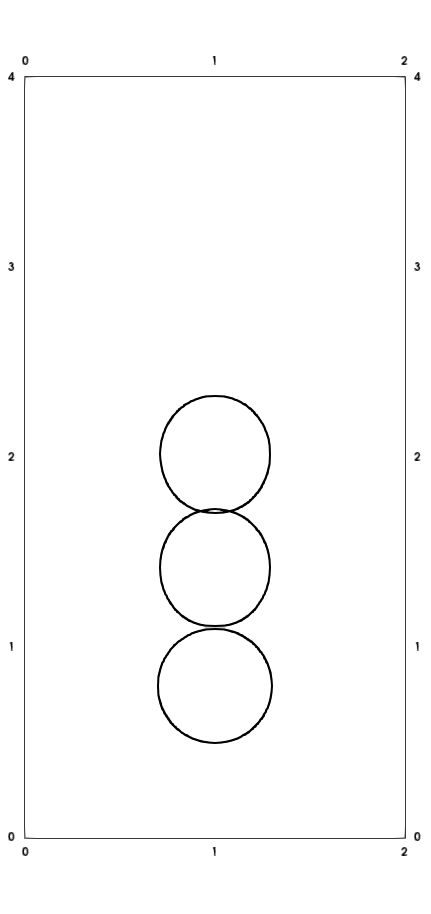}
\caption{\cblue{Time evolution of the rising bubble with variable shear modulus for the case $c_+=19.5$, $c_-=1$ (left), and for $c_+=1$, $c_-=19.5$ (right) at times $t\in\{0,0.5,1\}$.}}
\label{fig:test94_Gvar_simulation}
\end{figure}

\begin{figure}
\centering
\includegraphics[width=0.32\textwidth]{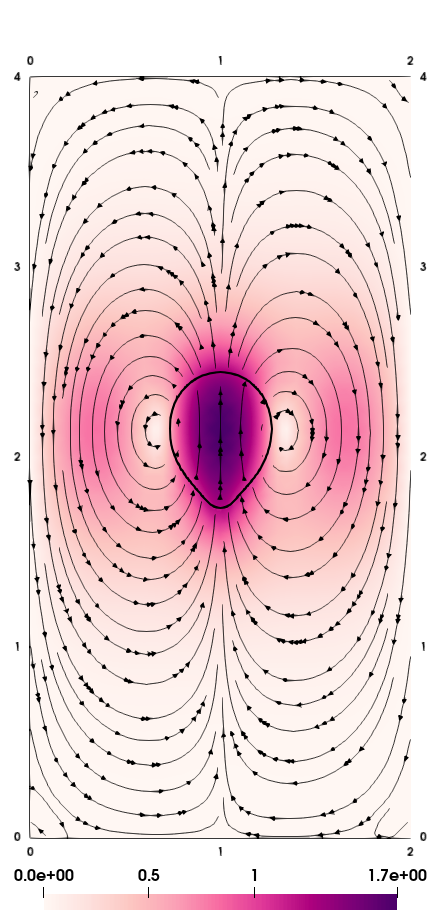}
\includegraphics[width=0.32\textwidth]{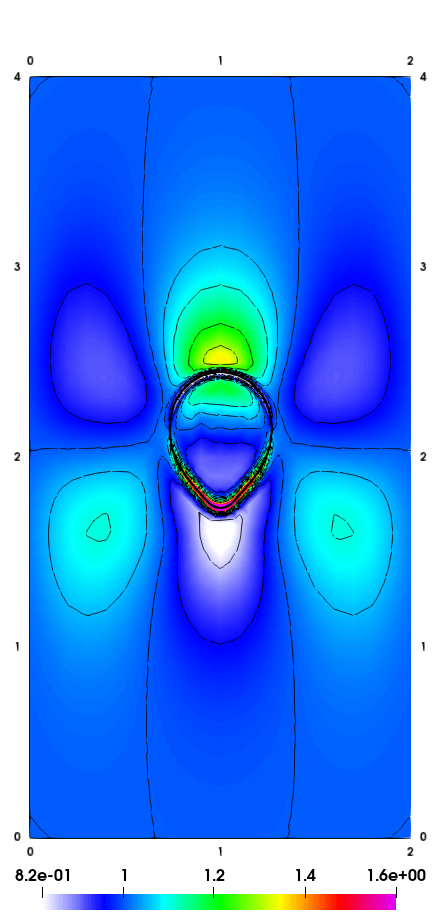}
\includegraphics[width=0.32\textwidth]{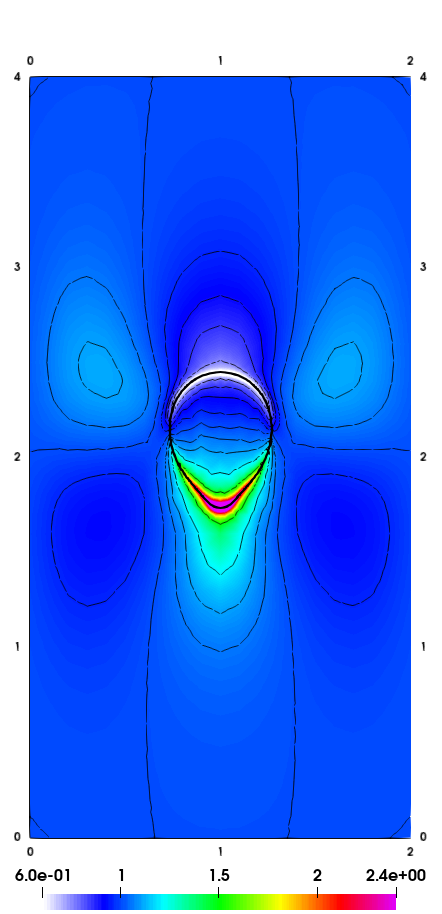}
\caption{\cblue{Solution with variable shear modulus at time $t=1$ for the first case ($c_+=19.5$, $c_-=1$). Left: magnitude and streamlines of the velocity $\bu^n$. Center and right: diagonal components $\bbB^n_{11}$ and $\bbB^n_{22}$, respectively, of the tensor $\bbB^n$.
}}
\label{fig:test94_Gvar1_steamlines_isovalues_t=1}
\end{figure}

\begin{figure}
\centering
\includegraphics[width=0.32\textwidth]{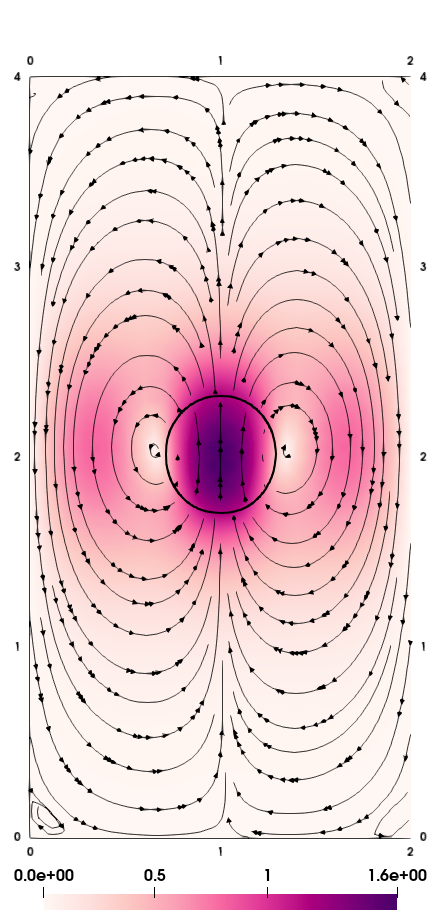}
\includegraphics[width=0.32\textwidth]{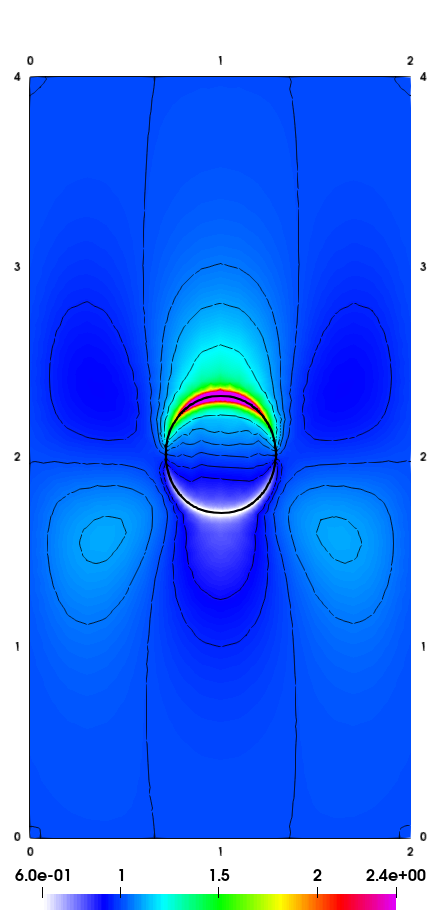}
\includegraphics[width=0.32\textwidth]{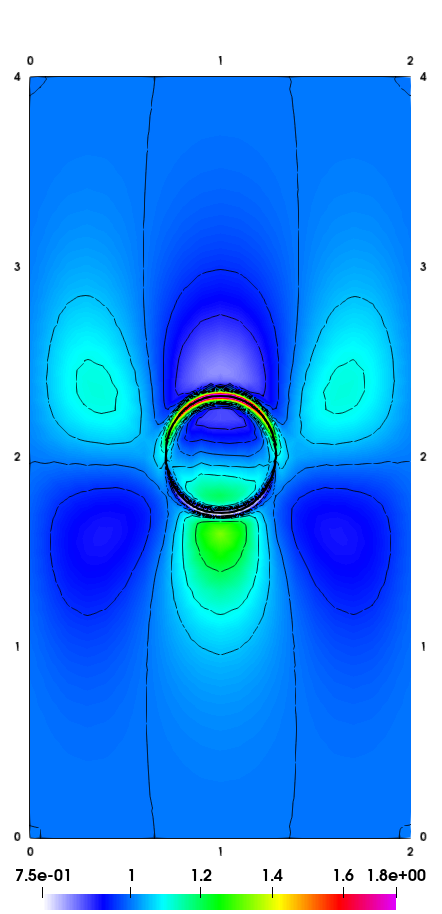}
\caption{\cblue{Solution with variable shear modulus at time $t=1$ for the second case ($c_+=1$, $c_-=19.5$). Left: magnitude and streamlines of the velocity $\bu^n$. Center and right: diagonal components $\bbB^n_{11}$ and $\bbB^n_{22}$, respectively, of the tensor $\bbB^n$.
}}
\label{fig:test94_Gvar2_steamlines_isovalues_t=1}
\end{figure}

In the following, we show two numerical results for the discrete system \eqref{eq:system_FE_G} with a variable shear modulus $G(\cdot,t) = G_\pm$ in $\Omega_\pm(t)$. The discrete setting is the same as in Section \ref{sec:numerics_rising_bubble}, but we now use the following model parameters:
\begin{align*}
    \rho_+ &= 1, \quad
    \rho_- = 0.1, \quad
    \mu_+ = \frac{10.25}{1+c_+}, \quad 
    \mu_- = \frac{10.25}{1+c_-}, \quad
    \gamma = 10, \\
    \lambda_+ & = 0.05, \quad
    \lambda_- = 0.05, \quad
    G_+ = \frac{c_+ \mu_+}{\lambda_+}, \quad
    G_- = \frac{c_- \mu_-}{\lambda_-}, \quad
    \alpha = 10^{-2}.
\end{align*}
Here, the constants $c_\pm \geq 0$ characterise the elastic contribution in the two phases $\Omega_\pm(t)$ of the viscoelastic fluid. Adjusting the values of $c_\pm$ allows us to change the viscosity ratios
\begin{align*}
     \frac{\mu_\pm}{\mu_\pm + G_\pm \lambda_\pm } 
    = \frac{1}{1+c_\pm} \in (0, 1],
\end{align*}
while maintaining the total shear viscosities $\mu_+ + G_+ \lambda_+ = \mu_- + G_- \lambda_- = 10.25$ in both phases.
We consider two different scenarios: In the first scenario, the outer phase is highly viscoelastic with $c_+ = 19.5$, while the inner phase is less viscoelastic with $c_- = 1$. In the second scenario, the roles are reversed: the inner phase is highly viscoelastic with $c_- = 19.5$, and the outer phase is less viscoelastic with $c_+ = 1$.

The time evolution of the rising bubble for both cases is shown in Figure \ref{fig:test94_Gvar_simulation}.
Figure \ref{fig:test94_Gvar1_steamlines_isovalues_t=1} illustrates the streamlines of the velocity $\bu^n$ at time $t=1$ and the isovalues of $\mathbb{B}_{11}^n$ and $\mathbb{B}_{22}^n$ for the first scenario ($c_+=19.5$, $c_-=1$). Similarly, Figure \ref{fig:test94_Gvar2_steamlines_isovalues_t=1} presents the corresponding results for the second scenario ($c_+=1$, $c_-=19.5$).
In the first scenario, where the outer phase is more viscoelastic, the bubble develops a small cusp at its bottom end due to the localized elastic stresses. In contrast, in the second scenario, where the inner phase is more viscoelastic, the bubble shape remains close to its initial circular shape.

The energy contributions for both cases, including kinetic energy, elastic energy, and interface length, are displayed in Figure \ref{fig:test94_Gvar_energy}. Additionally, Figure \ref{fig:test94_Gvar_energy_full} shows the total discrete energy and its changes (excluding the contribution from body forces) computed from the discrete energy inequality \eqref{eq:stability_FE_G} using
\begin{align}
    \label{eq:dissipation_Gvar} \nonumber
    - \mathcal{D}^{n+1} &\coloneqq 
    \frac12 \norm{\sqrt{\rho^n} \bu^{n+1}}_{L^2}^2 
    + \frac12 \skp{G^n}{\trace(\bbB^{n+1} - \ln\bbB^{n+1} - \bbI)}_{\calT^n}^h 
    + \gamma \calH^{d-1}(\Gamma^{n+1})
    \\ \nonumber
    &-\frac12 \norm{\sqrt{\mathrm{I}_0^n \rho^{n-1}} \bI_2^n \bu^n}_{L^2}^2
    - \frac12 \skp{G^n}{\trace(\bbB^{n} - \ln\bbB^{n} - \bbI)}_{\calT^n}^h 
    - \gamma \calH^{d-1}(\Gamma^n)
    \\
    &\quad
    - \Delta t \skp{\rho^n \mathbf{f}_1^{n+1} + \mathbf{f}_2^{n+1}}{\bu^{n+1}} .
\end{align}
The term $\mathcal{D}^{n+1}$ remains nonnegative throughout, confirming the dissipation of energy in the system. Similarly to above, the full energy is not monotonically decreasing due to the presence of a nonzero body acceleration.
}

\section{Conclusion}
In this work, we considered a parametric finite element approximation for viscoelastic two-phase flow, combining ideas from the two-phase Navier--Stokes model without viscoelasticity \cite{BGN_2015_navierstokes} and the one-phase viscoelastic Oldroyd-B model \cite{barrett_boyaval_2009}.
First, we presented a two-phase Navier--Stokes model that incorporates viscoelastic effects, and we derived jump conditions across the interface $\Gamma(t)$ and boundary conditions on $\partial\Omega$ using an energy identity. This approach provides the basis for a variational formulation, which is essential for our numerical approximation.
On the formal level, we emphasised physical properties such as energy conservation and volume preservation, which is the foundation for a stable finite element approximation. In our approach, we obtain a positive definite tensor $\bbB$, which is important not only for physical reasons but also to ensure finite energy. Moreover, the energy exhibits a dissipative behaviour in cases without external forces.
We used a regularised system to prove the unconditional stability and existence of solutions. Furthermore, we introduced a semi-discrete scheme that maintained a good mesh quality of the parametric mesh, incorporating an extended finite element method (XFEM) to guarantee exact conservation of phase volumes, which is a crucial aspect validated through numerical experiments.
In scenarios involving a variable shear modulus, adjustments to our numerical method were necessary. We outlined the reasons for these adjustments and presented an approach that can be used instead.

To validate the practicability of our numerical approach, we presented simulations in two dimensions, focussing on a scenario of a rising bubble surrounded by a viscoelastic fluid. By varying the intensity of viscoelasticity (ranging from none to low and high levels), we replicated typical behaviours observed in previous experimental studies (e.g.~\cite{hassager_1979_viscoelastic}). Our simulations demonstrated very good volume and parametric mesh properties, confirming the robustness of our numerical approach.


\section*{Acknowledgments}
H.G. and D.T. gratefully acknowledge the support by the Graduiertenkolleg 2339 IntComSin of the Deutsche Forschungsgemeinschaft (DFG, German Research Foundation) -- Project-ID 321821685.


\footnotesize\setlength{\parskip}{0cm}

\addcontentsline{toc}{section}{References}
\bibliography{main}

\end{document}